\documentclass[11pt,english]{article}
\usepackage{mathptmx}

\usepackage{a4wide}
\usepackage{amssymb,amsthm}
\usepackage{amsmath}
\usepackage{amsfonts}
\usepackage{latexsym}
\usepackage{amsxtra}
\usepackage{graphicx}
\usepackage{wrapfig}
\usepackage{mathrsfs}
\usepackage{epsfig}
\usepackage{times}
\usepackage{makeidx}
\usepackage[titles]{tocloft}

\newcommand{\heikodetail}[1]{}

\def\yy{\mbox{$\spadesuit$}}
\newcommand{\invisible}[1]{\par ($\spadesuit$ \emph{hidden comments in \TeX{}
file\/})} 
\newcommand{\reallyinvisible}[1]{}      

\newfont{\thickmath}{msbm10 scaled \magstephalf}%
\newfont{\smallthickmath}{msbm7 scaled \magstephalf}%
\newfont{\footnotethickmath}{msbm8}%
\newfont{\footnotesmallthickmath}{msbm6}%

\newtheorem{lemma}{\bf Lemma}[section]
\newtheorem{theorem}[lemma]{\bf Theorem}
\newtheorem{proposition}[lemma]{\bf Proposition}
\newtheorem{corollary}[lemma]{\bf Corollary}
\newtheorem{conjecture}[lemma]{\bf Conjecture}
\newtheorem{definition}[lemma]{\bf Definition}
\newtheorem{REMARK}[lemma]{\bf Remark}
\newtheorem{example}[lemma]{\bf Example}
\newcommand{\Fo}{\,\,\,\text{for }\,\,}
\newcommand{\Foa}{\,\,\,\text{for all }\,\,}

\newcommand{\INT}{\,\text{int}\,}

\newcommand\Reals{{\mathbb R}}
\newcommand\R{{\mathbb R}}
\newcommand\Z{{\mathbb Z}}

\newcommand\N{{\mathbb N}}
\renewcommand\S{{\mathbb S}}
\newcommand\Sphere{{\mathbb S}^{n-m-1}}

\newcommand{\bbbr}{\Reals}

\newcommand\dist{\mathop{\rm dist}\nolimits}
\newcommand\diam{\mathop{\rm diam}\nolimits}
\newcommand\subsp{\mathop{\rm span}\nolimits}
\newcommand\disk{D^{n-m}}
\newcommand\lk{\mathrm{lk}_2\,}
\renewcommand\deg{{\rm deg}\,}

\newcommand\ang{\mathop{\mbox{$<\!\!\!)$}}\nolimits}

\newcommand{\xx}{\mbox{$\clubsuit$}}

\newcommand{\E}{\mathcal{E}}

\newcommand{\A}{\mathcal{A}}

\def\mbbbr{\mbox{\smallthickmath R}}

\renewcommand{\R}{\bbbr}
\newcommand{\mR}{\mbbbr}

\newcommand{\eps}{\varepsilon}
\renewcommand{\H}{\mathscr{H}}

\def\osc{\mathop{\rm osc\,}}         
\def\rtp{R_{\rm tp}}

\begin{document}

\renewcommand{\thefigure}{\arabic{figure}}

\title{\large\bf Tangent-point repulsive potentials \\
for a class of non-smooth $m$-dimensional sets in $\bbbr^n$.\\
Part I: Smoothing and self-avoidance effects}

\author{\normalsize Pawe\l{} Strzelecki, Heiko von der Mosel}

\date{\normalsize version of \today}

\maketitle


\frenchspacing

\begin{abstract}
We consider repulsive potential energies $\E_q(\Sigma)$, whose integrand
measures tangent-point interactions,  on a large class of non-smooth
$m$-dimensional sets $\Sigma$ in $\R^n.$ Finiteness of the energy $\E_q(\Sigma)$
has three sorts of effects for the set $\Sigma$: topological effects
excluding all kinds of (a priori admissible) self-intersections, 
geometric and measure-theoretic effects, providing large projections of
$\Sigma$ onto suitable $m$-planes and therefore large $m$-dimensional
Hausdorff measure of $\Sigma$ within small balls up to a uniformly
controlled scale, and finally, regularizing effects culminating
in a geometric variant of the Morrey-Sobolev embedding theorem:
Any admissible set $\Sigma$ with finite $\E_q$-energy, for any exponent
$q>2m$, is, in fact, a $C^1$-manifold whose tangent planes vary
in a H\"older continuous manner with the optimal H\"older exponent
$\mu=1-(2m)/q$. Moreover, the patch size of the local $C^{1,\mu}$-graph 
representations is uniformly controlled from below
only in terms of the energy value
$\E_q(\Sigma)$. 

\vspace{2mm}


\end{abstract}

\setcounter{tocdepth}{1}
\tableofcontents
      
\bigskip 


\renewcommand\theequation{{\thesection{}.\arabic{equation}}}
\def\setnumbers{\setcounter{equation}{0}}


\section{Introduction}\label{sec:1}

This paper grew out of a larger project, devoted to the investigation of 
so-called \emph{geometric curvature energies} which include various types of geometric integrals, measuring the degree of smoothness and bending for objects that do not, at least a priori, have to be smooth. Here, we study the energy functional 
\begin{equation}
	\label{energy-intro}
\E_q(\Sigma) = \int_{\Sigma}\int_{\Sigma} \frac{1}{\rtp^q(x,y)} \,  d\H^m(x)\, d\H^m(y)
\end{equation}
defined for a class $\A$ of admissible, $m$-dimensional sets 
in $\R^n$. The  precise definition of $\A$ is given in Section \ref{sec:2}; 
we just mention now that for each $\Sigma\in \A$ a weak counterpart of the classic tangent plane  is defined 
almost everywhere with respect to the $m$-dimensional
Hausdorff measure $\H^m$ on $\Sigma$.  In other words,
for $\H^m$-a.e. $x\in\Sigma$ there is an $m$-plane $H_x$ such that
the portion of $\Sigma$ near the point $x$ is close to the
affine plane $x+H_x\subset\R^n$.
The quantity 
\begin{equation}\label{rtp}
\rtp (x,y):=\frac{|y-x|^2}{2\dist(y,x+H_x)}
\end{equation} 
in the integrand is referred to as the \emph{tangent-point radius} and denotes the radius of the smallest%
\heikodetail{

\bigskip

\xx are there more (larger) tangent spheres in higher codimension??\xx 

\yy As far as I can see, for $n-m>1$ the condition that an $(n-1)$-sphere be tangent to $x+H_x$ and pass through $y$ does not specify that sphere uniquely. Think of $2$-spheres in $\R^3$ tangent to the $x$-axis at $0$, say, and passing through $(0,1,0)$. There are many of them and they have various radii. \yy

\bigskip

}            
sphere tangent to the affine plane $x+H_x$  and passing through $y$. (If $y$ happens to be contained
in $x+H_x$, then we set $1/\rtp (x,y)=0$.) 
Thus, $1/\rtp(x,y)$ is defined a.e. on $\Sigma\times \Sigma$ with respect to the product measure $\H^m\otimes\H^m$. 
Notice that for any compact embedded manifold of class $C^{1,1}$
this repulsive potential $\E_q$ is finite. For two-dimensional surfaces in $\R^3$, i.e. $n=3$, $m=2$,
Banavar et al. \cite{banavar} suggested, in fact, the use of such tangent-point
functions to construct self-interaction energies with non-singular
integrands that do not require any sort of ad hoc regularization, in
contrast to standard repulsive potentials. The latter would 
penalize any two
surface points that are close in Euclidean distance, no matter
whether these points are adjacent on the surface (leading to singularities)
or belong to different
sheets of the same surface.
Our aim here is to show that for the infinite range of exponents $q>2m$
finiteness of $\E_q(\Sigma)$ has three sorts of consequences for any admissible set $\Sigma\in\A$: measure-theoretic, topological, and analytical. To see them in a proper perspective, let us give a plain description of the surfaces we work with.  

Our class $\A$ consists of $m$-dimensional sets  $\Sigma\in\R^n$ with finite measure $\H^m(\Sigma)<\infty$  on which  we impose (1) a certain degree of flatness in the neighbourhood of many (but a priori not all!) points of $\Sigma$, and  (2) some degree of connectivity.
A priori, we allow for various self--intersections of $\Sigma$, and for singularities along low dimensional subsets.  For the purposes of this introduction,
however, it is enough to keep in mind the following examples of admissible surfaces (more general examples are presented in Section   \ref{sec:2.3}):    
\medskip
\begin{enumerate}
\renewcommand{\labelenumi}{(\roman{enumi})}
	\item If $\Sigma_0=M_1\cup \ldots \cup M_N$, where $N\in\N$ is arbitrary and all $M_i\subset\R^n$ are compact, closed, embedded $m$-dimensional submanifolds of class $C^1$ such that $\H^m(M_i\cap M_j)=0$ whenever $i\not=j$, then $\Sigma_0$ is admissible; 
	\item If $\Sigma_0$ is as above, then $\Sigma_1=F(\Sigma_0)$ is admissible whenever $F$ is a bilipschitz homeomorphism of $\R^n$. 
\end{enumerate}
\medskip\noindent	
The dimension $m$ and the codimension $n-m$ of $\Sigma$ in $\R^n$ are fixed throughout the paper but otherwise arbitrary. The reader may adopt for
now the temporary definition
\[
\A\colon=\{\Sigma\subset\R^n\ \colon \ \Sigma=F(\Sigma_0),\ \Sigma_0\mbox{ as in (i) above, } F\colon \R^n\to\R^n \mbox{ bilipschitz}\}.
\]

It is easy to see that $q_0=2m=\dim (\Sigma\times\Sigma)$ is a critical exponent here: for $q=q_0$ the energy $\E_q(\Sigma)$ is scale invariant, and for each $q\ge q_0$ a surface $\Sigma$ with a conical singularity at one point must have $\E_q=\infty$. We prove in this paper that for $q>q_0=2m$ all kinds of singularities are excluded. In fact, upper bounds for $\E_q(\Sigma)$ lead to three kinds of
effects.  Firstly,
measure-theoretic effects:
the measure of $\Sigma$ contained in a ball or radius $r$ is
comparable to $r^m$ on small scales that  depend
solely on the energy. Secondly topological effects: 
an admissible surface $\Sigma$ with finite $\E_q$-energy has no self-intersections, it must be an embedded manifold, and finally,
far-reaching analytical consequences: we have precise $C^{1,\mu}$ bounds for the charts in an atlas of $\Sigma$.
         
%
%
%
   
\begin{figure}[!t] 
\begin{center}
\includegraphics*[totalheight=9cm]{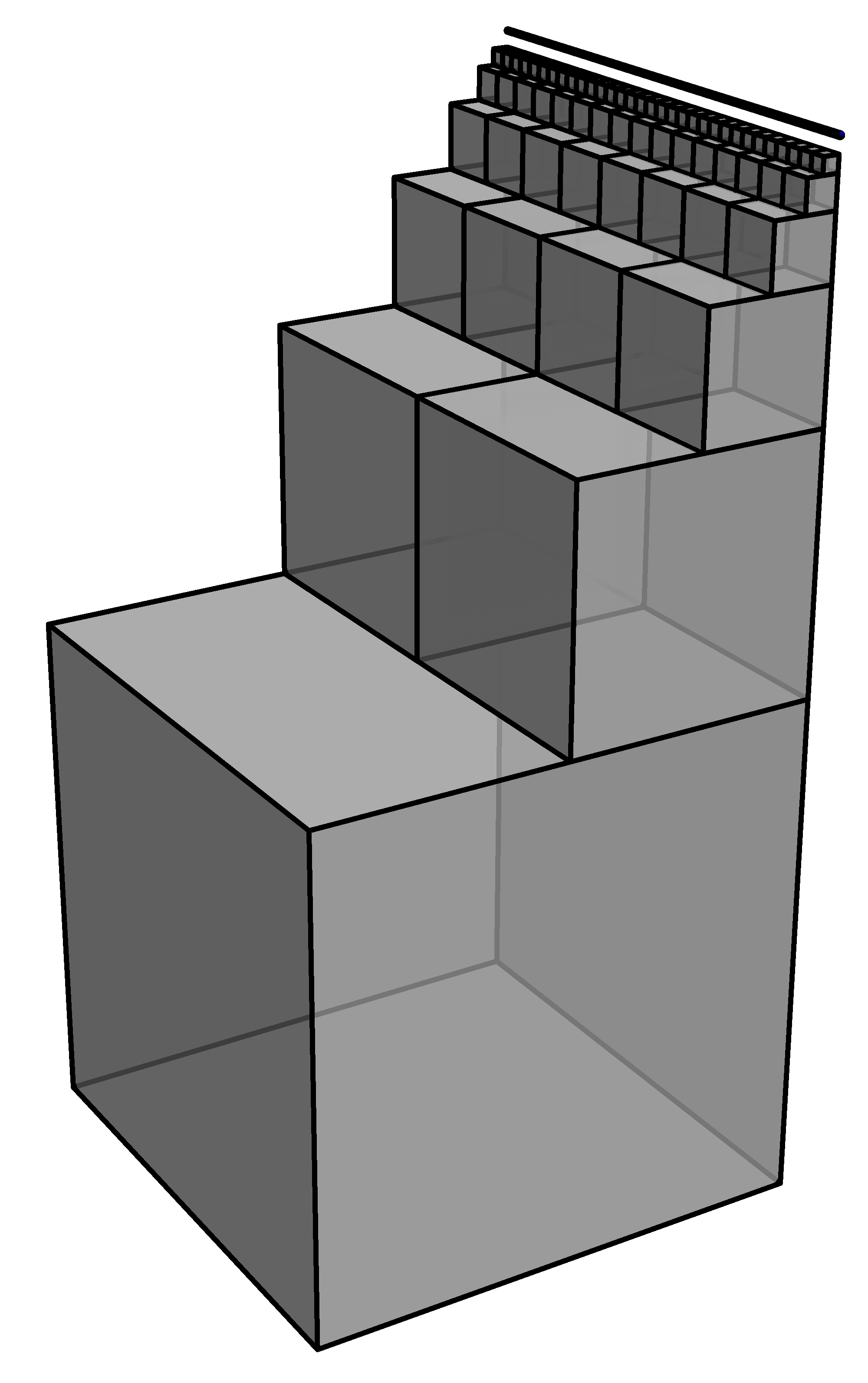}
                   
\medskip
\end{center}      

\caption{An infinite stack of cubes, also an admissible surface
(which turns out to have infinite $\E_q$-energy for all $q>2m=4$).
}
\label{fig:1} 

\end{figure}

Let us first state the results precisely and then comment on the proofs and discuss the relations of this paper to existing research.      

\begin{REMARK} Keep in mind, though, that all results stated in the introduction will be proved for a more general class $\A(\delta)$ of admissible sets much
larger than the preliminary class $\A$ defined above; 
see Section \ref{sec:2.3}. To get a first impression of  other
admissible sets have a look at Figure \ref{fig:1}.
\end{REMARK}          
                        
\medskip

\begin{theorem}[\textbf{Uniform Ahlfors regularity}]\label{thm:UAR}
Assume that $\Sigma\in\A$ is an admissible $m$-dimen\-sio\-nal surface in $\R^n$ with $\E_q(\Sigma)\le E$, $q>2m$. There exists a constant $a_1=a_1(q,n,m)>0$, depending \emph{only} on $q,n$ and $m$,  such that
\[
\H^m(\Sigma \cap B(x,r))\ge \frac 12 \omega(m) r^m
\]
for all $x\in \Sigma$ and all radii 
\[
0< r< R_1\equiv R_1(q,n,m,E):=\frac{a_1}{E^{1/(q-2m)}}\, .
\]
\label{conj1}
(Here, and throughout the paper,  $B(x,d)$ denotes 
the closed ball of radius $d$ centered
at $x$.)
\end{theorem}

In other words: if $\Sigma\in \A$ has finite energy for some $q>2m$, then up to the length scale given by $R_1$ -- which depends \emph{only\/} on the energy bound $E$ and the parameters $m,n,q$, but \emph{not} on
$\Sigma$ itself -- isolated \emph{thin fingers\/}, narrow tubes, and the like\
cannot form on $\Sigma$. The
measure of the portion of  $\Sigma$ inside the ball $B(x,r)$ is at least as large as half of the measure of the $m$-dimensional equatorial cross-section of $B(x,r)$.    
A similar lower estimate on the 
Ahlfors regularity was proven by L. Simon for smooth
two-dimensional surfaces with finite Willmore
energy \cite[Corollary 1.3]{Si93}; see also the work of P. Topping
\cite{topping} which even contains {\it sharp} lower bounds for the sum of local $L^2$-norm of the classic mean curvature and the area of the surface in a small ball.  Mean curvature at a particular point $x$ on a
smooth surface in $\R^3$ may be viewed as the arithmetic mean
of minimal and maximal normal curvature at $x$. $\rtp^{-1}$, on the other hand,
is a two-point function taking non-local interactions into account as well,
but if one looks at the coalescent limits $\lim_{y\to x}\rtp^{-1}(x,y)$
one obtains absolute values of intermediate normal curvatures at $x$
depending
on the direction of approach as $y$ tends to $x$ (cf. \cite[Section 3.2]{banavar}). So, the local portion 
of our energy $\E_q$ near $x$ may be regarded as another kind of 
averaging normal curvatures at $x$, leading to density estimates
as does the Willmore functional.

\heikodetail{

\bigskip

{\tt\raggedright \xx Remark: absolute values of normal curvatures is what Rumpf 
had talked about when I met him and discussed global curvature
6-7 years ago!!

\yy I do realize. *Very good* that you mention Banavar's et al.  here

}}

The next result gives a quantitative description of flatness of $\Sigma$, in terms of the so-called $\beta$-numbers introduced by P. Jones.

\begin{theorem}[\textbf{Uniform decay of $\beta$-numbers}] 
\label{thm:betas} 
Let $\Sigma\in\A$ be an admissible $m$-di\-men\-sional surface in $\R^n$ with $\E_q(\Sigma)<E$ for some $q>2m$. 
There exist two constants $a_2(q,n,m)>0$ and
$A_2(q,n,m)\linebreak[0]<\infty$, both depending \emph{only} on $n,m$ and $q$, such that whenever the radius
\begin{equation}
d\le R_2\equiv
R_2(q,n,m,E):=\frac{a_2(q,n,m)}{E^{1/(q-2m)}}
\end{equation}
and the bound $\eps>0$ satisfy the balance condition
\begin{equation}  
	\label{eq:1.3}
\eps^{4m+q}d^{2m-q} \ge A_2(q,n,m) E\, ,
\end{equation}
then we have
\begin{equation}\label{betasmall}
\beta_\Sigma(x,d) :=\inf_{P\in G(n,m)}\left(\sup_{y\in B(x,d)\cap \Sigma} \frac{\dist(y,x+P)}{d}\right)\le \eps\, , \qquad x\in \Sigma\, ,
\end{equation}
where $G(n,m)$ denotes the Grassmannian of all $m$-dimensional subspaces
of $\R^n$.

\end{theorem}
Thus, for small $d$ we have
\[
\beta_\Sigma(x,d) \lesssim E^{1/(4m+q)} d^\kappa\, ,\qquad
\kappa:= \frac{q-2m}{q+4m}>0\, .
\]
It is known that this condition alone does not suffice to conclude that $\Sigma$ is a topological manifold. D. Preiss, X. Tolsa and T. Toro \cite{ptt}, extending an earlier work of G. David, C. Kenig and T. Toro \cite{davidkenigtoro}, study Reifenberg flat sets $\Sigma$ whose $\beta$-numbers satisfy such estimates, see e.g. \cite[Prop. 2.4]{ptt} where it is proved that a decay bound for $\beta$'s combined with Reifenberg flatness\footnote{We do not define that condition here 
since we will not work with it directly; let us just mention  that Reifenberg flatness means that the rescaled \emph{Hausdorff distance\/} between $\Sigma \cap B(x,d)$ and an $m$-plane $P$ in $B(x,d)$ is uniformly controlled, and small.} implies that $\Sigma$ must be a submanifold of class $C^{1,\kappa}$. 

Since $\Sigma\in \A$ might, at least a priori, have transversal self--intersections, we do not have Reifenberg flatness here, and a quick direct use of the results of \cite{davidkenigtoro,ptt} is impossible. However, we are able to use the energy estimates and the information given by Theorem~\ref{thm:betas} iteratively. Extending the ideas from our earlier work
\cite[Section~5]{svdm-surfaces} devoted to surfaces in $\R^3$, we prove  here
that at every point $x\in\Sigma$ there exists the classic tangent plane $T_x\Sigma$,  and that the oscillation of tangent planes along $\Sigma$
satisfies uniform
H\"older estimates. 
This implies that each $\Sigma\in\A$ with $\E_q(\Sigma)<\infty$ must be  an embedded $m$-dimensional manifold of class $C^{1,\kappa}$. Later on, working with graph patches of $\Sigma$, we  use slicing techniques and a bootstrap reasoning to improve and sharpen this information. The following theorem is the main result of this paper.

\begin{theorem}[\textbf{Geometric Sobolev--Morrey imbedding}]                    
\label{thm:bootstrap}
Let $\Sigma\in\A$ and $\E_q(\Sigma)<+\infty$ for some $q>2m$. Then $\Sigma$ is an embedded submanifold of class $C^{1,\mu}$, where $\mu = 1 - 2m/q$. 

In fact, there exist constants $a_3, A_3>0$, depending only on $m,n,$ and $q$, with the following property:  For each $x\in \Sigma$ and each $r\le R_3=a_3\E_q(\Sigma)^{1/(q-2m)}$ there exists an $m$-plane $P\in G(n,m)$ and a  function $f\colon P\simeq  \R^m\to P^\perp\simeq\R^{n-m}$ of class $C^{1,\mu}$ such that
\[
\Sigma\cap B(x,r) = \Sigma \cap \mathrm{graph}\,  f\, ,
\]     
where $\mathrm{graph}\,  f \subset P\times P^\perp =\R^n$ denotes the graph of $f$, and 
\begin{equation}
|\nabla f(z)-\nabla f(w)|\le A_3  E(x,r)^{1/q}  |z-w|^\mu,  \qquad z,w\in P \cap B(0,r), 
\label{optimal}  
\end{equation}
where 
\[ 
E(x,r):= \int_{B(x,r)\cap \Sigma}\int_{B(x,r)\cap \Sigma} \left(\frac{1}{\rtp (u,v)
}\right)^q \,d\H^m(u)\, d\H^m(v)\, .
\] \label{conj:2.3}      
\end{theorem}

We believe that the exponent $\mu=1-{2m}/{q}$, strictly larger than $\kappa=(q-2m)/(q+4m)$, is optimal here. It is clear that finiteness of $\E_q$ does not lead to $C^2$ regularity: consider  a rotational cylinder  closed with two hemispherical caps as an  admissible surface $\Sigma$ of class 
$C^{1,1}$ but not in $C^2$. For this particular surface inequality \eqref{optimal} is qualitatively optimal and, due to the factor $E(x,r)^{1/q}$ and boundedness of $1/\rtp$, yields in fact \emph{Lipschitz} estimates for the gradient of local graph representations of $\Sigma$.

Please note two more things. First, the exponent $\mu=1-2m/q$ is computed according to the recipe used in the classic Sobolev--Morrey imbedding theorem in the supercritical case. Here, the dimension of the domain of integration, i.e. of   $\Sigma\times\Sigma$, equals $2m$. We have $\mu\to 1$ as $q\to \infty$; 
for two-dimensional surfaces, the limiting case $q=\infty$ has been treated earlier in our papers \cite{StvdM1} and~\cite{StvdM2}. 

Second, what we have learnt about $\Sigma$ is not limited to 
embeddedness and purely qualitative $C^{1,\mu}$ estimates. 
It is clear that the bounds given 
by Theorem~\ref{thm:bootstrap} are uniform in any class of surfaces with uniformly bounded energy $\E_q$. In other words, if  $\mathcal{K}=\{\Sigma_i\colon i\in I\}\subset \A$ satisfies 
\begin{equation}
\sup_{i\in I} \E_q(\Sigma_i) \le M < \infty\, ,      
\label{energybound}
\end{equation}
then we can find two constants $A,\delta > 0$, depending only on $M,m,n$ and $q$, such that each $\Sigma_i\cap B(x,\delta)$, where $i\in I$ and $x\in \Sigma_i$, is obtained by a rigid motion of $\R^n$ from a graph of a function $f\colon \R^m\to (\R^m)^\perp \simeq \R^{n-m}$ which satisfies the uniform estimate $\| f\|_{C^{1,\mu}}\le A$, no matter how $i\in I$ and $x\in \Sigma_i$ have been chosen. Thus, a uniform upper bound on $\E_q$ allows us to fix a uniform size of charts for all  $\Sigma_i\in \mathcal K$, and forces  the equicontinuity of gradients of local graph representations of the surfaces $\Sigma_i$. In a forthcoming paper \cite{svdm-tp2} we show how to use this idea to obtain finiteness theorems for classes of $C^{1}$ embedded manifolds $\Sigma_i$ in $\R^n$ satisfying a volume constraint and a uniform energy bound \eqref{energybound}. 

We do not know what happens in the critical case $q=2m$. Let us mention here one plausible conjecture that we cannot prove at this stage.

\begin{conjecture}
Every immersed $m$-dimensional $C^1$-manifold in $\Sigma\subset\R^n$ with finite $\E_{2m}$-energy is embedded.
\end{conjecture} 

Another, probably more difficult, question that we cannot handle at present
is the following: how regular are the minimizers of $\E_q$ (say, with 
upper bounds for the total measure, to prevent the decrease of energy 
caused by rescaling) in isotopy classes of $C^1$ embedded manifolds? Are they $C^{1,1}$ 
(this is optimal for \emph{ideal links\/} \cite{CKS} -- corresponding to the case $q=\infty$ in dimensions
$n=3$ and $m=1$ -- where contact phenomena are present)? Or maybe $C^\infty$, as minimizers of a M\"obius invariant knot energy in \cite{FHW}, \cite{He};
see also \cite{reiter-phd}, \cite{reiter-th1}?  In addition, S. Blatt
\cite{blatt-preprint} characterized all curves with finite M\"obius energy as embeddings in certain
Sobolev-Slobodecki\u \i{} classes; such a {\it characterization} of finite energy submanifolds for the tangent-point energy $\E_q$ is presently not known.

\bigskip

Our interest in this topic has been triggered by several factors. They include manifold applications of Menger curvature in harmonic analysis and geometric measure theory (see e.g. the survey articles of P. Mattila \cite{Ma98}, \cite{Ma04}, G. David \cite{david-survey} and X. Tolsa \cite{tolsa-survey}, and the literature cited therein, including J.C. L\'{e}ger \cite{leger} and the relation between $1$-rectifiability and $L^2$-integrability of Menger curvature). There are also works of different origin, investigating another geometric concept, the so-called global curvature introduced by Gonzalez and Maddocks \cite{GM}. The second author of the present paper took part in laying out the strict mathematical foundations for global curvature of rectifiable loops and its variational applications to elastic curves and rods with positive thickness;
see \cite{GMSvdM}, \cite{heiko1}, \cite{heiko2}, \cite{heiko3}, 
\cite{gerlachvdm1}, \cite{gerlachvdm2}. 
Part of this work, in turn, has been a starting point for our subsequent joint research devoted to  various energies that, roughly speaking, interpolate between global curvature and Menger curvature. Finiteness of these energies; see e.g. \cite{ssvdm-triple}, \cite{svdm-surfaces}, \cite{svdm-tpcurves}, analogously to the case that we consider here,  leads to an increase of regularity, to compactness effects, and yields a tool to control the amount of bending of non-smooth objects in purely geometric terms.\footnote{In \cite{svdm-tpcurves} we treat the toy case $m=1$ of the present paper, along with a few  knot-theoretic applications of $\E_q$ for curves.} The novelty
in the present paper is that we work in full generality, overcoming the difficulty that both the dimension and the codimension may be arbitrary.  In an ongoing research \cite{skol1,skol2} S. Kolasi\'{n}ski obtains analogues of our results for basically the same admissible class of surfaces that we consider here, but for a different integral energy, defined as an $(m+2)$-fold integral (with respect to $\H^m$) over the set of all simplices with vertices on $\Sigma$,  directly extending our 
results in \cite{svdm-surfaces} to surfaces of arbitrary dimension and codimension.

\heikodetail{

\bigskip

{\tt just a remark similar to what I found in Banavar et al.: an $m$-plane 
is spanned by $m+1$ points, in that sense also $\rtp $ is a function
depending on $m+2$ points like Menger curvature for $m$-dim. surfaces...\xx\xx}  
}

Closely related research includes also G. Lerman and J.T. Whitehouse \cite{LW08a}, \cite{LW08b}, who investigate a number of ingenious high-dimensional curvatures of Menger type and obtain rectifiability criteria for $d$-dimensional subsets of Hilbert spaces. Last but not least, the deep and classic paper of W. Allard \cite{allard}   sets forth a regularity theory for  $m $-dimensional  varifolds whose first variation (roughly: the distributional counterpart of mean curvature) is in $L^p$ for some $p>m $. Our regularity results bear some resemblance to his Theorem 8.1. There are many differences, though, that remain to be fully understood. It is clear that without some extra topological assumptions on $\Sigma$ finiteness of \eqref{energy-intro} cannot lead to the
conclusion that $\Sigma$ is locally (on a scale depending only on the energy!) homeomorphic to a disc; one could punch an arbitrary number of `holes' in a smooth surface and this would just decrease the energy we work with. In Allard's case, once we fix a ball where appropriate density estimates hold and the weight $\|V\|$ of the varifold  $V$ is close to the Hausdorff measure of a disk, then the `lack of holes' is built into his assumption on the first variation $\delta V$ of $V$. On the other hand, $\E_q$ -- as a non-local energy in 
contrast to the locally defined distributional mean curvature -- 
averages over all global tangent-point interactions, which leads
to self-avoidance and control over topology of the given surface. Admissible sets with finite $\E_q$-energy are differentiable manifolds, which Allard's result
cannot guarantee for varifolds with distributional mean curvature in $L^p$, $p>m$: there is a remaining (small) singular set, such that there is no control
on the topology of the support of the varifold measure. To possibly
bridge the
apparent gap between Allard's work and our results
we should note that versions of $\E_q$ can be defined for general $m$-dimensional varifolds $V$, via double integrals:  the integrand $1/\rtp$ can be treated as a function on \emph{points and planes}. It is an intriguing question whether finiteness of such integrals for some $q$'s lead to rectifiability criteria or to an improved regularity in the case of varifolds.

\heikodetail{

\bigskip

{\tt\xx\raggedright
Somehow I refrained from talking about the disadvantage of
non-local energies regarding lack of pde, lack of geometric flows,
should we mention that at all? 

\yy\yy No, I think.

}}

\bigskip

Let us now informally sketch the main thread of our reasoning, and describe the organization of the paper in more detail. We want to exclude self intersections and to have a quantitative description of flatness; for this, Theorem~\ref{thm:betas} would be a good starting point. The main idea behind its proof is pretty straightforward: if the $\beta$-numbers were too large, i.e. if $x\in \Sigma$ but $\Sigma\cap B(x,d)$ were not confined to a narrow tube $B_{\eps d}(x+P)$ around some affine $m$-plane $x+P$, then,  a simple argument shows that we would have two much smaller balls $B_1,B_2\subset B$, say with 
\[
\diam B_1 = \diam B_2 \approx \eps^2 d  
\]
such that for all $y\in \Sigma \cap B_2$ and a nonzero proportion of $z\in \Sigma \cap B_1$ the distance $\dist (y,z+T_z\Sigma)$ would be comparable to $\eps d$. This yields
$1/\rtp (y,z)\gtrsim \eps/d$, and a lower bound for the energy follows easily, leading to a contradiction, \emph{if} the $B_1,B_2$ and the bound for the $\beta$'s are chosen in a suitable way which happens to be precisely the balance condition \eqref{eq:1.3}. There is only one serious catch here: in order to make the resulting estimate uniform, and to be able to iterate it later on, we must guarantee that 
\[
\H^m(\Sigma\cap B_r)  \ge c\cdot r^m   \qquad\mbox{for all $r<r_0=r_0(\text{energy})$,}
\]
with some \emph{absolute} constant $c$. And we want \emph{both} $r_0$ and $c$ independent of a particular $\Sigma$.

For this, we need Theorem~\ref{thm:UAR} which serves as the backbone for all the later constructions and estimates of the paper. The overall idea here is somewhat similar to an analogous result in our work \cite{svdm-surfaces} on Menger curvature for surfaces in $\R^3$. The main difference, however, leading to crucial difficulties, is that the \emph{codimension\/} of $\Sigma$ may be arbitrary.

The proof of Theorem~\ref{thm:UAR} has two stages. First, for a fixed generic point $x\in \Sigma$ and all radii $r$  
below a stopping distance $d_s(x)$, 
we control the size of projections of $\Sigma\cap B(x,r)$ onto some $m$-plane $H(r) $ (which may vary as $r$ varies). Here, topology 
comes into play. To grasp the essence of our idea, it is convenient to think of $\Sigma=M_1\cup\ldots \cup M_N$ as in  Example (i) at the beginning of the introduction. For $x\in M_i\setminus \bigcup_{j\not=i} M_j$ and for infinitesimally small radii $r$ we start with the  tangent planes $P=T_xM_i$, 
and note that   small $(n-m-1)$-spheres 
that are perpendicular to $T_xM_i$ are nontrivially linked with $M_i$. 
Then, for a sequence of growing radii $\rho$, we rotate $P$ 
if necessary
by a controlled angle to a new position $P_\rho$  in order to keep the projections large. At the same time, we construct a growing connected \emph{excluded region\/} $S_\rho$ which 
does not contain any point of $\Sigma$ in its interior.
The size of the projections is controlled via a topological argument, involving the homotopy invariance of the linking number mod 2 of submanifolds. 
%
The  construction stops at  some stopping distance
$r=d_s(x)$ , and yields another point $y\in \Sigma$ with $|y-x|\approx d_s(x)$ and two smaller balls $B(x,cd_s(x) )$, $B(y,cd_s(x) )$, where $c\in (0,1)$ is an explicit absolute constant, such that 
\[
\frac{1}{\rtp(z,w)}  \gtrsim \frac{1}{d_s(x)}
\]                                  
for all $w\in B(y,cd_s(x))$ and a significant proportion 
of $z\in B(x,cd_s(x) )$. In the second 
stage we use the energy bounds to show that 
$d(\Sigma):=\inf_{x\in \Sigma} d_s(x)$ is positive and satisfies $d(\Sigma)
\ge R_1$, where $R_1$ is the uniform constant given in Theorem~\ref{thm:UAR}. 
The details of that part are given in Section \ref{sec:4}; see Lemma~\ref{1/R-est}, Lemma~\ref{mainlemma} and their corollaries.

Sections \ref{sec:2} and 
\ref{sec:3} contain all the necessary prerequisites and are included for the sake of completeness. In Section~\ref{sec:2} we gather elementary estimates of angles between planes spanned by nearby almost orthogonal bases, and introduce the class 
of admissible sets whose definition is designed so that the above sketchy idea can be made precise. In Section~\ref{sec:3} 
we explain how the linking number mod 2 can be used for elements of $\A$, providing specific statements (and short proofs) for sake of further reference.  

Once Theorem~\ref{thm:UAR} is proved, we use the Hausdorff convergence of excluded regions defined for generic points $x\in \Sigma$ to obtain a corollary which, roughly speaking, ascertains that for \emph{every\/} $x\in\Sigma$ and $r\le R_1$ there is some plane $H=H_{x,r}\in G(n,m)$ such that $\Sigma\cap B(x,r)$ has large projection onto $H$ and is contained either in $B(x,r/2)$ (where, a priori, at this stage of the reasoning, $\Sigma$ might behave in a pretty wild way) or in a narrow tubular region $B_\delta(x+H)$, for some specific constant $\delta \ll 1$. A use of energy bounds yields now Theorem~\ref{thm:betas}, and an iterative argument implies that in fact $\Sigma$ must locally be a $C^{1,\kappa}$ graph. All this is done in Section~\ref{sec:5}. Embeddedness of $\Sigma$ is established here, too.   

Finally, in Section~\ref{sec:6}, we prove Theorem~\ref{thm:bootstrap} and sharpen the H\"older bounds. To this end, we show that if $\Sigma \cap B$ is a graph of $f\in C^{1,\kappa}$, then $\nabla f$ satisfies an improved estimate,
\begin{equation}       
\label{improved}
|\nabla f(a_1)-\nabla f(a_2)| \le 2 \Phi^\ast
(|a_1-a_2|/N) + C \, E^{1/q} \, |a_1-a_2|^\mu,
\end{equation}   
where $\Phi^\ast(s)$ stands for the supremum of oscillations of $\nabla f$ over all possible balls of radius $s$, and $E$ is the portion of energy coming from some ball containing $a_1,a_2$. The point is that \eqref{improved} holds for some $N=N(q)\gg 1$, so that for $f\in C^{1,\kappa}$ the first term of the right hand side can be viewed as an unimportant, small scale perturbation. The main idea behind \eqref{improved} is that when the integral average of $ (1/\rtp)^q$ is bounded by $K$, then there are numerous points $u_i$ in small balls around the $a_i$, $i=1,2$, where $(1/\rtp)^{q}\lesssim K$. A geometric argument implies that for such points $|\nabla f(u_1)-\nabla f(u_2)|$ can be controlled by the second term in the right hand side of \eqref{improved}, and a routine iterative reasoning, with a certain Morrey--Campanato flavour, allows us to get rid of the $2\Phi^\ast$ and finish the whole proof.


\medskip \noindent {\bf Acknowledgement.}\, The  authors  would like to
thank the Deutsche Forschungsgemeinschaft, Polish Ministry of
Science and Higher Education, and the Alexander von Humboldt
Foundation, for generously supporting this research. Substantial parts
of this work have been written while the first author has been staying
at the RWTH Aachen University in the fall of 2009; he is very grateful
to his German colleagues for their hospitality.


\section{Bases, projections,   angle estimates, and the class of admissible sets}\label{sec:2}   

\setnumbers

\medskip

\subsection{Balls, slabs, planes}
\label{sec:2.1}

We write $B(x,r)$ to denote the \emph{closed\/} ball in $\bbbr^n$,
with center $x$ and radius $r>0$. The volume of the unit ball in
$\bbbr^k$ is denoted by $\omega(k)$.

For a closed set $F$ in $\bbbr^n$ we set
\[
U_\delta(F) :=\{ x\in \bbbr^n\, \colon \, \dist(x,F)<\delta\},
\qquad \delta>0.
\]

$G(n,m)$ denotes the Grassmannian of all $m$-dimensional linear
subspaces of $\bbbr^n$. If $P\in G(n,m)$, then $\pi_P$ denotes the
orthogonal projection of $\bbbr^n$ onto $P$, and $Q_P$ is the
orthogonal projection onto $P^\perp \in G(n,n-m)$.

For two planes $P_1,P_ 2\in G(n,m)$ we define their distance (or
angle)
\[
\ang(P_1,P_2)\equiv d(P_1,P_2) := \|\pi_{P_1}-\pi_{P_2}\|\, ,
\]
where the right hand side is the usual norm of the linear map
$\pi_{P_1}-\pi_{P_2}\colon \bbbr^n\to \bbbr^n$. The Grassmannian
$G(n,m)$ equipped with this metric is compact.

Finally, we use the following variant of P. Jones' beta-numbers (see David's and Semmes' monograph \cite[Chapter 1, Sec. 1.3]{davidsemmes} for a discussion):
\begin{equation}
\label{Jones-b}
 \beta_\Sigma(x,r) := \inf_{L\in G(n,m)} \left( \sup_{y\in \Sigma \cap B(x,r)}\
 \frac{\dist (y,x+L)}{r}\right) , \qquad x\in \Sigma,\quad r>0.
\end{equation}

\subsection{Nearby planes: bases, projections, angle estimates}

\label{sec:linear}
\label{sec:2.2}

Throughout most of the paper, we shall work with estimates of
various geometric quantities related to two planes in $G(n,m)$
that form a small angle.  For sake of further reference, we
gather here several such estimates. We also fix specific constants
(which in all cases are far from being optimal) that are needed
later, in more involved computations in
Sections~\ref{ahlfors}--\ref{slicing}. All proofs are elementary,
but we provide them to make the exposition complete.

\begin{lemma}\label{sequences} Assume that $a,b>0$ and a
sequence of nonnegative numbers $s_k$ satisfies $s_1\le 1$,
\[
s_{k+1} \le a k + b\sum_{j=1}^k s_j, \qquad k\ge 1.
\]
Then for each $A\ge 1+\max (2a,2b)$ we have $s_k< A^k$,
$k=1,2\ldots$.
\end{lemma}

\smallskip\noindent\textbf{Proof.} One proceeds by induction.
Clearly, for $k=1$ we just need $s_1\le 1<A$. For each $A>1$ the
recursive condition for $s_{k+1}$ yields, under the inductive
hypothesis,
\begin{equation}
\label{skAk} s_{k+1}< a k + \frac{Ab}{A-1} (A^k-1)
\end{equation}
Now, $A\ge 1+\max(2a,2b)$ guarantees that $2ak < (1+2a)^k \le
A^k\le A^{k+1}$ and$\frac{b}{A-1}\le \frac 12$. Thus, \eqref{skAk}
yields $2s_{k+1}<  A^{k+1}+ A^{k+1}- A< 2A^{k+1}$. \hfill $\Box$

\begin{lemma}\label{bases-angle}
If $X,Y\in G(n,l)$ have orthonormal bases
$(e_j)\subset X$ and $(f_j)\subset Y$ such that $|e_j-f_j|\le
\alpha$ for each $j=1,\ldots, l$, then $\ang (X,Y)\le 2l\alpha$.
\end{lemma}

\smallskip\noindent\textbf{Proof.} Take an arbitrary
unit vector $v\in \R^n$ and estimate $|\pi_{X}(v)-\pi_{Y}(v)|$,
expressing both projections in orthonormal bases $(e_j)$ and
$(f_j)$. \hfill $\Box$

\begin{lemma}\label{GramSchmidt} Assume that $1\le l\le m\le n$.
If $e_1,\ldots, e_l$ is an orthonormal basis of
a subspace $X\in G(n,l)$ and $h_1,\ldots, h_l\in \R^n$ satisfy $|h_i-e_i|<
\eps < \eps_1:=10^{-1}(10^m+1)^{-1}$, then $(h_i)_{i=1,\ldots,l}$ are linearly independent. Moreover, the Gram-Schmidt orthogonalization process 
\[
u_i:=\frac{v_i}{|v_i|}, \qquad\mbox{where} \quad v_1=h_1, \quad v_{k+1}= h_{k+1}- \sum_{j=1}^k \frac{\langle
h_{k+1},v_j \rangle}{|v_j|^2} v_j\, , \quad  k+1\le l,
\]
yields vectors $v_i,u_i$ $(i=1,\ldots, l)$ that satisfy
\begin{gather}
|v_k-h_k| <  10^k \eps, \qquad \bigl||v_k|-1\bigr| <  (10^k+1) \eps 
<  \frac 1{10}\quad\mbox{for all 
$k=1,\ldots,l$,}\label{vkhk}\\
|u_k-e_k|< c_1\eps <\frac 12  \quad\mbox{for all $i=1,\ldots,
l$,}
\end{gather}
where $c_1:=2(10^m+1)$. If $Y=\mathrm{span}\, (h_1,\ldots, h_l)$, then 
\begin{equation}
\label{angXY}
\ang (X,Y) \le c_2\eps\, ,
\end{equation}
with $c_2:=2mc_1 = 4m(10^m+1)$.
\label{bases}
\end{lemma}

\reallyinvisible{
\begin{lemma} Assume that $l\le m$.
If $X,Y\in G(n,l)$, $e_1,\ldots, e_l$ is an orthonormal basis of
$X$ and $Y=\mathrm{span}\, (h_1,\ldots, h_l)$ with $|h_i-e_i|<
\eps < \eps_1:=10^{-1}(10^m+1)^{-1}$, then:
\begin{enumerate}
\item[{\rm (i)}] there is an orthonormal basis $(u_1,\ldots, u_l)$ of
$Y$ such that $|u_i-e_i|< c_1\eps <\frac 12$ for $i=1,\ldots,
l$;
\item[{\rm (ii)}] we have
\[
\ang (X,Y) \le c_2\eps\, ,
\]
\end{enumerate}
where $c_1:=2(10^m+1)$ and $c_2:=2mc_1 = 4m(10^m+1)$.
\end{lemma}
}

\medskip\noindent\textbf{Proof.}
As $|h_j-e_j|<\eps$ for all $j$, we have $|\langle h_i, h_j\rangle-\langle e_i, e_j\rangle| <  3\eps$. Therefore,
$|\langle h_{k+1},v_j\rangle| <  3\eps +(1+\eps)|h_j-v_j|$ for
$j=1,\ldots,k$ and $k\le l-1$. Using this observation, one proves \eqref{vkhk} by induction; assuming \eqref{vkhk} for $k$ and all $j<k$, we obtain
\begin{eqnarray*}
|v_{k+1}- h_{k+1}| & \le & \sum_{j=1}^k \frac{|\langle h_{k+1},v_j
\rangle|}{|v_j|} 
<   \sum_{j=1}^k \frac{3\eps +
(1+\eps)|v_j-h_j|}{|v_j|}\\
& <  & \frac{10}{9} 3\, k\eps + \frac{11}9 \sum_{j=1}^k10^j\eps < 10^{k+1}\eps,
\end{eqnarray*}
where the last inequality follows from elementary computations (the estimate is not sharp). This yields the first part of \eqref{vkhk} for $k+1$; the second one follows from the triangle inequality.

In particular, we also have $\dist \bigl(h_{k+1}, \subsp(h_1,\ldots, h_k)\bigr)=|v_{k+1}|>0$, and therefore $h_1,\ldots, h_l$ are linearly independent.

Setting $u_i:=v_i/|v_i|$, we easily conclude the proof of the whole lemma. (To check inequality \eqref{angXY}, apply Lemma~\ref{bases-angle} and note that $l\le m$.) \hfill $\Box$

\begin{lemma}\label{bases2} Let $\eps_1$ be the constant
defined in Lemma~\ref{bases} above. Assume that 
\begin{enumerate}
\item[{\rm (i)}] there exist orthonormal $e_1,\ldots, e_m\in \R^n$ such that
$h_i\in B^n(e_i,\delta)$ for $i=1,\ldots,m$, and $\delta
<\eps_1/2$;
\item[{\rm (ii)}] $w_i\in B^n(h_i,\eps)$ for all $i=1,\ldots,m$, and $\eps
<\eps_1/2$.
\end{enumerate}
Then the subspaces
$H=\subsp (h_1,\ldots,h_m)$ and $W=\subsp (w_1,\ldots, w_m)$ belong to $G(n,m)$, and we have $\ang (H,W) \le c_3 \eps$ with $c_3=14m\cdot 20^m$.
\end{lemma}

\medskip\noindent\textbf{Proof.} It follows from Lemma~\ref{bases} that $\dim H=\dim W = m$. We use again 
the Gram-Schmidt algo\-rithm and set $v_1=h_1$, $u_1=w_1$,
\[
\qquad v_{k+1}= h_{k+1}- \sum_{j=1}^k \frac{\langle h_{k+1},v_j
\rangle}{|v_j|^2} v_j\, ,\qquad u_{k+1}= w_{k+1}- \sum_{j=1}^k
\frac{\langle w_{k+1},u_j \rangle}{|u_j|^2} u_j\, , \qquad k+1\le
m.
\]
Then, $v_i$ and $u_i$ form orthogonal bases of $H$ and $W$,
respectively.  Inequality \eqref{vkhk} yields $t^{-1} <
|u_i|, |v_i| <  t$ with $t=10/9$. 
We now show that $s_i = \eps^{-1} |u_i-v_i|$
satisfies the assumptions of Lemma~\ref{sequences} with $a=1$ and
$b= 8$. For $k=1$ we have $s_1=
\eps^{-1}|h_1-w_1| <  1$.

Let $\phi(x)= |x|^{-2}x$.  For all $x,y$ in the annulus $\{t^{-1} \le |z|\le
t\}$ we have $|\phi(x)|\le t$ and,  for $x\not=y$,
\begin{eqnarray*}
|\phi(x)-\phi(y)|& \le & \frac{|x-y|}{|x|^2} +
|y|\biggl|\frac{1}{|x|^2}-\frac{1}{|y|^2}\biggr| 
\le  t^2|x-y| + t \biggl|\int_{|x|}^{|y|} \frac{2}{\tau^3}\,
d\tau \biggr| \\& \le & t^2(1+2t^2)|x-y|  \ <  \ 5|x-y| \, , \qquad\mbox{as $t=10/9$.}
\end{eqnarray*}
Thus, since $w_j,u_j,v_j\in \{t^{-1} \le |z|\le
t\}$, we obtain
\begin{eqnarray*}
|u_{k+1}-v_{k+1}| & \le & |h_{k+1}-w_{k+1}| +\biggl|\sum_{j=1}^k
\langle h_j,\phi(v_j)\rangle v_j - \langle w_j,\phi(u_j)\rangle
u_j \biggl| \\
& \le &  \eps +\sum_{j=1}^k \Bigl(\eps + t^2 \bigl(
|\phi(v_j)-\phi(u_j)|+|v_j-u_j|\bigr)\Bigr)\ \le \  \eps +\sum_{j=1}^k \Bigl(\eps + 6t^2 |v_j-u_j|\Bigr)\,.
\end{eqnarray*}
Hence, $$s_{k+1}=\eps^{-1} |u_{k+1}-v_{k+1}| \le (k+1) +
b\sum_{j=1}^k s_j$$ for each $b\ge 6t^2$, in particular for $b=8$. 
Therefore certainly $s_k\le 20^k$, $k=1,\ldots, m$, by
Lemma~\ref{sequences}. Keeping in mind that $t^{-1}=\frac 9{10}\le |u_j|,|v_j|\le \frac{10}9=t$, we
obtain
\[
\biggl|\frac{u_j}{|u_j|}- \frac{v_j}{|v_j|} \biggr| =
\bigl||u_j|\phi(u_j)-|v_j|\phi(v_j)\bigr| <  6t |u_j-v_j| <  
7\cdot 20^m\eps\, .
\]
The inequality $\ang (H,W)\le c_3\eps$ follows now from
Lemma~\ref{bases-angle}. \hfill $\Box$

\medskip

The next two lem\-mata are concerned with the set
\begin{equation}
\label{SHH} S(H_1,H_2):= \{y\in \R^n\, \colon \dist (y,H_i)\le
1\quad \mbox{for $i=1,2$}\},
\end{equation}
where $H_1\not= H_2\in G(n,m)$ form a small angle so that
$\pi_{H_1}$ restricted to $H_2$ is bijective. Since $\{y\in \R^n\,
\colon \dist (y,H_i)\le 1\}$ is convex, closed and centrally
symmetric\footnote{The term \emph{central symmetry} is used here for central symmetry with respect to $0$ in $\R^n$.} for each $i=1,2$, we immediately obtain the
following:

\begin{lemma}
\label{cccs} $S(H_1,H_2)$ is a convex, closed and centrally
symmetric set in $\bbbr^n$; $\pi_{H_1}(S(H_1,H_2))$ is a convex,
closed and centrally symmetric set in $H_1\cong \R^m$.
\end{lemma}

The next lemma and its corollary provide a key tool for  bootstrap estimates
in Section~\ref{slicing}.

\begin{lemma}\label{proj-slab} Let $\eps_1>0$ and $c_2>0$ denote
the constants defined in Lemma~\ref{bases}. If $H_1,H_2\in G(n,m)$
satisfy $0<\ang(H_1,H_2) = \alpha < \eps_1$, then there exists an
$(m-1)$-dimensional subspace $W\subset H_1$ such that
\[
\pi_{H_1}\bigl(S(H_1,H_2)\bigr) \subset \{y\in H_1\colon \dist
(y,W)\le 5c_2/\alpha\}\, .
\]
\end{lemma}

\medskip\noindent\textbf{Proof.} Let $H:=H_1\cap H_2$; we have
$k:=\dim H <m$. For $i=1,2$ set $X_i =\{x\in H_i\, \colon x\perp
H\}$. Then, $H_i$ is the orthogonal sum of $H$ and $X_i$. Let $X:=X_1\oplus X_2$; by construction, $X\perp H$.
Finally, let $L$ be the orthogonal complement of $H\oplus
X=H_1\oplus H_2$ in $\R^n$, so that $\R^n$ is equal to $H\oplus
X\oplus L$, and the spaces $H,X,L$ are pairwise orthogonal.  It is now easy to see, directly by
definition, that $\ang (H_1,H_2)=\ang(X_1,X_2)$. 

\smallskip\noindent\textbf{Step 1.} We shall first show that there
exists a vector $x_1\in X_1$ such that
\begin{equation}
\label{goodray} |x_1|=5c_2 /\alpha, \qquad x_1\not \in
\pi_{H_1}(S(H_1,H_2))\, .
\end{equation}
Fix an orthonormal basis $e_1,\ldots, e_{m-k}$ of $X_1$. Since
$\ang (X_1,X_2)=\alpha$, we have $|e_j-\pi_{X_2}(e_j)|\le \alpha$.
Applying Lemma~\ref{bases} with $l=m-k$ to $X_1$ and $X_2$, we
check that the $\pi_{X_2}(e_j)$, $j=1,\ldots,m-k,$
form a basis of $X_2$, and $|e_j-\pi_{X_2}(e_j)|\ge \alpha/c_2$ for at least one
$j\in \{1,2,\ldots, m-k\}$. Assume w.l.o.g. that this is the case
for $j=1$. Thus, there are no points of $X_2$ in the interior of
$B:=B^n(e_1, \alpha/c_2)$, and therefore there are no points of
$X_2$ in the interior of the cone
\[
K:=\{y\in \R^n \, \colon\ y=tv,\ t\in\R,\ v\in B\}\, .
\]
Set $\lambda := 5 c_2/\alpha$. Then $\lambda B=B(\lambda
e_1,5)\subset K$, so that the closed ball $B(\lambda e_1,
4)\subset \mathrm{int}\, K$. Hence,
\begin{equation}
\label{V2far} B(\lambda e_1, 3) \cap \{y\in \R^n \, \colon \dist
(y,X_2)\le 1\} = \emptyset\, .
\end{equation}
Let $I$ denote the segment $\{se_1 \, \colon |s-\lambda|\le 1\}$;
we claim that $\pi_{H_1}(S(H_1,H_2))\cap I=\emptyset$. To check
this, we argue by contradiction. If $y\in S(H_1,H_2)$ and
$\pi_{H_1}(y) \in I$, then, decomposing $y=h+x+l$ where $h\in H$,
$x\in X$, and $l\in L=(H\oplus X)^\perp$, we have
\[
\pi_{H_1} (y) = h + \pi_{H_1}(x) = h +\pi_{X_1} (x) = se_1
\]
for some $s$, $|s-\lambda|\le 1$. As $H\perp X$ and $e_1\in
X_1\subset X$, this yields $h=0$ and $y=se_1 +\beta$ for some
$\beta \perp X_1$. Now, \eqref{V2far} shows that if $y\in
S(H_1,H_2)$, then we must have $|\beta|^2\ge 3^2-1^2=8$. This,
however, yields $ \dist (y, X_1) = |\beta| > 2$, which contradicts
the assumption $y\in S(H_1,H_2)$. Thus, $x_1:=\lambda e_1 =
(5c_2/\alpha) e_1$ satisfies \eqref{goodray}.

\smallskip  Now, in order to prove the existence of the desired
subspace $W\subset H_1$, consider the function $$w\mapsto g(w):=
\inf\{t>0 : tw \not \in \pi_{H_1}(S(H_1,H_2))\}\in \bbbr_+\cup\{\infty\}$$  defined
on the unit sphere in $H_1$. If $w\in H$, then $g(w)=\infty$.
Since $B^n(0,1)\subset S(H_1,H_2)$, we have $g\ge 1$ everywhere.
Note that if $g(w)=s$ then $sw\in \pi_{H_1}(S(H_1,H_2))$ and $tw \not\in \pi_{H_1}(S(H_1,H_2))$ for
every $t>s$. (Thus, $g(w)$ is the `exit time' that we need to
leave $\pi_{H_1}(S(H_1,H_2))$, travelling with unit speed from $0$ in the direction
given by $w$.)

\smallskip\noindent\textbf{Step 2.} We shall first show that
there exists a vector $w_0\in H_1$, $|w_0|=1$, such that
\[
g(w_0) = r= \inf g  <  \frac{5c_2}\alpha.
\]
Since, by Step 1, we have $g(e_1)< \lambda = 5c_2/\alpha$, it is of course enough to
show that $1\le r=\inf g$ is achieved on the unit sphere of $H_1$.
Take a sequence of unit vectors $w_i\in H_1$ such that $g(w_i)\to
\inf g$; passing to a subsequence, we can assume $w_i\to w_0$ as
$i\to \infty$. Suppose now that $g(w_0)> \inf g$. Then, for some
fixed $\eps>0$ we have $g(w_0) > g(w_i)+\eps>r=\inf g$ for all
$i\gg 1$. Consider the points $p_0=g(w_0)w_0$ and $p_i=g(w_i)w_i$  in
$\pi_{H_1}(S(H_1,H_2))$. Then, $\frac{p_i-p_0}{|p_i-p_0|}\to - w_0$ as $i\to \infty$. By definition of $r$ and convexity,
\[
\pi_{H_1}(S(H_1,H_2)) \ \supset \ \mathrm{conv}\, \Bigl(\{p_0\}\cup (B^n(0,r)\cap H_1)\Bigr)\, .
\]
Thus, for all $i$ such that 
$g(w_i)<r+(\eps/2)$ and 
$\ang(w_i,w_0)<\arccos (r/(r+\eps))-\arccos(r/(r+\frac\eps 2))$ the point $p_i$ is \emph{in the interior of} $ \mathrm{conv}\, \bigl(\{p_0\}\cup (B^n(0,r)\cap H_1)\bigr)$. Then, however, by definition of $g$ we obtain $g(w_i)>|p_i|=g(w_i)$, a contradiction which shows that $g(w_0)=\inf g$.

\smallskip\noindent\textbf{Step 3.} $W=\{y\in H_1\colon y\perp w_0\}$
satisfies the desired condition. Note that $W$ is
chosen so that the set $F:=\{y\in H_1\colon \dist (y,W)\le r=\inf
g\}$ is the `narrowest strip in $H_1$' containing $\pi_{H_1}(S(H_1,H_2)) $.

Indeed, if there was a point $y\in \pi_{H_1}(S(H_1,H_2)) \setminus F$, then, taking
the straight line through $y$ and $y_0= rw_0\in \partial F\cap
\pi_{H_1}(S(H_1,H_2)) $ (with $r=\inf g$), one could easily reach a contradiction:
take a unit vector $v$ in $\mathrm{span}\, (y,y_0)$, $v\perp
y-y_0$, and use convexity of $\pi_{H_1}(S(H_1,H_2)) $ to show that then
$g(v)<r=g(w_0) =\inf g$ 
(for otherwise the straight segment connecting the points $y$ and $g(v)v$ contained in $\pi_{H_1}(S(H_1,H_2)) $ would intersect
the ray $\{y_0+tw_0: t>0\}$ contradicting the definition of $g(w_0)$).

This completes the proof of Lemma~\ref{proj-slab}. \hfill $\Box$

\smallskip

The next Lemma is practically obvious.

\begin{lemma} \label{strip-ball} Suppose that $H\in G(n,m)$ and
a set $S'\subset H$ is contained in $\{y\in H\colon \dist(y,W)\le
d\}$ for some $d>0$, where $W$ is an $(m-1)$-dimensional subspace
of $H$. Then
\[
\H^m\bigl(S'\cap B^n(a,s) \bigr) \le 2^m s^{m-1}d\,
\]
for each $a\in H$ and each $s>0$.
\end{lemma}

\smallskip\noindent\textbf{Proof.} Decomposing each $y\in S'\cap
B^n(a,s)$ as $y= \pi_W(y) + (y-\pi_W(y))$, one sees that $S'\cap
B^n(a,s)$ is contained in a rectangular box with $(m-1)$ sides
parallel to $W$ and of length $2s$ and the remaining side
perpendicular to $W$ and of length $2d$.\hfill $\Box$

\begin{lemma}\label{proj-meas}
If two planes $H_1,H_2\in G(n,m)$ satisfy $\ang
(H_1,H_2)\le \eps<m^{-1}2^{-m}$, then
\begin{equation}
\H^m(\pi_{H_1}(A))\ge (1 - m\eps 2^m) \H^m(A)
\end{equation}
for every $\H^m$-measurable set $A\subset H_2$.
\end{lemma}

\smallskip\noindent\textbf{Proof.} It is enough to prove the
inequality when $A$ is the $m$-dimensional unit cube in $H_2$, and
$\H^m(A)=1$. Fix an orthonormal basis $e_1,\ldots, e_m$ of $H_2$
and let $f_i:=\pi_{H_1}(e_i)$ for $i=1,\ldots, m$. Then, by
Hadamard's inequality,
\begin{eqnarray*}
\H^m(\pi_{H_1}(A)) & = & |f_1\wedge \ldots \wedge f_m| \\
& \ge & |e_1\wedge \ldots \wedge e_m| - \sum_{j=1}^m
|f_j-e_j|\prod_{j<i\le m} (1+|f_i-e_i|) \\
& \ge & 1 - \eps \sum_{j=1}^m (1+\eps)^{m-j} \ge 1-m\eps2^m\, .
\end{eqnarray*}

\subsection{The class of admissible sets}
\label{sec:2.3}

Let us now give a precise definition of the class of admissible
surfaces. Intuitively speaking, the energy functional $\E_q$ can
be defined for all $\Sigma\subset \bbbr^n$ compact,
$\H^m(\Sigma)<\infty$, which are a union of continuous images of
$m$-dimensional closed manifolds of class $C^1$, satisfying two
additional conditions. One of them ensures that $\Sigma$ is pretty flat
near $\H^m$--almost all its points $x$, so that we have, in a
sense, a `mock' tangent plane $H_x$ to $\Sigma$ at $x$. A priori,
$H_x$ does not even have to coincide with the classic tangent
plane. The second condition guarantees, as we shall see later,
that small $(n-m-1)$-dimensional spheres centered at $x$ and
parallel to $(H_x)^\perp$ are nontrivially linked with the surface $\Sigma$.

As we have already said in the introduction, it might be convenient to think of the following example. Assume
that $\Sigma_1,\ldots, \Sigma_N$ are embedded, compact, closed
$m$-dimensional $C^1$-submanifolds of $\R^n$. They might intersect
each other but only along sets of $m$-dimensional measure zero, so
that $\H^m(\Sigma_i\cap\Sigma_j)=0$ whenever $i\not=j$. Then, for
any bilipschitz homeomorphism $f\colon \R^n\to \R^n$,
\[
\Sigma:=f(\Sigma_1\cup \ldots \cup \Sigma_N)
\]
is an admissible surface.

\medskip

The definition of admissible surfaces involves the notion of
degree modulo 2; here are its relevant properties.

\def\deg2{\mathrm{deg}_2\, }

\begin{theorem}[Degree modulo 2] \label{degmod2}
Let $M,N$ be compact manifolds of class $C^1$  without boundary and of the same dimension $k$. Assume that $N$ is connected. There exists a unique
map
\[
\deg2 \colon C^0(M,N) \to \Z_2=\{0,1\}
\]
such that:
\begin{enumerate}
\item[{\rm (i)}] If $\deg2 g=1$, then $g\in C^0(M,N)$ is
surjective;

\item[{\rm (ii)}] If $H\colon M\times [0,1]\to N$ is continuous,
$f=H(\cdot, 0)$ and $g=H(\cdot, 1)$, then $\deg2 f =\deg2 g$;

\item[{\rm (iii)}] If $f:M\to N$ is of class $C^1$ and $y\in
N$ is an  arbitrary regular value of $f$, then
\[
\deg2 f = \# f^{-1}(y) \mod 2\, .
\]
\end{enumerate}
\end{theorem}

For a proof, see e.g. the monograph of M.W. Hirsch
\cite[Chapter 5]{hirsch}, Theorem 1.6 and the surrounding comments. Blatt gives
a detailed presentation of degree modulo 2  (even for noncompact manifolds) in his thesis
\cite{blatt-phd}.

\medskip

Now, let $\delta\in (0,1)$ and let $I$ be a finite or countable
set of indices.

\begin{definition}\label{admissible} We say that a compact set $\Sigma\subset \R^n$
is an admissible ($m$-dimensional) surface of class $\A(\delta)$
if the following conditions are satisfied.
\begin{description}
\item[(H1) Ahlfors regularity.] $\H^m(\Sigma)<\infty$ and there exists a constant
$K=K_\Sigma$ such that
\begin{equation}
\H^m(\Sigma\cap B^n(x,r)) \ge K_\Sigma r^m \qquad\mbox{for all
$x\in \Sigma$, $0<r\le \diam \Sigma$.}
\end{equation}
\item[(H2) Structure.] There exist compact, closed $m$-dimensional
manifolds $M_i$ of class $C^1$ and continuous maps $f_i\colon
M_i\to \R^n$, $i\in I$, where $I$ is at most countable, such that $\Sigma = \bigcup_{i\in I}
f_i(M_i) \cup Z$, where $\H^m(Z) = 0$.

\item[(H3) Mock tangent planes and $\delta$-flatness.]
There exists a dense subset $\Sigma^\ast\subset \Sigma$ with the
following property: $\H^m(\Sigma\setminus \Sigma^\ast)=0$ and for
each $x\in \Sigma^\ast$ there is  an $m$-dimensional plane $H=H_x\in G(n,m)$ and a
radius  $r_0=r_0(x)>0$ such that
\begin{equation}\label{d-flat}
|y-x-\pi_{H}(y-x)| < \delta |y-x| \qquad\mbox{for each $y\in
B^n(x,r_0)\cap \Sigma$,  $y\not=x$.}
\end{equation}

\item[(H4) Linking.] If $x\in \Sigma^\ast$ and $r_0(x)$ is given by 
{\rm (H3)}
above, then there exists an $i\in
I$ such that the map\footnote{Note that $\Phi_i$ is well defined, as $f_i(w)\in \Sigma$, and $z\not \in \Sigma$  by virtue of (H3).}
\[
\Phi_i\colon M_i\times \Sphere(x,r_0(x); (H_x)^\perp)\, \ni\,
(w,z) \ \mapsto\ \frac{f_i(w)-z}{|f_i(w)-z|}\, \in\,   \S^{n-1}
\]
satisfies the condition $\deg
2 \Phi_i=1$.\
(Here we use the notation $\S^{l}(\xi,\rho; P)
:=\xi +\{v\in P\colon |v|=\rho\}$ for $\xi\in \R^n$, $\rho>0$ 
and $P\in G(n,l)$.) 

\end{description}
\end{definition}

\begin{example}\rm If $\Sigma$ is a compact, connected manifold of
class $C^1$ without boundary, embedded in $\R^n$, then $\Sigma \in
\A(\delta)$ for every $\delta \in (0,1)$.

We can take $Z=\emptyset$, $I=\{1\}$, $f_1=\mathrm{id}_{\mR^n}$, and
$\Sigma^\ast=\Sigma$; (H1) and (H2) follow. It is clear that
Condition (H3) is satisfied if we choose $H_x=T_x \Sigma$ for
$x\in \Sigma$. Condition (H4) is then satisfied, too. In this
simple model case (H4) ascertains that small $(n-m-1)$-dimensional
spheres centered at the points of an embedded manifold $\Sigma$
and contained in planes that are normal to $\Sigma$ are linked
with that manifold; see e.g. \cite[pp. 194-195]{pontryagin} for
the definition of linking coefficient. (We do not assume
orientability of $M_i$; this is why degree modulo 2 is used.)

Note that if $\delta >0$ is fixed, then we are not forced to set
$H_x\equiv T_x \Sigma$; conditions (H3) and (H4) in this example
would be satisfied also if $H_x$ were sufficiently close to
$T_x\Sigma$. Thus, for given $\Sigma$ satisfying (H1) and (H2) the
choice of $H_x$ \emph{does not have to be unique.\/}
\end{example}

\smallskip
The next two examples show that we can allow $\Sigma$ to have
several $C^1$-pieces that intersect along sets of $m$-dimensional
measure zero, and are embedded away from those sets.

\begin{example}\label{C1union}
\rm If $\Sigma$ is connected, $\Sigma=\bigcup_{i=1}^N
\Sigma_i$, where $\Sigma_i$ are compact, connected manifolds of
class $C^1$ without boundary, embedded in $\R^n$, and moreover
\[
\H^m(\Sigma_i \cap \Sigma_j) = 0\qquad\mbox{for $i\not= j$,}
\]
then $\Sigma \in \A(\delta)$ for every $\delta \in (0,1)$.

The set $I$ is now equal to $\{1,\ldots, N\}$ and we set
\begin{equation}
\Sigma^\ast := \Sigma \setminus S, \qquad S:= \bigcup_{1\le i<j\le
N} (\Sigma_i \cap \Sigma_j)\, ;
\label{SigstarS}
\end{equation}
for each $x\in \Sigma^\ast$ there is a unique $i$ such that $x\in
\Sigma_i$ and  we take $H_x:= T_x\Sigma_i$. Conditions (H1) and
(H2) are clearly satisfied  with $Z=\emptyset$ and $f_i=\mathrm{id}_{\mR^n}$
for $i=1,\ldots,N$, and the verification of (H3) and (H4)
is similar to the previous example; one just has to ensure that
for $x\in \Sigma^\ast \cap \Sigma_i$ the radius $r_0=r_0(x)$ is
chosen so that $r_0< \dist \bigl(x, S)$.

\end{example}

\begin{example}\rm   Let the $M_i$,
$i\in I=\{1,\ldots, N\}$, be  compact, connected $m$-dimensional $C^1$-manifolds
without boundary. Let $f_i\colon M_i\to \R^n$ be   $C^1$-immersions, and
let $\Sigma_i=f_i(M_i)$ for $i=1,\ldots,N$. If $\Sigma =\bigcup \Sigma_i$ is connected,
\[
\H^m(\Sigma_i \cap \Sigma_j) = 0 \qquad\mbox{for $i\not=j$,}
\]
and
\[
\H^m(\{y\in \Sigma_i\colon \, \# f_i^{-1}(y) > 1\}\, ) = 0
\qquad\mbox{for all $i=1,\ldots, N$,}
\]
then $\Sigma \in \A(\delta)$ for every $\delta \in (0,1)$. We
leave the verification to the reader.

\end{example}

It is also clear that the condition that all maps $f_i$ in  the previous example be
of class $C^1$ is too strong. We can allow $\Sigma_i =f_i(M_i)$ to
have large intersections with other $\Sigma_j$ as long as the
flatness condition in (H3) is satisfied, and we need $H_x$ only
for a.e. $x\in \Sigma$. Thus, it is relatively easy to give more
examples of admissible surfaces.

\begin{example}\rm If $h\colon \R^n \to \R^n$ is a bilipschitz homeomorphism,
and we take $\Sigma$ as in Example~\ref{C1union}, then
\[
\widetilde{\Sigma}\equiv h(\Sigma) \subset \bigcap_{\delta\in (0,1)} \A(\delta)\, .
\]
Indeed, we 
then set $f_i=h\circ \mathrm{id}_{\Sigma_i}$. Let $S$ be given by \eqref{SigstarS}. To define $\widetilde \Sigma^\ast$, we use compactness and smoothness of the $\Sigma_i$ to fix a radius $r>0$ with the following property: for each $i=1,2,\ldots, N$ and each point $a\in \Sigma_i$ there is a function $g_a\colon P_a\to P_a^\perp$, $P_a=T_a\Sigma_i\in G(n,m)$,  such that $|Dg_a|\le 1$ and 
\[
\Sigma_i\cap B^n(a,r) = \mathrm{graph}\,g_a  \cap B^n(a,r).
\]
We also let $G_a (\xi )=(\xi,g_a(\xi))$ for $\xi\in P_a$.
Now, a point $y\in h(\Sigma)$ is in $\widetilde \Sigma^\ast$ if  $y\not \in h(S)$ (we exclude the intersections), 
and moreover there is $i\in\{1,\ldots,N\}$ and
an $a\in \Sigma_i$ such that 
$y=h(\xi,g_a(\xi))$ for some $\xi\in P_a\cap B^n(a,r)$ which 
is a point where $F:=h\circ G_a\colon P_a\to \R^n$ is differentiable.

It follows from Rademacher's theorem that $\widetilde\Sigma^\ast$ has full measure and is dense in $\widetilde\Sigma=h(\Sigma)$.

\heikodetail{

\bigskip

{\tt \xx\raggedright Why no ``hairs'' here, because $h$ is bilipschitz?\xx 

\yy Yes, the `bi' works here.

 }}
Condition (H1) is also satisfied, 
since bilipschitz maps distort the measure  $\H^m$ at 
most by a constant factor. To check (H3), 
one notes that as $F=h\circ G_a \colon P_a \to \R^n$ is bilipschitz, 
its differential $DF$ must have maximal rank $m$ at 
all points where it exists; it is then a simple 
exercise to check that for $y=h(x)\in\widetilde\Sigma^\ast$ 
the plane $H_y=DF(x)(P_a )$ satisfies  all requirements of Condition (H3)
for all $\delta\in (0,1).$ 
To check (H4), one can use the homotopy invariance of the degree; we leave the details to the reader.

\end{example}

\smallskip

We do not have a simple characterization of the class of admissible surfaces. However, it contains weird countably rectifiable sets, too.
\begin{example}[Stacks of spheres or cubes]\rm
(a) For $i=0,1,2,\ldots$ let $p_i=(2^{-i},0,0)\in \R^3$, $c_i=(p_i+p_{i+1})/2$, $r_i=2^{-i-2}>0$, $M_i=\Sigma_i=\S^2(c_i,r_i)\subset\R^3$ (so that the spheres $\Sigma_i$ and $\Sigma_{i+1}$ touch each other at $p_{i+1}$), and let $f_i=\mathrm{id}_{M_i}$. Set $\Sigma=\bigcup_{i=0}^\infty\Sigma_i \, \cup \, \{0\}$. Then, $\Sigma$ is an admissible surface, belonging to $\A(\delta)$ for each $\delta>0$. All points of $\Sigma$ except $0$ and the $p_i$ for $i\ge 1$ belong to $\Sigma^\ast$. For  $x\in \Sigma^\ast$, one verifies (H3) and (H4)  
just as as in Example~\ref{C1union}. Moreover, (H1) is also valid. To see this, fix $x\in \Sigma\setminus\{0\}$. If $x\in \Sigma_i$ and $r\le 2r_i$, then 
\[
\H^2(\Sigma\cap B(x,r))\ge \H^2(\Sigma_i\cap B(x,r))= \pi r^2\, ,
\]  
by the standard formula for the area of a spherical cap.
If $r>2r_i$ but $r\le \diam \Sigma=1$, then it is possible to check that the largest of all $\Sigma_j$ completely
contained in $B(x,r)\cap \Sigma$ has $r_j\in [r/6,r/2]$. Estimating $\H^2(\Sigma\cap B(x,r))$ from below by $\H^2(\Sigma_j\cap B(x,r))$, we obtain (H1) for $x\not=0$; a similar argument works for $x=0$.

%
%

(b) A modification of the above example yields the following (see Figure \ref{fig:1}): set
\[
\Sigma=\bigcup_{i=0}^\infty \biggl(\bigcup_{k=1}^{2^i} \Sigma_{i,k} \biggr)\ \cup \ Z,
\]
where $\Sigma_{i,k}$ is the surface of a cube of side length $2^{-i}$, and $Z$ is a segment of length $1$. To be more specific, $\Sigma_{0,1}=\partial\bigl( [0,1]^3\bigr)\subset \R^3$, and we let $\Sigma_{i,k}$ be a translated copy of $2^{-i}\cdot \Sigma_{0,1}=\partial [0,2^{-i}]^3$,
\[
\Sigma_{i,k} : =\partial\bigl ( [0,2^{-i}]^3\bigr) + (1-2^{-i})(e_1+2e_3) + (k-1) 2^{-i}e_2, 
\]
so that, for fixed $i$, the $\Sigma_{i,k}$ with $k=1,\ldots, 2^i$ form a layer of touching cubes stacked on top of the union of all the previous   $\Sigma_{j,s}$, $0\le j<i$ and $1\le s\le 2^j$. Finally, set $Z=\{(1,t,2)\colon t\in [0,1]$; we add this segment to the union of all $\Sigma_{i,k}$ to make $\Sigma$ closed. 

It is possible to check that if $\Sigma^\ast$ is equal to the union of the interiors of all the faces of the cubes (which is a dense set of full surface measure in $\Sigma$), then (H3) and (H4) are satisfied. Ahlfors regularity of $\Sigma$ can be checked as in (a) above.

\end{example}

\begin{REMARK}\rm 
The mock tangent planes $H_x$ are not unique in the definition 
of the class $\A(\delta)$ but \emph{a posteriori} it 
follows from Theorem~\ref{thm:bootstrap} that if 
$\Sigma\in \A(\delta)$, then  for any $x\in\Sigma^*$
there is at most {\it one} choice of the $H_x$ (up to a set of zero measure)  if one wants $\E_q(\Sigma)$ to be finite. 
Thus, finiteness of the energy is a very strong assumption: it forces us to abandon the apparent freedom of choice of the $H_x$, and forces $\Sigma$ to be a single embedded manifold, with a controlled amount of bending at a given length scale, depending only on the energy. 
\end{REMARK}

\section{Topological prerequisites}
\label{sec:3}

\setnumbers

To guarantee the existence of big projections later on, we shall
need a topological invariant, which is a version of the linking
number modulo 2.

\begin{definition}[Linking number modulo 2] Assume that $\Sigma$
is an admissible surface of class $\A(\delta)$ and $N^{n-m-1}$ is
a compact, closed $(n-m-1)$-dimensional manifold of class $C^1$,
embedded in $\R^n$ and such that $N^{n-m-1}\cap \Sigma=\emptyset$.

For each $i\in I$ and for the manifolds $M_i$ which satisfy (H2)
and (H4) of Definition~\ref{admissible}, let
\begin{equation}
\label{Gi} G_i\colon M_i\times N^{n-m-1}\, \ni\, (w,z) \ \mapsto\
\frac{f_i(w)-z}{|f_i(w)-z|}\, \in\,   \S^{n-1}.
\end{equation}
We set
\[
\lk (\Sigma^m, N^{n-m-1}) := \left\{
\begin{array}{ll}
1 & \quad\mbox{if $\ \deg2 G_i=1$ for some $i\in I$,}\\ [3pt] 0 &
\quad\mbox{if $\ \deg2 G_i=0$ for all $i\in I$.}
\end{array}\right.
\]
\end{definition}
We shall use this definition mostly in the case where $N^{n-m-1}$
is a round sphere (or an ellipsoid with ratio of axes very close
to $1$) contained in some $(n-m)$-affine plane in $\R^n$.

We need the following four properties of this invariant.

\begin{lemma}[Homotopy invariance]\label{H-inv-link}
Let $\Sigma\in \A(\delta)$ and let $N$ be a compact, closed
$(n-m-1)$-dimensional manifold of class $C^1$, and let
$N_j:=h_j(N)$ for $j=0,1$, where $h_j$ is a $C^1$ embedding of $N$
into $\R^n$ such that $N_j\cap \Sigma=\emptyset$. If there is a
homotopy
\[
H\colon N\times [0,1]\to \R^n\setminus \Sigma
\]
such that $H(\cdot,0)=h_0$ and $H(\cdot, 1)=h_1$, then
\[
\lk (\Sigma,N_0)=\lk(\Sigma,N_1)\, .
\]
\end{lemma}

\smallskip\noindent\textbf{Proof.} Note that the mappings
\[
g_{i,j} \colon M_i\times N \ni (w,z) \longmapsto
\frac{f_i(w)-h_j(z)}{|f_i(w)-h_j(z)|}\in \S^{n-1}, \qquad i\in I,
\quad  j=0,1,
\]
are such that  $g_{i,0}$ is homotopic to  $g_{i,1}$ for each $i\in
I$. Thus, the lemma follows directly from
Theorem~\ref{degmod2}~(ii). \hfill $\Box$

\begin{lemma}[Small spheres in `mock' normal planes are linked with $\Sigma$]
\label{lk=1} Assume that $\Sigma\in  \A(\delta)$ and~$x\in \Sigma^\ast$.
Then for all $r\in
(0,r_0(x))$ and for $V_x=(H_x)^\perp$ we have
\begin{equation}\label{lk2SS}
\lk (\Sigma, \Sphere(x,r;V_x)) =1,
\end{equation}
where $r_0(x)$ is the constant in Condition (H3) of 
Definition~\ref{admissible}.
\end{lemma}

\smallskip\noindent\textbf{Proof.} Due to condition (H3), each
sphere $\Sphere (x,r;V_x)$ with $r\in (0,r_0(x))$ can be deformed
homotopically to $\Sphere(x,r_0(x);V_x)$; we simply adjust the
radius, changing it linearly. Since the image of that homotopy is
disjoint from $\Sigma$, the lemma follows from (H4) and
Lemma~\ref{H-inv-link}. \hfill $\Box$

\begin{lemma}[Distant spheres are not linked]\label{lk=0}
If $\Sigma\in \A(\delta)$, $0< \eps <r < 2\eps$ and
$\dist(y,\Sigma)  > 3\eps$, then
\[
\lk (\Sigma, \Sphere(y,r;V)) = 0
\]
for each plane $V\in G(n,n-m)$.
\end{lemma}

\smallskip\noindent\textbf{Proof.} Fix an arbitrary $i\in I$. Set
$N^{n-m-1}\equiv \Sphere(y,r;V)$ and let $G_i \colon M_i\times
N^{n-m-1}\to \S^{n-1}$ be defined by \eqref{Gi}. We shall prove
that $\deg2 G_i=0$.  To this end, consider the homotopy
\[
H\colon M_i \times N^{n-m-1} \times [0,1] \to \S^{n-1}
\]
given by
\begin{equation}
H(w,z,t) =  \frac{f_i(w)-\bigl(y+(1-t)(z-y)\bigr)}
{\left|f_i(w)-\bigl(y+(1-t)(z-y)\bigr)\right|}, \qquad w\in M_i,\
z\in N^{n-m-1},\ t\in [0,1]\, .
\label{Hwzt}
\end{equation}
It is easy to see that $H$ is well defined and continuous; we have
$H(w,z,0)=G_i(w,z)$. Thus, by Theorem~\ref{degmod2}~(ii), $\deg2
G_i=\deg2 H(\cdot,\cdot, t)$ for each $t\in (0,1]$.

If $m< n-1$, then  the image of $H(\cdot, \cdot, 1)$ in $\S^{n-1}$ is the same as image of $\Sigma_i$ under the map $\xi\mapsto (\xi-y)/|\xi-y|$ which is Lipschitz in a neighbourhood of $\Sigma_i$. Since $\H^m(\Sigma_i)<\infty$, $H(\cdot, \cdot, 1)$ cannot be surjective, since the $\H^{n-1}$-measure of its image is zero.
Thus, we obtain $\deg2 H(\cdot, \cdot, 1)= 0 = \deg2 H(\cdot,\cdot,0)
=\deg2 G_i$.

If $m=n-1$, we first approximate $f_i$ by a smooth map $\tilde f_i\colon M_i\to \R^n$, so that $\|f_i-\tilde f_i\|_\infty<\eps/2$. Then, $$\widetilde G_i(w,z):=(\tilde f_i(w)-z)/|\tilde f_i(w)-z|, \qquad (w,z)\in M_i\times N^{n-m-1},$$ satisfies\footnote{Just move $f_i(w)$ to $\tilde f_i (w)$ along a segment, which avoids $N^{n-m-1}$, as $\|f_i-\tilde f_i\|_\infty<\eps/2$ and $\dist(\tilde f_i(w),N)>\eps/2$.} 
$\deg2 G_i=\deg 2\widetilde G_i$.
Next, we define $\widetilde H$ by \eqref{Hwzt} with $f_i$'s replaced by $\tilde f_i$'s. If $\widetilde H(\cdot, \cdot, 1)$ has no regular points (= points where the differential has rank equal to $n-1$), then its Jacobian is zero, and $\widetilde H(\cdot, \cdot, 1)$ is not surjective. If $\widetilde H(\cdot, \cdot, 1)$ has at least one regular point, then since $N^{n-m-1}$ consists of two distinct points $z_1,z_2$ and $\widetilde H(w,z_1,1)=\widetilde H(w,z_2,1)$, we see
each regular value of $\tilde H(\cdot, \cdot, 1)$ has an even number of preimages in $M_i\times N^{n-m-1}$. 
Hence, in either case $\deg2 \widetilde H(\cdot, \cdot, 1)= 0 = \deg2 
\widetilde H(\cdot,\cdot, 0) =\deg2 \widetilde G_i=\deg2 G_i$.

\begin{lemma}\label{link->proj}
If $\Sigma\in \A(\delta)$ and for some $y\in \R^n$, $r>0$ and
$V\in G(n,n-m)$ we have
\[
\lk (\Sigma,\Sphere (y,r;V)) =1
\]
then the disk $$\disk (y,r;V):=y+\{v\in V\colon |v|\le r\}$$ contains at least one point of
$\Sigma$.
\end{lemma}

\smallskip\noindent\textbf{Proof.} Suppose this were not the case.
Then $\dist (\Sigma, \disk (y,r;V)) > 3\eps$ for some $\eps>0$. We
deform continuously the sphere $\Sphere (y,r;V)$ to $\Sphere (y,3\eps /2; V)$,
staying all the time in $y+V$, at the distance at least $3\eps$ to
$\Sigma$. This yields
\[
\lk (\Sigma,\Sphere (y,r,;V))= \lk (\Sigma, \Sphere(y,3\eps/2; V))
 = 0
\]
by Lemma~\ref{H-inv-link} and Lemma~\ref{lk=0}, a contradiction.
\hfill $\Box$

\section{Uniform Ahl\-fors regularity}

\label{ahlfors}
\label{sec:4}

\setnumbers

\subsection{Good couples of points}

\label{sec:4.1}

We introduce here the notion of a \emph{good couple}. It expresses in a quantitative way the following rough idea: if there are two points $x,y\in \Sigma$ such that the distance from $y$ to a substantial portion of the affine planes $z+H_z$ (where $z$ is very close to $x$) is comparable to $|x-y|$, then a certain portion of energy comes \emph{only from the neighbourhood of points forming such a configuration}. Quantifying this, and iterating the resulting information in the next section, we eventually are able to pinpoint some of the local and global properties of the surface.

Recall that $Q_{H_z}$ stands for the orthogonal projection onto $(H_z)^\perp$.

\begin{definition}[Good couples]\label{couples} We say that $(x,y)\in
\Sigma \times \Sigma$ is a {\rm $(\lambda,\alpha,d)$--good couple}  if
and only if the following two conditions are satisfied:
\begin{enumerate}
\item[{\rm (i)}] $d/2 \le |x-y|\le 2d$;
\item[{\rm (ii)}] The set
\[
S(x,y;\alpha,d) :=\{z\in B^n(x,\alpha^2 d) \cap \Sigma^\ast\colon
|Q_{H_z}(y-z)|\ge \alpha d\}
\]
satisfies
\[
\H^m(S(x,y;\alpha,d)) \ge \lambda \H^m(B^m(0,\alpha^2 d)) =
\lambda \omega(m) \alpha^{2m}d^m\, .
\]
\end{enumerate}
\end{definition}
We shall be using this definition for fixed $0<\alpha,\lambda\ll
1$ depending only on $n$ and $m$. Intuitively, good couples force the energy to be large. Once we have a
$(\lambda,\alpha,d)$--good couple, then $1/\rtp$ must be $\gtrsim
\alpha d^{-1}$ on a  set in $\Sigma\times \Sigma$ of  $\H^m\otimes\H^m$-measure roughly $d^{2m}$. Thus, for $q>2m$, one cannot have $d$ small and $\E_q(\Sigma)$ small simultaneously.  We quantify that in Lemma~\ref{low-dE-bounds}.

\begin{lemma}\label{1/R-est} If $(x,y)\in \Sigma \times
\Sigma$ is a $(\lambda,\alpha,d)$--good couple with $\alpha <
\frac 12$ and an arbitrary $\lambda\in (0,1]$, then
\begin{equation}
\frac{1}{\rtp(z,w)} > \frac 19 \frac \alpha d
\end{equation}
for all $z\in S(x,y;\alpha,d)$ and $w\in B^n(y,\alpha^2 d)$.
\end{lemma}

\smallskip\noindent\textbf{Proof.} For $z,w$ as above we have
\begin{eqnarray*}
|Q_{H_z}(w-z)| & = & |Q_{H_z}(y-z) + Q_{H_z}(w-y)| \\
& \ge & \alpha d - |w-y| \qquad \mbox{by Def.~\ref{couples}~(ii)}
\\
&  > & \frac{\alpha d}{2} \qquad\mbox{as $\alpha<1/2$.}
\end{eqnarray*}
Moreover, $|w-z|\le |x-y|+|x-z|+|w-y|\le 2d + 2\alpha^2d <  3d$.
Thus, by \eqref{rtp}, 
\[
\frac{1}{\rtp(z,w)} = \frac{2\dist(w,z+H_z)}{|w-z|^2}   =
\frac{2|Q_{H_z}(w-z)|}{|w-z|^2} > \frac{\alpha
d}{(3d)^2} = \frac 19 \frac \alpha d\, .
\]

\subsection{Finding good couples and large projections}

\label{sec:4.2}

To prove uniform Ahlfors regularity, we shall demonstrate that
each $\Sigma$ with finite energy cannot penetrate certain conical
regions of $\R^n$. The construction of those regions will
guarantee that in a neighbour\-hood of each point $x\in
\Sigma^\ast$ the projections of $\Sigma$ onto suitably chosen
$m$-planes passing through $x$ are large, and a bound on the
energy will allow us to prove that such neighbour\-hoods have to
be uniformly large,   independent of the particular point  $x\in \Sigma^\ast$ we have chosen.

For a plane $H\in G(n,m)$ and $\delta\in (0,1)$ we set
\begin{eqnarray}
C(\delta,H) &:=& \{ z\in \bbbr^n \colon |Q_H(z)|\ge \delta |z|\}\, , \label{CH}\\
C_r(\delta,H) &:=& C(\delta,H) \cap B^n(0,r)\, .\label{CrH}
\end{eqnarray}
(These are closed `double cones' with `axis' equal to $H^\perp$.
Note that if $n>m+1$, then the interior of $C(\delta,H)$ and of
$C_r(\delta,H)$ is connected.) We shall also use the intersections
of cones with annuli,
\begin{equation}
A_{R,r}(x,\delta,W):=x+\mathrm{int}\, \Bigl(C_R(\delta,W)\setminus
B^n(0,r)\Bigr)\, . \label{ann-cone}
\end{equation}

\begin{lemma}[Stopping distances, good couples and large
projections]\label{mainlemma}$\phantom{a}$

There exist constants $\eta=\eta(m), \delta=\delta(m),
\lambda=\lambda(n,m) \in (0,\frac 19) $ which depend only on
$n,m$, and have the following property.

For every $\Sigma\in \A(\delta)$ and every $x\in \Sigma^\ast$
there exist $d\equiv d_s(x)>0$ and  $y\in \Sigma$ such that
\begin{enumerate}
\renewcommand{\labelenumi}{{\rm (\roman{enumi})}
}

\item $(x,y)$ is a $(\lambda,\eta,d)$--good couple;

\item for each $r\in (0,d]$ there exists a plane $H(r)\in G(n,m)$
such that
\[
\pi_{H(r)}(\Sigma\cap B^n(x,r)) \ \supset\ H(r)\cap
B^n\bigl(\pi_{H(r)}(x),r\sqrt{1-\delta^2}\bigr)\, ,
\]
and therefore $ \H^m(\Sigma\cap B^n(x,r)) \ge
(1-\delta^2)^{m/2}\omega(m) r^m$ for all $0<r \le d_s(x)$;

\item the plane $W=H(d)\in G(n,m)$ is such that
$\Sigma\cap A_{d,d/2}(x,\delta,W)=\emptyset$.

\item Each disk $D^{n-m}(z,r;W^\perp)=  z +\{v\in W^\perp\colon
|v|\le r\}$ with $z\in x+W$, $|z-x|\le d\sqrt{1-\delta^2}$, and
radius $r$ such that
\begin{equation}
 \S^{n-m-1}(z,r;W^\perp):= z +\{v\in W^\perp\colon |v|= r\} \
\subset \ A_{d,d/2}(x,\delta, W) \label{s-in-cone}
\end{equation}
contains at least one point of $\Sigma$.
\end{enumerate}
\end{lemma}
The number $d_s(x)$ is referred to as the \emph{stopping distance}. It
can be checked that the condition \eqref{s-in-cone} for the radii
of disks containing points of $\Sigma$ is equivalent to
\begin{gather}
\frac{\delta^2}{1-\delta^2} |z-x|^2\le r^2\le  d^2 - |z-x|^2
\qquad\mbox{if} \quad \frac d2 \sqrt{1-\delta^2} < |z-x| \le d
\sqrt{1-\delta^2}, \label{cone-radius1}\\
\left(\frac d2\right)^2 -|z-x|^2\le  r^2 \le  d^2 - |z-x|^2
\qquad\mbox{if} \quad |z-x|\le \frac d2 \sqrt{1-\delta^2}\,
.\label{cone-radius2}
\end{gather}

\begin{lemma}\label{low-dE-bounds} Let $\delta(m)$ be the constant of Lemma~\ref{mainlemma}. If $\Sigma\in \A(\delta)$ for some $\delta\in (0,\delta(m)]$ and $\E_q(\Sigma) < \infty$
for some $q  > 2m$, then the numbers $d_s(x)$ satisfy
\begin{equation} 
\label{low-d-bound}
d(\Sigma) := \inf_{x\in \Sigma^\ast} d_s(x) >0 \, .
\end{equation}
Moreover, we have 
\begin{equation}
\label{low-E-bound}
d(\Sigma)  \ge \left(\frac c{ \E_q(\Sigma)}\right)^{1/(q-2m)} =:R_1
\end{equation}                      
where 
\begin{equation}
c = (2\cdot 9^q)^{-1}\omega(m)^2\lambda \eta^{4m+q} 
\label{little-c}
\end{equation}   
for $\lambda=\lambda(n,m)$ and $\eta=\eta(m)$ as in Lemma~\ref{mainlemma}.
\end{lemma}
The rest of this Section is organized as follows. We prove Lemma~\ref{mainlemma} in the next subsection. Then, in subsection~\ref{sec:4.4}, 
we derive Lemma~\ref{low-dE-bounds} from Lemma~\ref{mainlemma}, and prove Theorem~\ref{thm:UAR}.

\subsection{The proof of Lemma~\ref{mainlemma}}


\label{sec:4.3}

The proof of Lemma~\ref{mainlemma} is similar to the proof of
Theorem~3.3 in our paper \cite{svdm-surfaces}. It has algorithmic
nature. Proceeding iteratively, we construct an increasingly
complicated set $S$ which is centrally symmetric with respect to
$x$ and its intersection with each sphere $\partial B^n(x,r)$ is
equal to the union of two or four spherical caps. The size of
these caps is proportional to $r$ but their position may change as
$r$ grows from $0$ to the desired stopping distance $d_s(x)$. The
interior of $S$ contains no points of $\Sigma$ but it contains
numerous $(n-m-1)$-dimensional spheres which are nontrivially
linked with $\Sigma$. Eventually, this ensures parts (ii)--(iv) of
the lemma. To find a good couple $(x,y)$, we construct $S$ so that
$\partial S\cap (\Sigma\setminus\{x\})$ is nonempty, and one of
the points in this intersection, or one of nearby points of
$\Sigma$ will be good enough for our purposes.

\medskip

The rest of this subsection is organized as follows. We first list
the conditions that have to be satisfied by $\eta$, $\delta$, and
$\lambda$. Then, we set up the plan of the whole inductive
construction and describe the first step in detail. Next, we give
the stopping criteria. Analyzing them, we demonstrate that when
the iteration stops, then (i)--(iv) of the lemma are satisfied. If
the stopping criteria do not hold, then we perform the iterative
step. Due to the nature of stopping criteria the total number of
steps in the iteration must be finite, since $\Sigma$ is compact.

 We fix a sufficiently small $\delta >0$ (to be specified
soon) and assume that $\Sigma$ belongs to the class $\A(\delta)$
of all admissible surfaces defined in Section~\ref{sec:2}. For the sake of
simplicity, we assume throughout the whole proof that
$0=x\in \Sigma^\ast$.

\medskip\noindent\textbf{The constants.} We fix the three constants $\eta=\eta(m), \delta=\delta(m),
\lambda=\lambda(n,m)$ in $(0,\frac 19) $ so that several
conditions are satisfied. We first pick $\eta$ and $\delta$ so
small that
\begin{equation}
6c_2(\delta+\eta) <  6c_3 (\delta+\eta) < \eps_1,
\label{d+e}
\end{equation}
where $\eps_1$, $c_2$ and $c_3$ denote the constants (depending
only on $m$) introduced in Lemma~\ref{bases} and
Lemma~\ref{bases2} in Section~\ref{sec:linear}.\footnote{The
stronger inequality, involving $c_3$, is needed later, in
applications in Section~\ref{sec:5}. Here, in this proof, just the condition
$6c_2(\delta+\eta)<\eps_1$ would be sufficient.}
Without loss of
generality we can also assume that
\begin{equation}
\label{delta-1} (1-\delta^2)^{m/2} > \frac 12
\qquad\mbox{and}\qquad \frac 9{10} (1-\delta^2)^{1/2}> \frac 23,
\end{equation}
and
\begin{equation}
\delta\ge 5\eta\, . \label{d5e}
\end{equation}

Next, we let $J$ be the minimal number such that there exist $J$
balls
\[
B_k:=B_{G(n,m)}(P_k,\eta^2)=\{H\in G(n,m) \colon\ang(H,P_k)\le
\eta^2\}, \qquad k=1,\ldots, J, \quad P_k\in G(n,m),
\]
that form a covering of the whole Grassmannian $G(n,m)$. Since
$\eta$ depends only on $m$, this number $J$ depends in fact only
on $n,m$. Finally, once $J$ is fixed, we let
\begin{equation}
\label{lambda} \lambda = \frac{1}{3J}.
\end{equation}

\smallskip\noindent\textbf{The construction.} Proceeding
iteratively, we shall construct three finite sequences:
\begin{itemize}
\parskip -2pt

\item of compact, connected, centrally symmetric
sets $S_0\subset T_1\subset S_1 \subset T_2\subset S_2\subset
\cdots \subset S_{N-1}\subset T_N\subset S_N\subset \R^n$,

\item of $m$-planes $H_0,\ldots,H_N$ and $H_0^\ast,\ldots,H_{N-1}^\ast\in G(n,m)$
such that the angle $\ang(H_i,H_{i}^\ast) <\eps_1$ for each
$i=0,\ldots, N-1$, where $\eps_1$ is the small constant of
Lemma~\ref{bases},
\item and of radii $\rho_0<\rho_1 < \cdots < \rho_N$, where
$\rho_N=:d_s(x)$, so $\rho_N$ will provide the desired stopping
distance for $x$ as claimed in the statement of
Lemma~\ref{mainlemma}.

\end{itemize}
Everywhere below in this subsection, we write $V_i:=H_i^\perp$ and
$V_i^\ast:= (H_i^\ast)^\perp$.

These sequences will be shown to satisfy the following properties:
\newcounter{cond}
\newenvironment{A-conditions}{\begin{list}
{{\rm (\Alph{cond})}}{\usecounter{cond}
               \setlength{\labelwidth}{4.5em}
               \setlength{\labelsep}{1em}
               \setlength{\leftmargin}{1.8cm}
               \setlength{\rightmargin}{1em}
}}{\end{list}}

\begin{A-conditions}

\item \textbf{(Diameter of $S_i$ grows geometrically).}
We have $S_i \subset B^n_{\rho_i}\!\equiv\! B^n(0,\rho_i)$ and $\diam
S_i = 2 \rho_i$ for $i=0,\ldots,N$. Moreover
\begin{equation}\label{rho-cond}
\rho_{i}> 2\rho_{i-1}\quad\Fo i=1,\ldots,N.
\end{equation}

\item \textbf{(Large `conical caps' in $S_i$ and $T_i$).}\quad
\begin{equation}\label{conical_caps1}
S_i\setminus B_{\rho_{i-1}}=C_{\rho_i} (\delta,H_i) \setminus
B_{\rho_{i-1}}\quad\Fo i=1,\ldots,N,
\end{equation}
and
\begin{equation}\label{conical_caps2}
T_{i+1}\subset B_{\rho_i}, \qquad T_{i+1}=S_i\, \cup \,
\overline{A_{\rho_i,\rho_i/2}(0,\delta, H_i^\ast)} \quad\Fo
i=0,\ldots,N-1.
\end{equation}

\item \textbf{($\Sigma$ does not enter the interior of
$S_i$ or $T_{i+1}$).}\quad
\begin{eqnarray}\label{interior_cond1}
\Sigma \cap \INT S_i& =& \emptyset\quad\Fo i=0,\ldots,N,\\
\label{interior_cond2} \Sigma\cap\INT T_{i+1}& =&
\emptyset\quad\Fo i=0,\ldots,N-1.
\end{eqnarray}
Moreover, we have 
\begin{equation}
\label{jump} \Sigma \cap
\partial B_{r} \cap C(\delta,H^\ast_{i}) = \emptyset
\qquad\mbox{for $\rho_i\le r \le 2\rho_i$,\quad $i=0,\ldots,N-1$.}
\end{equation}

\item \textbf{(Points of $\Sigma\setminus\{x\}$ on $
\partial S_i$).}\quad
The intersection $\Sigma \cap \partial B_{\rho_i} \cap
\partial S_i$ is nonempty for each $i=1,\ldots,N$.

\item \textbf{(Linking).} If $z\in H_i$ satisfies  $|z|<\rho_i
\sqrt{1-\delta^2}$ and the radius $r>0$ is chosen  such that the
$(n-m-1)$-dimensional sphere
\[
\Sphere(z,r; V_i) =z+\{v\colon\, v\in V_i, \ |v|=r\}
\]
is contained in the interior of $S_i \cap\bigl(B^n_{\rho_i}\setminus
B^n_{\rho_i/2}\bigr)$, then 
\begin{equation}
\label{link-Hi} \lk (\Sigma^m, \Sphere(z,r; V_i)) =1
\end{equation}
for $i=1,\ldots, N$.

\item \textbf{(Big projections of $B^n_{\rho_i}\cap \Sigma$ onto $H_i$).}\quad
For $t\in [\rho_{i-1},\rho_i]$, $i=1,\ldots,N$,  we have
\begin{equation}\label{bigproj}
\pi_{H_i} (\Sigma\cap B_t^n) \ \supset \ H_i \cap
B^n_{t\sqrt{1-\delta^2}}\, .
\end{equation}

\end{A-conditions}

\medskip \noindent\textbf{Start of the iteration.} We set
$S_0:=\emptyset$, $T_1:=\emptyset$, $\rho_0:=0$ and
$H_0=H_0^\ast=H_1:=H_x\in G(n,m)$,  where $H_x$ stands for the
\emph{mock tangent\/} plane  at $x=0\in \Sigma^\ast$, satisfying (H3) of
Definition~\ref{admissible}.

Moreover, we use the convention that our closed balls are defined
as
$$
B^n_r=B^n(0,r):=\overline{\{y\in\R^n:|y|<r\}}
$$
so that the closed ball $B_0$ of radius zero is the empty set.

Notice that for a complete iteration start we need to define
$\rho_1$ and $S_1$ in order to check Conditions \eqref{rho-cond}
in (A), \eqref{conical_caps1} in (B), \eqref{interior_cond1} for
$i=1$, and \eqref{link-Hi}--\eqref{bigproj} constituting
Conditions (E) and (F). All the other conditions within the whole
list are immediate for $i=0.$

We set
\begin{equation}
\label{K1} K^1_t := C_t(\delta,H_1)\,.
\end{equation}
With  growing radii $t$ the sets $K^1_t$ describe larger and
larger double cones with `axis' perpendicular to $H_1$ and fixed
opening angle which is very close to $\pi$ when $\delta$ is small.
Now we define
\begin{equation}\label{stop_1}
\rho_1:=\inf\{t>\rho_0=0:\Sigma\cap K^1_t\cap\partial
B_t\not=\emptyset\},
\end{equation}
and notice that  since $\Sigma^\ast$ satisfies \eqref{d-flat} of condition (H3) by definition, one has $\rho_1>
r_0(x)>0=2\rho_0$. This yields \eqref{rho-cond} in (A) for $i=1$.
Set $S_1:=K^1_{\rho_1}$; in other words we have
$S_1=C_{\rho_1}(\delta,H_1)\subset B^n_{\rho_1}$ with $\diam
S_1=2\rho_1$, so that all properties mentioned in (A) are
satisfied for $i=1$. Moreover, since we have adopted the
convention that $B_0$ is an empty set and $\rho_0=0$, condition
\eqref{conical_caps1} in (B) does hold for $i=1$. The definition
of $\rho_1$ guarantees that there are no points of $\Sigma$ in
$\INT S_1$, implying \eqref{interior_cond1} in (C)  for $i=1$. Condition (D)
for $i=1$ follows from the definition of $\rho_1$, as $\Sigma$ is
a closed subset of $\R^n$.

Let us now take care of (E) and (F) for $i=1$. To check (E), note
that by Lemma~\ref{lk=1} we have
\[
\lk(\Sigma^m, \Sphere (0,  r_1; V_1))  = 1
\]
for every  $ r_1>0$, $ r_1<r_0(x)=r_0(0)$.
Any sphere $\Sphere (z,r;V_1)$ with $z$ and $r$ specified in (E) for $i=1$ 
which is contained in $\INT S_1$ can be homotopically deformed to,
say, $\Sphere (0,  r_1; V_1)$ with $ r_1=r_0(x)/2$; to this
end, we just first move the base point $z$ to $0$ along the
segment $\{tz \colon\, t\in [0,1]\}$ in $H_1$, and then adjust the
radius. Notice that all  $(n-m-1)$-spheres used to define such
a homotopy are contained in $\INT S_1$ and therefore stay away
from $\Sigma$  by \eqref{interior_cond1} in (C) for $i=1$.

Thus, by Lemma~\ref{H-inv-link}, (E) follows for $i=1$. (Note that
in this first step we have even proved more. In fact, every sphere
$\Sphere (z,r;V_1)$ with $z\in H_1$, $|z|<\rho_1\sqrt{1-\delta^2}$
and radius $r$ such that $$\Sphere (z,r;V_1)\subset \INT S_1$$ is
nontrivially linked with $\Sigma$; for $i=1$ we do not have to
restrict ourselves to spheres in  $\INT S_1$ intersected with the
annulus. This restriction, however, will be necessary at later
steps.)

Invoking Lemma~\ref{link->proj}, we conclude that each 
$(n-m)$-dimensional disk $\disk (z,r; V_1)$, with $z$ and $r$ as
in (E) for $i=1$, must contain at least one point of $\Sigma$. 
Therefore,    
$$
\pi_{H_1}\bigl( \Sigma\cap \disk (z,r; V_1)\bigr)=\{z\}\Foa z\in H_1 \textnormal{
\,with $|z|<\rho_1\sqrt{1-\delta^2}$.}
$$
Since all disks $\disk (z,r; V_1)$ are contained in $B^n_{\rho_1}$ we
conclude
$$
H_1 \cap B^n_{\rho_1\sqrt{1-\delta^2}}\subset \pi_{H_1} (B^n_{\rho_1}\cap \Sigma).
$$
This is the big projection property (F) for $i=1$.

To summarize this first step,  we have defined the sets
$S_0\subset T_1\subset S_1 \subset\R^n$, and the planes
$H_0$, $H_0^\ast$ and $H_1$ which, up to now, are all identical, so
that the desired estimate for the angle $\ang(H_i,H_i^\ast)$ holds trivially  for $i=0$. We also
have defined $\rho_1>2\rho_0=0$, postponing the decision whether
$N>1$ or $N=1$. Note that we have not defined $H_1^\ast$ yet.
However, (E)--(F) do hold for $i=1$, and all those items in the
list (A)--(D) for $i=1$ which do not involve statements about
$T_2$ or $H_1^\ast$ also do hold.

We shall now discuss the stopping criteria and show how to pass to
the next step of the iteration when it is necessary.

\medskip \noindent\textbf{Stopping criteria and the iteration step.}
For the decision  whether to stop the iteration or to  continue it
with step number  $j+1$ for $j\ge 1$,  we may now assume that the
sets
$$
S_0\subset T_1\subset S_1\subset T_2\subset S_2\subset \cdots
\subset T_{j}\subset S_{j}\subset\R^n,
$$
and the $m$-planes $H_0,\ldots,H_{j}$,
$H_0^\ast,\ldots,H_{j-1}^\ast$ with $\ang(H_i,H_i^\ast)< \eps_1$
for $i=0,\ldots,j-1$, have already been defined. We also have at
this point a sequence of radii $ \rho_0=0<\rho_1 < \cdots <
\rho_{j} $ satisfying the growth condition \eqref{rho-cond} for
$i=1,\ldots, j$.

The first two conditions in (A) may be assumed to hold for
$i=0,\ldots, j$. In (B) we may suppose \eqref{conical_caps1} for
$i=1,\ldots,j$, in contrast to  \eqref{conical_caps2} which holds
only for $i=0,\ldots,j-1.$ Similarly, we may now work with
\eqref{interior_cond2} in (C) and \eqref{jump} in (D) for all
$i=0,\ldots, j-1,$ whereas \eqref{interior_cond1} in (D) can be
assumed for $i=0,\dots,j$.  The statements in (E) and (F) can be
used for $i=1,\ldots,j$.

\medskip

We are going to study the geometric situations that allow us to
stop the iteration right away; if this is the case, then we set
$N:=j$ and $d_s(x):=\rho_j=\rho_N$. Basically, there are two cases
when we can stop the construction  because  then there is a point
$y\in\bigl(B_{\rho_j}\setminus \mathrm{int}\, B_{\rho_{j}/2}\bigr)
\cap\Sigma $ such that
$(x,y)$ form a $(\lambda,\eta,\rho_j)$--good couple. In the third
case it turns out that $\Sigma \cap B_{\rho_j}\setminus \mathrm{int}\,
B_{\rho_{j}/2}$ is contained in a thin tubular neighbourhood of
some plane $H_j^\ast$, which is close to $H_j$ and \emph{very\/}
close to many of the mock tangent planes $H_z$ for points $z$ in
$B^n(x,\eta^2\rho_j)\cap \Sigma^\ast$ --- a priori, possibly even
to all of these tangent planes. When this happens, then we set
$H_{j+1}:=H_j^\ast$, define a new radius $\rho_{j+1}$, new sets
$T_{j+1}\subset S_{j+1}$ containing $S_j$,
 and finally check all the properties listed in (A)--(F).

\bigskip

The different geometric situations depend on the position of the
point where the surface hits the current centrally symmetric set
$S_j$.

\begin{description}
\item [Case 1. (First hit immediately gives a good couple.)] This occurs
if there exists at least one point $y\in \partial B_{\rho_j}\cap
C(\delta, H_j)\cap \Sigma$ such that the set $S(x,y;\eta,\rho_j)$, cf.
Definition~\ref{couples}~(ii), satisfies
\begin{equation}
\label{case1cond} \H^m(S(x,y;\eta,\rho_j)) \ge \lambda
\omega(m)\eta^{2m} \rho_j^m\, .
\end{equation}

\end{description}

If Case 1 holds, then, directly by definition, $(x,y)$ is a
$(\lambda,\eta,\rho_j)$--good couple. We then set $N:=j$,
$d_s(x)=\rho_N$, and stop the construction. It is easy to see that
all conditions of Lemma~\ref{mainlemma} are satisfied.

If Case 1 fails, then we define the new plane $H_j^\ast$ which,
roughly speaking, gives a very good approximation of a significant
portion of the mock tangent planes $H_z$ for $z$ close to $x$, and
examine the portion of $\Sigma$ contained in the closed set
\begin{equation}
\label{annulus-j} F_j:=B^n(0,2\rho_j)\setminus\mathrm{int}\,
B^n(0,\rho_j/2)
\end{equation}
to distinguish two more cases. In one of them the iteration can be
stopped in a similar way. In the second one, the whole
intersection $\Sigma\cap F_j$ might be very close to all mock tangent
planes $H_z$ so that there is no chance of finding a good
couple; we have to continue the iteration then.

We begin with the definition of $H_j^\ast$. The choice of
$\lambda$ in \eqref{lambda} comes into play here. In one of the two remaining cases
$H_j^\ast$ will become the new $H_{j+1}$. In the other case we can stop the iteration, setting $j=N$.

Fix $y\in \partial B_{\rho_{j}}(x) \cap \Sigma \cap
C(\delta,H_{j})$. Cover the Grassmannian $G(n,m)$ by finitely many
balls
\[
B_k=\{H\in G(n,m)\colon \ang(H,P_k) \le \eta^2\}, \qquad
k=1,2,\ldots, J(n,m), \quad P_k\in G(n,m).
\]
Let
\[
Y_j:= B^n_{\eta^2\rho_j}\cap \Sigma^\ast\, .
\]
Since we already can use the big projection
property\eqref{bigproj} of Condition (F) for all $i\le j$, it
follows that
\begin{equation}
\label{largemeas} \H^m (\Sigma\cap B^n_r) \ge \omega(m)
(1-\delta^2)^{m/2} r^m\qquad\mbox{for all $r\le \rho_j$.}
\end{equation}
Thus, we can estimate
\begin{eqnarray*}
\H^m(Y_j) & = & \H^m\bigl(B^n_{\eta^2\rho_j}\cap
\Sigma^\ast\bigr) \\
& \ge & \omega(m) (1-\delta^2)^{m/2}  \eta^{2m} \rho_j^m  \\
& >   & \frac 12 \omega(m) \eta^{2m} \rho_j^m \qquad\mbox{by
\eqref{delta-1}}.
\end{eqnarray*}
Now, let
\[
G_{k}:=\{z\in Y_j\colon \ang(H_z,P_k)\le \eta^2\}, \qquad
k=1,2,\ldots, J\, .
\]
Since the $G_k$ cover $Y_j$, there exists at least one $k_0\in
\{1,2,\ldots, J\}$ such that
\begin{eqnarray}
\label{k0} \H^m(G_{k_0}) & \ge & \frac{1}{J}\, \H^m(Y_j) \\
&  > &\frac 1{2J}\, \omega(m) \eta^{2m} \rho_j^m \nonumber \\
&  > & \lambda \omega(m) \eta^{2m} \rho_j^m\qquad\mbox{by
\eqref{lambda}.} \nonumber
\end{eqnarray}
We set $H_j^\ast:=P_{k_0}$, and distinguish two more cases.

\begin{description}

\item [Case 2. (Some points of $\Sigma\cap F_j$ are far from $H_j^\ast$.) ]
By this we mean that there exists a point $y \in \Sigma\cap F_j$
such that
\begin{equation}
\label{case2} |y-\pi_{H_j^\ast}(y)|\equiv |Q_{H_j^\ast}(y-x)| \ge
2\eta \rho_j\, .
\end{equation}

\end{description}
If \eqref{case2} holds, then, as in Case 1, we set $N:=j$,
$d_s(x)=\rho_N$, and stop the iteration. It remains to check that
$(x,y)$ is a $(\lambda,\eta,\rho_j)$--good couple. Condition~(i)
of Definition~\ref{couples} is clearly satisfied. To check (ii) of
that definition    we estimate for 
 each $z\in G_{k_0}\subset B_{\eta^2\rho_j}\cap
\Sigma^\ast$, using the triangle inequality,
\begin{eqnarray*}
|Q_{H_z}(y-z)| & = & |y-z-\pi_{H_z} (y-z)| \\
& = & |y-\pi_{H_j^\ast}(y) + \pi_{H_j^\ast}(y)-\pi_{H_z}(y) -z
+\pi_{H_z}(z)|\\
& \ge & 2\eta \rho_j - \ang(H_j^\ast,H_z)|y| - 2|z|
\qquad\mbox{by definition of the angle between $m$-planes}\\
& \ge & 2\eta \rho_j - \eta^2 |y| - 2 \eta^2\rho_j
\qquad\mbox{by choice of $G_{k_0}$ and $H_j^\ast$}\\
&  > & \eta\rho _j\, .
\end{eqnarray*}
(For the last inequality we just use $|y|\le 2\rho_j$ and  $\eta <
1/4$.)  Therefore, $G_{k_0}\subset S(x,y,\eta, \rho_j)$. Moreover, \eqref{k0} guarantees that $\H^m(G_{k_0})$ is
large enough. It follows that $(x,y)$ is a
$(\lambda,\eta,\rho_j)$--good couple. As before in Case~1, it is
easy to see now that all conditions of Lemma~\ref{mainlemma} are
satisfied  with $H(r)=H_i$ for all $r\in (\rho_{i-1},\rho_i]$.

\smallskip

If neither Case 1 nor Case 2 occurs, then we have to deal with
\begin{description}
\item [Case 3. (Flat position; the whole $\Sigma \cap F_j$ is
very close to $H_j^\ast$.)] This happens if and only if for each
point $y \in  \Sigma\cap F_j$ we have
\begin{equation}
\label{flat} |y-\pi_{H_j^\ast} (y)|\equiv |Q_{H_j^\ast}(y-x)| <
2\eta \rho_j\, .
\end{equation}

\end{description}

Intuitively, Case~3 corresponds to the following situation: most
points of $\Sigma\cap B_{\rho_j}$ are close to some fixed
$m$-plane which is a very good approximation of $H_z$ for
many (possibly all!) points $z\in \Sigma$ close to $x$. We then
set $H_{j+1}:=H_j^\ast$ and have to continue the iteration.

\medskip\noindent\textbf{Flat position and the passage to the next
step.} We shall first check that if Case~3 has occurred, then
\begin{equation}
\label{anglejj+1} \ang(V_j,V_{j+1})\equiv\ang(H_j,H_{j+1})\equiv
\ang(H_j,H_j^\ast) \le 3c_2(\delta+\eta) < \eps_1\, .
\end{equation}
In order to prove that this is indeed the case, we shall check
that
\begin{equation}
\label{ring-jj*} B^n(w,3(\delta+\eta)) \cap H_j^\ast
\not=\emptyset \qquad \mbox{whenever $w\in H_j$ and $|w|=1$.}
\end{equation}
Indeed, assume \eqref{ring-jj*} were false. Fix a unit vector
$w\in H_j$ such that $B^n(w,3(\delta+\eta)) \cap H_j^\ast$ is
empty. Let $z=sw$ for
\begin{equation}
s:= \frac{9}{10}(1-\delta^2)^{1/2} \rho_j
\stackrel{\eqref{delta-1}}{ >} \frac 23 \rho_j\, .\label{ineq-s}
\end{equation}
Pick
\begin{equation}
\label{ineq-r}
r:= \frac{10}9 \, \frac{\delta}{(1-\delta^2)^{1/2}}
|z| = \delta \rho_j  < \frac 19 \rho_j.
\end{equation}
Then, by \eqref{cone-radius1}, the sphere $\Sphere(z,r;V_j)$ is
contained in the interior of the intersection of $C(\delta,H_j)$
and the annulus $F_j$. Thus, we may use Condition (E),
\eqref{link-Hi} for $i=j$,  and Lemma~\ref{link->proj} to conclude that the disk $\disk
(z,r;V_j)$ contains at least one point $y_1\in \Sigma$. We also
have $y_1\in F_j$; this follows from the choice of $z$ and $r$.
Invoking \eqref{ineq-r} and \eqref{ineq-s} above, we have
\[
|y_1-z|\le r =\delta \rho_j  < 2 s\delta\, .
\]
Since $B^n(w,3(\delta+\eta)) \cap H_j^\ast=\emptyset$ and $z=sw$,
by scaling we have also
\begin{equation}
\label{Bz3} B^n(z,3s(\delta+\eta))\cap H_j^\ast = \emptyset,
\end{equation}
so that the triangle inequality gives, by \eqref{ineq-s},
\[
|y_1-\pi_{H_j^\ast}(y_1)|  > 3s(\delta+\eta) - 2 s\delta > 3s\eta
 > 2\eta\rho_j\, .
\]
This, however, is a contradiction to condition \eqref{flat} which
holds in Case~3. Hence, \eqref{ring-jj*} holds too, and for
every orthonormal basis $(e_i)\subset H_j$ the vectors
$f_i:=\pi_{H_{j+1}}(e_i)$ form a basis of $H_{j+1}$ which
satisfies $|e_i-f_i|\le 3(\delta+\eta)<\eps_1$. Lemma~\ref{bases} 
implies that
\[
\ang(H_j,H_{j+1})\equiv \ang (H_j,H_j^\ast) < c_2 \cdot
3(\delta+\eta) \ \stackrel{\eqref{d+e}}{<} \ \eps_1\, ,
\]
 which is \eqref{anglejj+1}.

As the angle $\ang(H_j,H_{j+1})=\ang(V_j,V_{j+1})$ is small, the
cones $C(\delta, H_j)$ and $C(\delta,H_{j+1})$ have a large
intersection. Indeed, for any unit vector $v\in \R^n$ with
$|\pi_{H_{j+1}}(v)|\le \theta$ we have $|\pi_{H_j}(v)| < \theta +
3c_2(\delta+\eta)$ by definition of $\ang(H_j,H_{j+1})$. Thus,
\[
|Q_{H_j}(v)|\ge |v|-|\pi_{H_j}(v)| > 1 - \theta -
3c_2(\delta+\eta) > \delta
\]
whenever $\theta < 1-\delta-3c_2(\delta+\eta) < 1-\eps_1$. In
particular, every unit vector $v\in V_{j+1}$ belongs to the
interior of $C(\delta,H_j)$.

We now define
\begin{equation}
\label{deftj+1} T_{j+1}:= S_j \ \cup \ \bigl(
C_{\rho_j}(\delta,H_{j+1}) \setminus \INT B^n_{\rho_j/2}\bigr)\, .
\end{equation}
According to \eqref{flat}, this immediately gives the missing conditions
\eqref{conical_caps2} in (B)  and \eqref{interior_cond2} in (C)
for $i=j$. To check \eqref{jump} in (C) for $i=j$, note that in
Case 3 we have
\[
|Q_{H_j^\ast}(y)|  < 2 \eta \rho_j \le 4 \eta |y|
\]
for each point of $\Sigma$ in the annulus $F_j$. However, when
$y\in C(\delta, H_j^\ast) \cap \partial B_r$ for some $\rho_j\le
r\le 2\rho_j$, then
\[
|Q_{H_j^\ast}(y)| \ \stackrel{\eqref{CH}}{\ge} \ \delta |y| \
\stackrel{\eqref{d5e}}{\ge}\ 5\eta|y|,
\]
so that $y$ cannot be a point of $\Sigma$. This gives \eqref{jump}
for $i=j$.

Now the crucial thing is to define the next radius $\rho_{j+1}$
and take care of the linking condition \eqref{link-Hi} for
$i=j+1$.


\medskip\noindent\textbf{The next radius and homo\-topies from large spheres to smaller tilted
ones.} Set \begin{equation}\label{Kj+1}
K^{j+1}_t:=C_t(\delta,H_{j+1}),
\end{equation}
and define
\begin{equation}\label{stopj+1}
\rho_{j+1}:=\inf\{t>\rho_{j}:\Sigma\cap K_t^{j+1}\cap\partial
B_t\not=\emptyset\}.
\end{equation}
 Notice that condition \eqref{jump} guarantees that $\rho_{j+1}
> 2\rho_j$. This  verifies \eqref{rho-cond} in Condition (A) for
$i=j+1.$ Now we define
\begin{equation}\label{Sj+1}
S_{j+1}:= T_{j+1}\cup (K^{j+1}_{\rho_{j+1}}\setminus \INT
B_{\rho_j}),
\end{equation}
and check that Conditions (A)--(F) are satisfied.

Indeed, $S_{j+1}\subset S_j\cup K^{j+1}_{\rho_{j+1}}\subset
B_{\rho_j}\cup B_{\rho_{j+1}}$ by Condition (A) for $i=j,$ which
implies that (A) holds for $i=j+1$ as well. Next,
$$
S_{j+1} \setminus B_{\rho_j}=K^{j+1}_{\rho_{j+1}}\setminus
B_{\rho_j}=C_{\rho_{j+1}}( \delta,H_{j+1})\setminus B_{\rho_j},
$$
since $S_j\subset B_{\rho_j}$ by Condition (A) for $i=j$. Hence
\eqref{conical_caps1} holds for $i=j+1$. The inclusion
$T_{j+1}\subset S_j\cup B_{\rho_j}\subset B_{\rho_j}$ and other
conditions involving $T_{j+1}$ have already been checked; they
follow directly from the definition of $T_{j+1}$; see \eqref{deftj+1}. 
Using
\eqref{interior_cond1} for $i=j$ and the definition of
$\rho_{j+1}>2\rho_{j}$ in \eqref{stopj+1} we infer that
\eqref{interior_cond1} holds for $i=j+1$, and
\eqref{interior_cond2} for $i=j$. We also have
\begin{equation}
\label{newjump} \Sigma \cap \partial B_r \cap C(\delta,H_{j+1}) =
\emptyset \qquad\mbox{for each $r\in [\rho_j,\rho_{j+1})$.}
\end{equation}
Also Condition (D) follows directly from the definition of
$\rho_{j+1}$.

Now we turn to the proof of the linking condition,
\eqref{link-Hi}, for $i=j+1.$

The definition \eqref{stopj+1} of $\rho_{j+1}$ implies that each
sphere $\Sphere (z,r_0,V_{j+1})$ where $z\in H_{j+1}$,
$|z|<\rho_{j+1} \sqrt{1-\delta^2}$ and the radius $r_0$ is such
that
\[
M_0:=\Sphere (z,r_0,V_{j+1}) \subset \mathrm \INT S_{j+1} \cap
\{y\in \R^n \colon\, \rho_{j+1}/2<|y|<\rho_{j+1}\}
\]
can be homotopically deformed to
\[
M_1:=\Sphere (0,r_1;V_{j+1}) \subset \mathrm \INT S_{j+1} \cap
\{y\in \R^n \colon\, \rho_{j+1}/2<|y|<\rho_{j+1}\}, \qquad
r_1^2:=|z|^2+r_0^2,
\]
without meeting any points of $\Sigma$, so that the linking
invariant used in \eqref{link-Hi} is preserved. One of the
possible homotopies is to move the base point $z$ to 0 along the
segment $z(t)=(1-t)z$, $t\in [0,1]$, at the same time increasing
the radius from $r_0=r(0)$ to $r_{1}=r(1)$ so that
\[
\rho^2:= |z(t)|^2 + r(t)^2
\]
remains constant for all $t\in [0,1]$; in this way, we simply
slide the $(n-m-1)$-dimensional spheres along the surface of a
fixed $(n-1)$-sphere, staying all the time in the interior of
$S_{{j+1}}$ intersected with the annulus $\{y\in \R^n \colon\,
\rho_{j+1}/2<|y|<\rho_{j+1}\}$. By Lemma~\ref{H-inv-link} we have
\begin{equation}
\lk (\Sigma, M_0) = \lk(\Sigma, M_1)\, . \label{link-M0M1}
\end{equation}
Next, we may homotopically deform the sphere $M_1$ to another 
sphere of radius $r_2$,
\[
M_2:= \Sphere (0,r_2;V_{j+1}), \qquad r_2=\frac 89 \rho_j\in
(\rho_j/2, \rho_j)\, .
\]
We just shrink the radius linearly, staying all the time  in the
$(n-m)$-dimensional subspace $V_{j+1}$. It is clear that all the
flat spheres realizing this homotopy $M_1\sim M_2$ stay in the
interior of $S_{j+1}$ (by \eqref{interior_cond1} for $i=j+1$ and the definition of $T_{j+1}$ in \eqref{deftj+1}) and do not contain any points of $\Sigma$,
so that, again by Lemma~\ref{H-inv-link},
\begin{equation}
\lk (\Sigma, M_1)  = \lk(\Sigma, M_2)\, . \label{link-M1M2}
\end{equation}
But $M_2$ can be homotopied --- still in the interior of $S_{j+1}$
--- to another sphere,
\[
M_3:=\Sphere (0,r_2;V_{j}),
\]
which has the same radius  $r_2$ but is slightly tilted; 
therefore,
\begin{equation}
\lk (\Sigma, M_2)  = \lk(\Sigma, M_3)\, . \label{link-M2M3}
\end{equation}
 To check this, we perform two steps. First we move each
point $y$ of $M_2\subset{V_{j+1}}$ along the segment that joins
$y$ to its projection $\pi_{V_j}(y)$. This gives an ellipsoid
which is nearly spherical and has all axes at least
$(1-\eps_1)r_2$ because of the condition \eqref{anglejj+1} for the
angle between $V_j$ and $V_{j+1}$. Next, we continuously blow up
this ellipsoid, moving each of its points along the rays that
emanate from $0$ to points of $M_3$. Because of the smallness
condition \eqref{d+e} for the constants that we use, each segment
$I_y$ with one endpoint at $y\in M_2$, $|y|=\frac 89 \rho_j$,
and the other at $\pi_{V_j}(y)$ {} is certainly contained in
the interior of $S_j$ (i.e. far away from $\Sigma$), as
\[
|y-\pi_{V_j}(y)| \le \ang(V_j,V_{j+1}) |y| \le 3c_2(\delta+\eta)
r_2 < 3 c_2(\delta+\eta)\rho_j < \frac 13 \rho_j\, ,
\]
so that $|\pi_{V_j}(y)|> \frac 89 \rho_j-\frac 13 \rho_j>\frac 12 \rho_j> 2\rho_{j-1}$.
Thus, invoking Lemma~\ref{H-inv-link} one more time, and applying
the inductive assumption, i.e. the linking condition
\eqref{link-Hi} for $i=j$, we finally obtain
\[
\lk(\Sigma,M_0) = \lk(\Sigma,M_3) = 1
\]
This gives \eqref{link-Hi} of (E) for $i=j+1$.

It is now easy to establish the big projection property of (F) for
$i=j+1$. We do this as in the first step of the proof: invoking
Lemma~\ref{link->proj}, we conclude that each flat
$(n-m)$-dimensional disk $\disk (z,r; V_{j+1})$, with $z$ and $r$
as in (E) for $i=j+1$, must contain at least one point of
$\Sigma$. 
Therefore,  
$$
\pi_{H_{j+1}}(\Sigma\cap \disk(z,r; V_{j+1}))=\{z\}\Foa
z\in H_{j+1}\textnormal{\, with $|z|<\rho_{j+1}\sqrt{1-\delta^2}$}.
$$
All disks $\disk(z,r; V_{j+1})$ are contained in $B^n_{\rho_{j+1}}$ so 
that
$$
H_{j+1} \cap B^n_{\rho_{j+1}\sqrt{1-\delta^2}}\subset \pi_{H_{j+1}} (B^n_{\rho_{j+1}}\cap \Sigma).
$$
This gives \eqref{bigproj} in (F) for $i=j+1$, and finishes the
proof  of all conditions in the list (A)--(F) in the iteration
step.

Since we have established Condition (E) in the iteration step and
\eqref{rho-cond} holds, too, we can deduce that Case 3 can happen
only finitely many times, depending on the position $x$ on
$\Sigma$ and on the shape and size of $\Sigma$:
$$
 \diam\Sigma\ge\rho_i>2\rho_{i-1}>\cdots > 2^{i-1}\rho_1>2^{i-1}r_0(x),
$$
whence  the maximal number of iteration steps is bounded by $$
 1+\log(\diam\Sigma/r_0(x))/\log 2\, .
$$
This concludes the consideration of Case 3, and the whole proof of
Lemma~\ref{mainlemma} .\hfill $\Box$

\subsection{Bounds for the stopping distances and uniform Ahlfors regularity}

\label{sec:4.4}

We shall now derive Lemma~\ref{low-dE-bounds} and Theorem~\ref{thm:UAR} from Lemma~\ref{mainlemma}. This is a relatively easy task at this stage. We shall just relay on estimates for the $\E_q$-energy in the neighbourhooud of a good couple $(x,y)\in \Sigma\times\Sigma$.

\medskip\noindent\textbf{Proof of Lemma~\ref{low-dE-bounds}.} Since $\A(\delta) \subset \A(\delta')$ for $\delta\le \delta'$, we assume from now on that $\delta=\delta(m)$ is the constant of Lemma~\ref{mainlemma}.  

Fix $\eps>0$ small (to be specified later on). Assume that $d(\Sigma)=\inf_{\Sigma^\ast} d_s<\eps$ and select a point $x\in\Sigma^\ast$ such that $ 
d_s(x)<\eps$. Use Lemma~\ref{mainlemma} to select a $(\lambda,\eta,d)$--good couple $(x,y)\in \Sigma\times\Sigma$. Let
\[
S:=S(x,y;\eta,d_s(x))
\]
be as in Definition~\ref{couples}~(ii), and let $B:=B(y,\eta^2 d_s(x) )$. Applying~Lemma~\ref{1/R-est} we estimate
\begin{eqnarray*} 
\E_q(\Sigma) & \ge & \int_{S}\int_{\Sigma\cap B}	\left(\frac{1}{\rtp(z,w)}\right)^q\, d\H^m(w)\, d\H^m(z) \\
&  > & \H^m(S)  \H^m(\Sigma\cap B) \left(\frac \eta{9d_s(x)}\right)^q \qquad\mbox{by Lemma~\ref{1/R-est}} \\
&\ge & \lambda\omega(m) \eta^{2m}d_s(x)^m \cdot K_\Sigma \eta^{2m} d_s(x)^m  \left(\frac \eta{9d_s(x)}\right)^q \qquad\mbox{by Definitions~\ref{couples} and~\ref{admissible}}\\
& = & K_\Sigma 9^{-q} \lambda{\eta}^{4m+q}  d_s(x)^{2m-q},
\end{eqnarray*}
which implies
\begin{eqnarray*}
 \eps^{q-2m} & > & d_s(x)^{q-2m}\\
&  > & K_\Sigma \lambda {\eta}^{4m+q} 9^{-q} \E_q(\Sigma)^{-1},
\end{eqnarray*}
a contradiction for
\[
\eps:= \left( \frac 12 K_\Sigma \, \lambda\,  {\eta}^{4m+q} 9^{-q} \E_q(\Sigma)^{-1}  \right)^{1/(q-2m)}\, .
\] 
This proves the first part of the  lemma. 

Now, for an arbitrarily small $\sigma\in (0,1)$ pick $x_0\in\Sigma^\ast$ such that $d(\Sigma)\le d_0=d_s(x_0)<(1+\sigma)d(\Sigma)$. Select $y_0\in\Sigma$ so that $(x_0,y_0)$ is a $(\lambda,\eta,d_0)$--good couple. We have $d_s(y_0)\ge d(\Sigma)>d_0/(1+\sigma)$, so that by Lemma~\ref{mainlemma} (ii) 
\[
\H^m\bigl(\Sigma \cap B^n(y,r)\bigr) \ge (1-\delta^2)^{m/2} \omega(m) r^m \ge \frac 12 \omega(m)r^m
\]
certainly holds for $r=\eta^2d_0 <d_0/(1+\sigma)$ since $\eta\ll 1$ by \eqref{d+e}. Estimating the energy one more time, as before, we obtain
\begin{eqnarray*} 
\E_q(\Sigma) & \ge &  \int_{ S(x_0,y_0;\eta, d_0)}\int_{\Sigma\cap B(y_0,\eta^2 d_0)}	\left(\frac{1}{\rtp(z,w)}\right)^q\, d\H^m(w)\, d\H^m(z) \\
&> & \frac\lambda {2\cdot 9^q} \omega(m)^2 \eta^{4m+q}\, d_0^{2m-q} \qquad\mbox{by Lemma~\ref{1/R-est}} .
\end{eqnarray*}
Thus,
\[
(1+\sigma)^{q-2m}d(\Sigma)^{q-2m} > d_0^{q-2m} >  c \E_q(\Sigma)^{-1}\, , 
\]
where $c=(2\cdot 9^q)^{-1}\omega(m)^2\lambda \eta^{4m+q}$, as  in \eqref{little-c}. Letting $\sigma\to 0$, we obtain \eqref{low-E-bound} and conclude the whole proof.\quad $\Box$ 

\bigskip\noindent\textbf{Proof of Theorem~\ref{thm:UAR}.} By the lower bound  \eqref{low-E-bound} for stopping distances, the inequality
\[
\H^m(\Sigma\cap B(x,r)) \ge (1-\delta^2)^{m/2} \omega(m)r^m \ge \frac 12 \omega(m)r^m
\]
holds for each $x\in \Sigma^\ast$ and each $r\le d(\Sigma)\le d_s(x)$. By density of $\Sigma^\ast$ in $\Sigma$, we obtain    $$\H^m(\Sigma\cap B(x,r))\ge \frac 12 \omega(m)r^m$$ for \emph{all} $x\in \Sigma$ and $r\le d(\Sigma)$. This implies Theorem~\ref{thm:UAR} with
\[
a_1:=\left(\frac{\lambda\, \omega(m)^2\, \eta^{4m+q}}{2\cdot 9^q}\right)^{1/(q-2m)},
\]
where $\lambda=\lambda(n,m)$ and $\eta=\eta(m)$ are the constants introduced in Lemma~\ref{mainlemma}.

\begin{REMARK}\rm  As we have already mentioned in the introduction, the proof above yields a result 
which is stronger than the formal statement of 
Theorem~\ref{thm:UAR}. In fact, the result holds also for all 
$\Sigma\in\A(\delta)$ with $0<\delta\le  \delta(m) $, where $\delta(m) $ 
is the positive constant of Lemma \ref{mainlemma}, and this is a wider class of sets than the one we used in the introduction.
\end{REMARK}


\section{Existence of tangent planes}

\label{sec:5}

\setnumbers

In this section  we prove that  for each point $x\in \Sigma$ there exists a plane $T_x\Sigma\in G(n,m)$ such that $\dist(x', x+T_x\Sigma)=o(|x'-x|)$ for $x'\in \Sigma$, $x'\to x$. Moreover, the mapping $ x\mapsto
T_x\Sigma $ is of class $C^\kappa$, $\kappa=(q-2m)/(q+4m)>0$.
 A posteriori it turns out that if $\delta>0$ is small
enough and $\Sigma\in \A(\delta)$ is an admissible surface with
$\E_q(\Sigma)< \infty$ for some $q>2m$, then the mock tangent
planes $H_x$ defined a.e. on $\Sigma$ must coincide with the
classically understood $T_x\Sigma$.

The idea is to combine the results of the previous section with
energy bounds and show that the P.~Jones' $\beta$-numbers of
$\Sigma$ satisfy a decay estimate of the form
$\beta_\Sigma(x,r)\lesssim E^{1/(q+4m)} r^\kappa$. This alone
would not be enough, but we already know that $\Sigma$ has big
projections. Adding this ingredient, we are able to prove that
$\Sigma$ is in fact a $C^{1,\kappa}$-manifold.  Moreover, in each
ball of radius $\approx \E_q(\Sigma)^{-1/(q-2m)}$ centered at
$x\in \Sigma$ the surface $\Sigma$ is a graph of a $C^{1,\kappa}$
function $f\colon P\to P^\perp$ over $P=T_x\Sigma\in G(n,m)$.

The core of this  section is formed by an iterative construction,
presented in Section~\ref{iteration}, which yields the existence
of tangent planes and estimates for their oscillation. At each
iteration step, we need to check that the  $\beta$-numbers
decrease sufficiently fast as the length scale shrinks to zero. At
the same time, we have to guarantee that the linking conditions
which imply the existence of big projections are also satisfied.
To make the presentation of that proof easier to digest, we
introduce an ad-hoc notion of \emph{trapping boxes\/}
(Section~\ref{boxes}) and prove an auxiliary lemma which is then
used in the iteration.

\subsection{Trapping boxes}

\label{sec:5.1}

\label{boxes}

Everywhere in this section $R_1$ denotes the radius specified in
Theorem~\ref{thm:UAR} ascertaining the uniform Ahlfors
regularity of surfaces with bounded energy.

 For the rest of the whole section, we fix $\delta,\eta>0$
small so that \eqref{d+e} is satisfied and all claims of
Lemma~\ref{mainlemma} are fulfilled.

\begin{definition}\label{TB} Assume that $\Sigma\in \A(\delta)$, $x\in \Sigma$, $0< r <
R_1$, $\theta\in (0,\delta]$ and $H\in G(n,m)$. We say that a
closed set $F\subset B^n(x,r)$ is \emph{a $(\theta,H)$-trapping
box for $\Sigma$ in $B^n(x,r)$\/} if and only if the following
conditions are satisfied:

\begin{enumerate}

\item[{\rm (i)}]  $\Sigma\cap B^n(x,r) \, \subset\, F$;

\item[{\rm (ii)}] $\{y\in B^n(x,r)\colon \dist(y, x+H)\le \theta
r\}\, \subset\, F$;

\item[{\rm (iii)}] if $z\in x+H$ satisfies $|z-x|<
(1-\theta^2)^{1/2} r$, then there exists a $t>0$ such that
$t^2+|z-x|^2 < r^2$, the sphere $\Sphere (z,t;H^\perp)$ is
contained in the interior of $B^n(x,r)\setminus F$ and
\[
\lk (\Sigma,\Sphere (z,t;H^\perp))=1\, .
\]
\end{enumerate}
\end{definition}

Thus, informally, a trapping box is a subset of $B=B^n(x,r)$ which
is at least as large as a cylindrical neighbourhood of $x+H$ in
$B$ (of  size specified by the parameter $\theta$), and gives us
some control of the location of $\Sigma\cap B$ and of its projections
onto $H$.

If $\Sigma\in \A(\delta)$, $x\in \Sigma^\ast$ and $H= H_x$ is
given by Condition (H3) of Definition~\ref{admissible}, then
--- for radii $r<r_0(x)$
--- a simple example of a trapping box is provided by  the cylinder
\[
\{y\in B^n(x,r)\colon \dist(y, x+H)\le \delta r\}.
\]
It satisfies all conditions of Definition~\ref{TB} for $\theta=\delta$;
in particular, Lemma \ref{lk=1} guarantees Condition (iii).

Another example is given by the following.

\begin{proposition}
\label{startbox} Let $\delta(m)$ be the small constant of
Lemma~\ref{mainlemma}. 
Assume that $\Sigma\in \A(\delta)$, $\delta\in (0,\delta(m)]$,
$\E_q(\Sigma)\le E$, and $R_1$ denotes the radius specified in
Theorem~\ref{thm:UAR}.
 Then, for each $x\in \Sigma$ and each $r\in (0,R_1)$ there exists a plane
$H\in G(n,m)$ such that
\begin{equation}\label{startbox-delta}
F:= \{y\in B^n(x,r)\colon \dist(y, x+H)\le \delta r\}  \cup
B^n(x,r/2)
\end{equation}
is a $(\delta,H)$-trapping box for $\Sigma$ in $B^n(x,r)$.
\end{proposition}

\medskip\noindent\textbf{Proof.} One can check  that conditions (A)--(F) stated at the beginning of the proof of Lemma~\ref{mainlemma} combined with the lower bound for stopping
distances obtained in Lemma~\ref{low-dE-bounds} imply the
statement of Proposition~\ref{startbox} for all points $x\in
\Sigma^\ast$. (To see this, look at condition \eqref{flat} of Case~3, which is the only case when the iterative construction is continued. It has been designed in such a way that the union of $\{y\in B^n(x,r)\colon \dist(y, x+H_j^\ast)\le 2\eta \rho_j\}$ and $B^n(x,\rho_j/2)$ be a trapping box for $\Sigma$ in $B(x,2\rho_j)$; condition (E), cf. \eqref{link-Hi},  implies the existence of \emph{many} spheres linked with $\Sigma$ so that (iii) of Definition~\ref{TB} is also satisfied.  Since $\eta\le\delta/5$ by \eqref{d5e}, the claim of the proposition
holds for all $r\in [\rho_j,\rho_{j+1}]$ with $H=H_j^\ast$, and we can certainly increase $r$ up to the infimum  $d(\Sigma)$ of all stopping distances, which satisfies $d(\Sigma)\ge R_1$ by Lemma~\ref{low-dE-bounds}.)


Assume now that $x\not\in \Sigma^\ast$. Fix $r\in (0,R_1)$ and select a sequence $x_l\in \Sigma^\ast$, $x_l\to x$ as $l\to \infty$. For each $l=1,2,\ldots$, let $H_l$ whose existence is given by the statement of the proposition at points $x_l\in \Sigma^\ast$. Passing to a further subsequence, we may  assume that $H_l\to H\in G(n,m)$ as $l\to \infty$. The trapping boxes $F_l$ corresponding to $x_l$ and $H_l$ via \eqref{startbox-delta} converge then in Hausdorff distance to a closed set $F$ given by \eqref{startbox-delta} for $x$ and $H$. Since $\Sigma$ is closed, $\Sigma \cap B^n(x,r)$ must be contained in $F$. Condition (ii) of Definition~\ref{TB} is trivially satisfied, and condition (iii) is easily verified by using homotopical invariance of the linking number as we already did before (one has to slightly tilt the spheres in $B^n(x,r)\setminus F$ to obtain spheres in $B^n(x_l,r)\setminus F_l$). \hfill $\Box$

\medskip

The main idea of this  section is to show that once we have a
trapping box of the form \eqref{startbox-delta}, possibly with
$\delta$ replaced by some smaller number $\theta>0$, then, under a
certain balance condition  for $\varphi$, $r$ and the energy of
$\Sigma$, we can perturb the plane $H$ slightly to a new position
$H_1$ and find a smaller, cylindrical $(\varphi,H_1)$-trapping
box. We make this precise in the next subsection.

\subsection{Energy bounds and trapping boxes in small scales}

\label{sec:5.2}

We introduce two new constants
\begin{equation}\label{c4c5}
c_4:= 3(c_3 +1), \qquad c_5:=
\frac{16m \cdot 9^q}{\omega(m)^2}.
\end{equation}
Recall from Lemma \ref{bases2} 
that the constant  $c_3= 14m\cdot 20^m$ depends on $m$ only.

\begin{lemma}\label{shrink-box} Assume that $H\in G(n,m)$, $x \in
\Sigma$, $0<r<R_1$, $0<\theta\le \delta$,  $q>2m$. Let $\Sigma\in
\A(\delta)$, $\delta\in (0,\delta(m)]  $  
be an admissible surface with $\E_q(\Sigma)\le  E$.
Suppose that
\[
F_{\theta,r} (H) := \{y\in B^n (x,r) \colon \dist(y, x+ H)\le \theta r\}
\cup B^n(x,r/2)
\]
is a $(\theta,H)$-trapping box for $\Sigma$ in $B^n(x,r) $. If
$0<\varphi< 1/(6c_4)$ satisfies the \emph{balance condition\/}
\begin{equation}\label{balance}
\varphi^{4m+q}r^{2m-q} \ge c_5 E\, ,
\end{equation}
then there exists a plane $H_1\in G(n,m)$ such that
\begin{enumerate}
\item [{\rm (i)}] $\ang (H,H_1)\le 2c_2\theta$;

\item [{\rm (ii)}] The cylinder
\begin{equation}\label{newbox}
F:= \{y\in B^n_{2r} \colon \dist(y,  x+ H_1)\le c_4\varphi \cdot 2r\}
\end{equation}
is a $(c_4\varphi,H_1)$-trapping box for $\Sigma$ in $B^n(x,{2r})$.
\end{enumerate}

\end{lemma}
The main point is that once we fix  a finite energy level $E$,
and $r$ sufficiently small, then the condition $q>2m$ guarantees
that there are numbers $\varphi>0$ which satisfy the balance
condition \eqref{balance} and are such that $c_4\varphi$ is (much)
smaller than $\theta$. Since the angle $\ang( H,H_1)$ is
controlled due to (i), the lemma can be applied iteratively. This
will be done in the next subsection.

\begin{REMARK}\label{2pt-plane}
If we fix an arbitrary point $y\in ( B^n(x,r)\cap \Sigma) \subset F_{\theta,r}(H)$ such
that $ \frac{9}{10}(1-\theta^2)^{1/2} r\,\le |y-x|< r$, then the
plane $H_1$ in Lemma~\ref{shrink-box} can be chosen so that
$y-x\in H_1$,  as can be seen from the first step of the following proof.
\end{REMARK}

\medskip\noindent\textbf{Proof of Lemma \ref{shrink-box}.}
Fix an arbitrary orthonormal basis $(e_1,\ldots,e_m)$
of $H$ and let
\[
d:=\frac{9}{10}(1-\theta^2)^{1/2} r\, .
\]
Since $\theta\le \delta$, we have $d> \frac 23 r$ by
\eqref{delta-1}. Set $z_i=de_i$, $i=1,\ldots,  m$.

\smallskip\noindent\textbf{Step 1. Choice of $H_1$.}
Using Condition~(iii) of Definition~\ref{TB},
Lemma~\ref{H-inv-link} and Lemma~\ref{link->proj}, we conclude
that each disk
\[
D_i:=\disk (z_i,\theta r;H^\perp), \qquad i=1,\ldots,m,
\]
contains\footnote{There are points of $\Sigma$ in all disks with slightly larger radii, and $\Sigma $ is closed.} a point $y_i\in \Sigma$. Set $H_1=\subsp(y_1,\ldots,
y_m)$. Letting $h_i=d^{-1}y_i$, we use $\theta\le \delta$ and
 \eqref{d+e} to estimate
\[
|h_i-e_i|=d^{-1}|y_i-z_i|\le \frac{\theta r}{d} <   2\theta <
\frac{\eps_1}2,
\]
and invoke Lemma~\ref{bases} to obtain $\ang(H,H_1)<
2c_2\theta$. (This initial step of the proof shows why
Remark~\ref{2pt-plane} is satisfied. We can work with  an
orthonormal basis $e_i$ such that $e_1=\pi_H(y)/|\pi_H(y)|$.)

\smallskip
Now, set $\Lambda=1/4m$.

\smallskip\noindent\textbf{Step 2. For $z$ near 0, most of the $H_z$ are  
close to $H_1$.} We shall establish the following: for each $i=1,\ldots, m$, the couple
of points $x=0$ and $y_i$ is  not a $(\Lambda,\varphi,r)$--good
couple.

Assume that the opposite were true and for some $i=1,\ldots, m$ we
had a $(\Lambda,\varphi,r)$--good couple $(x,y_i)$. Then, using  the two estimates
\begin{eqnarray}
\H^m(S(0,y_i;\varphi,r)) & \ge & \Lambda \omega(m)\varphi^{2m}r^m,
\label{m-S}\\
\H^m(\Sigma\cap B^n(y_i,\varphi^2 r)) & \ge & \frac{1}{2}
\omega(m)\varphi^{2m}r^m,\label{m-B}
\end{eqnarray}
 where \eqref{m-B} comes from Theorem~\ref{thm:UAR},
and the inequality of Lemma~\ref{1/R-est} to estimate $1/\rtp$, we
would obtain a lower bound for the energy,
\begin{eqnarray}
E & \ge  & \int_{S(0,y_i;\varphi,r)}\int_{\Sigma \cap
B^n(y_i,\varphi^2 r)} \frac{1}{\rtp{}^q(z,w)}\, d\H^m (w)\otimes d\H^m(z) \nonumber \\
& \ge & \frac{\Lambda}2 \omega(m)^2 \varphi^{4m} r^{2m}
\left(\frac 19 \, \frac \varphi r\right)^q \label{E2E} \\
& = & \frac{\omega(m)^2}{8m\cdot 9^q} \, \varphi^{4m+q} r^{2m-q} \
\ge \ 2 E\nonumber
\end{eqnarray}
by \eqref{c4c5} and the balance condition \eqref{balance}; this
contradiction proves that the claim of Step 2 does hold. In
particular, since the condition $ r/2< |y_i|< 2r$ is satisfied
for each $i$, we have
\begin{eqnarray}
\label{small-Syi} \H^m \Bigl(\bigcup_{i=1}^m
S(0,y_i;\varphi,r)\Bigr)
& \le  & \sum_{i=1}^m \H^m \bigl(S(0,y_i;\varphi,r)\bigr) \\
& <   & m\Lambda \omega (m) \varphi^{2m}r^m = \frac 14 \omega (m)
\varphi^{2m}r^m\, . \nonumber
\end{eqnarray}

\smallskip\noindent\textbf{Step 3. The new box contains $\Sigma \cap
B^n(x,2r)$.} We shall show that the cylinder $F$ defined by
\eqref{newbox} contains $\Sigma \cap B^n_{2r}$.

Again, we argue by contradiction. Suppose that there exists
$\zeta\in \Sigma\cap B^n_{2r}$ such that $\zeta\not \in F$. Set
\[
G:= \bigl( \Sigma^\ast \cap B^n(0,\varphi^2 r)\bigr) \setminus
\bigcup_{i=1}^m S(0,y_i;\varphi,r)\, .
\]
By  Theorem~\ref{thm:UAR} and \eqref{small-Syi}, we have
\begin{equation}\label{large-G}
\H^m(G) \ge \frac 14 \omega (m) \varphi^{2m}r^m\, ,
\end{equation}
and due to the definition of $S(0,y_i;\varphi,r)$ we know that
\begin{equation}
\label{short-Qi} |Q_{H_z}(y_i-z)| < \varphi r, \qquad z\in G, \quad
i=1,\ldots, m.
\end{equation}
Fix $z\in G$. \eqref{short-Qi} yields $|Q_{H_z}(y_i)| < \varphi r +
|z| \le 2 \varphi r$. Thus,  the basis $v_1,\ldots, v_m$ of
$W:=H_z$ given by
\[
v_i= y_i-Q_{H_z}(y_i),\qquad i=1,\ldots, m,
\]
satisfies{} $|v_i-y_i|\le 2\varphi r$ for each $i$. Letting
$w_i:=d^{-1} v_i$, we check that
\[
|w_i-h_i|=d^{-1}|v_i-y_i|\le \frac{2\varphi r}{d} < 3\varphi \ll
\frac{\eps_1}{2},
\]
as $6\varphi < (c_4)^{-1} \ll 10^{-1}(1+10^m)^{-1}=\eps_1$.
Invoking Lemma~\ref{bases2}  for $H=H_1$ and $W=H_z$, we conclude that
\[
\ang(H_1,W) \equiv \ang (H_1,H_z) \le 3c_3\varphi\, .
\]
Now, since $\zeta \not\in F$, we have $|Q_{H_1}(\zeta)| >
2c_4\varphi r$, and
\[
|Q_{H_1}(\zeta)-Q_W(\zeta)|\le \ang(H_1,W) \, |\zeta| \le
6c_3\varphi r\, .
\]
Thus, {} for $w\in B^n(\zeta,\varphi^2 r)$ and $z\in
G\subset B(0,\varphi^2 r)$
\begin{eqnarray*}
|Q_{H_z}(w-z)|\equiv |Q_W(w-z)| & = & |Q_W(\zeta-z)-Q_W(\zeta-w)|\\
& \ge  & |Q_W(\zeta)|-|z| - \varphi^2r \\
& \ge & |Q_{H_1}(\zeta)| - |Q_{H_1}(\zeta)-Q_W(\zeta)| - 2\varphi^2r\\
&  > & 2c_4\varphi r - 6 c_3\varphi r - 2\varphi^2 r \ge 5\varphi
r,
\end{eqnarray*}
since $c_4$ satisfies \eqref{c4c5} and $2\varphi^2\le\varphi
$.{} On the other hand, we certainly have $|w-z|\le 3r$ for
every point $w\in B^n(\zeta,\varphi^2 r)$. This yields
\[
\frac{1}{\rtp(z,w)}=\frac{2|Q_{H_z}(w-z)|}{|w-z|^2}  >
\frac{2\cdot 5\varphi r}{(3r)^2}
> \frac\varphi r, \qquad \mbox{for $z\in G$, $w\in
B^n(\zeta,\varphi^2 r)$.}
\]
We may now estimate the energy analogously to \eqref{E2E} and
obtain would obtain a lower bound for the energy,
\begin{eqnarray}
E & \ge    & \int_{G}\int_{\Sigma \cap
B^n(\zeta,\varphi^2 r)} \frac{1}{\rtp{}^q}(z,w)\, d\H^m(w)\, d\H^m(z)  \nonumber \\
&  > & \frac{1}{4\cdot 2} \omega(m)^2 \varphi^{4m} r^{2m}
\left( \frac \varphi r\right)^q \label{E2E-again} \\
& > & \frac{2}{c_5} \varphi^{4m+q} r^{2m-q} \ \ge \ 2 E.\nonumber
\end{eqnarray}
This is  again a contradiction, proving that $\Sigma\cap
B^n_{2r}\subset F$.

\smallskip\noindent\textbf{Step 4. The linking condition.} Since
we have established $\ang(H,H_1)\le  2  c_2\theta \le 2 
c_2\delta
\stackrel{\eqref{d+e}}{<} \eps_1$ in the first step of the proof,
the sphere
\[
M_1:=\Sphere(0, {\textstyle \frac 89} r; H_1^\perp)
\]
is contained in the interior of $B^n_r \setminus F$ and we have
$\dist(M_1,\Sigma) \ge \frac 89 r - 2 c_4\varphi r > \frac 59 r$,
since all points of $\Sigma \cap B^n_{2r}$ are in the cylinder $F$
defined in \eqref{newbox}, and $\varphi < 1/6c_4$. Thus, we may
deform $M_1$ homotopically to
\[
M_0:= \Sphere (0, {\textstyle \frac 89} r; H^\perp),
\]
so that the whole family of spheres realizing the homotopy stays
in $B^n_r\setminus F$, i.e. far away from $\Sigma$. (This can be
done precisely as in the verification of \eqref{link-M2M3} at the
end of the proof of Lemma~\ref{mainlemma}: we move the points of
$M_1$ to their projections onto $H^\perp$, and then deform the
resulting ellipsoid to obtain the round sphere $M_0$.)

Thus,
\[
\lk (M_1,\Sigma) = 1
\]
by Lemma~\ref{H-inv-link}. Now, every other sphere $\Sphere
(z,t;H_1^\perp)$, with $z\in H_1$, $|z|<(1-(c_4\varphi)^2)^{1/2}
\cdot 2r$ and $c_4\varphi \cdot 2r<t< (2r)^2-|z|^2$, i.e. every
$(n-m-1)$-sphere parallel to $H_1^\perp$ and contained in the
interior of $B^n_{2r}\setminus F$, can  obviously be deformed homotopically to $M_1$ without hitting points of $\Sigma$, since
$\Sigma\cap B^n\subset F$. Thus, again by Lemma~\ref{H-inv-link},
we conclude that Condition~(iii) of Definition~\ref{TB} is
satisfied for $F$ in $B^n_{2r}$.

This completes the whole proof of the lemma.\hfill $\Box$

\subsection{The tangent planes arise: an iterative construction}

\label{sec:5.3}

\label{iteration}

In this subsection, we apply Lemma~\ref{shrink-box} iteratively
and prove the following.

\begin{theorem}
\label{Tx-short}  Let $\delta(m)$ be the constant of 
Lemma~\ref{mainlemma}. Assume that $\Sigma\in \A(\delta)$ 
for some $\delta\in (0,\delta(m)]  $,  $\E_q(\Sigma) \le   E$,
$q > 2m$. Then $\Sigma$ is an embedded $m$-dimensional submanifold
of class $C^{1,\kappa}$, $\kappa=(q-2m)/(q+4m)$.
\end{theorem}

In fact, Theorem~\ref{Tx-short} will be just a corollary of
another result, which gives a lot of more precise, quantitative
information.

\begin{theorem}
 Let $\delta(m)$ be the constant of Lemma~\ref{mainlemma}. 
Assume that $\Sigma\in \A(\delta)$ for some
$\delta\in (0,\delta(m)]  $, $\E_q(\Sigma) \le   E$, $q>2m$. Then 
for each $x\in\Sigma$ there exists a unique plane
$T_x\Sigma\in G(n,m)$ (which we refer to as \emph{tangent plane of $\Sigma$
at $x$}) such that
\begin{equation}\label{distance_tp}
\dist(x',x+T_x\Sigma)\le C(n,m,q,E)|x'-x|^{1+\kappa} \quad\Foa
x'\in\Sigma, \quad x'\to x,
\end{equation}
Moreover, there
exists a  constant $a_2=a_2(n,m,q)>0$ with the following property.

Whenever $x,y\in \Sigma$ are such that
\begin{equation}
\label{xy-R2} 0< d:=|x-y|< R_2 := a_2 E^{-1/(q-2m)},
\end{equation}
then
\begin{equation}
\label{osc-tan} \ang (T_x\Sigma, T_y\Sigma) < c_6 E^{1/(q+4m)}
|x-y|^\kappa, \qquad \kappa = \frac{q-2m}{q+4m}
\end{equation}
for some constant $c_6$ depending only on $n,m$ and $q$. Moreover,
$U:=\Sigma\cap \INT B^n(x,R_2)$ is an open $m$-dimensional
topological disk, the orthogonal projection $\pi_{T_x\Sigma}$ onto $T_x\Sigma$
restricted to $U$ is injective, and each cylinder
\begin{equation}
\label{cyl-KN} K_N :=\{w\in B^n(x,2d_N)\colon \dist(w, x+
T_x\Sigma) \le \beta_N\cdot 2 d_N\}, \qquad N=1,2,\ldots
\end{equation}
with
\begin{equation}
\label{dNbetaN} d_N:=\frac{d}{5^{N-1}}, \qquad \beta_N = c_6
E^{{1}/(q+4m)} d_N^\kappa < \frac{1}{20}
\end{equation}
is a $(\beta_N, T_x\Sigma)$-trapping box for $\Sigma$ in
$B^n(x,2d_N)$.
\label{thm:5.6}
\end{theorem}

\medskip\noindent\textbf{Proof.}  A rough plan of the proof is the following. We shall first show, using Lemma~\ref{shrink-box} iteratively, that for each
$x\in \Sigma$ there exists a plane $H^\ast_x\in
G(n,m)$ such that for $x,y$ sufficiently close we have $\ang(H_x^\ast,H_y^\ast)\lesssim |x-y|^\kappa$. As a byproduct, we shall obtain a sequence of trapping boxes around each $H_x^*$, allowing us to show that $H_x^\ast$ is in fact unique. Finally, we set $T_x\Sigma=H_x^\ast$ and verify the statements concerning $\pi_{T_x\Sigma}$.

\smallskip\noindent\textbf{Step 1.}
Fix $x,y\in \Sigma$ and assume that \eqref{xy-R2} does
hold for a sufficiently small positive constant $a_2$ that shall
be specified later on. Fix $r_1>0$ such that
\begin{equation}
\label{r1} \frac 23 r_1< \frac{9}{10} (1-\delta^2)^{1/2} r_1 \le
|x-y|=d< r_1<R_{ 2}.
\end{equation}
Invoking Proposition~\ref{startbox} for $x$ and $r=r_1$, we obtain
a plane $H\in G(n,m)$ such that
\[
F:= \{w\in B^n(x,r_1)\colon \dist(w, x+H)\le \delta r_1\} \cup
B^n(x,r_1/2)
\]
is a $(\delta,H)$-trapping box for $\Sigma$ in $B^n(x,r_1)$.

Now, for $N=1,2,\ldots$ we set
\begin{eqnarray}
r_N & := & \frac{r_1}{5^{N-1}}, \label{rN}\\
\varphi_N & := & c_5^{1/(q+4m)} E^{1/(q+4m)} r_N^\kappa, \qquad
\kappa = \frac{q-2m}{q+4m},\label{varphiN}\\
\theta_N & :=& 10 c_4\varphi_N. \label{thetaN}
\end{eqnarray}
We have  $\varphi_N\lesssim r_N^\kappa \to 0$ as $N\to \infty$; the constant $a_2$ will
be chosen later,  in \eqref{def-cons-a2} below, so small that  $\delta$ and  $r_1$ shall
satisfy the assumptions of Lemma~\ref{shrink-box}. The choice of
$r_1$ guarantees that
\begin{equation}
\label{rNdNkappa} r_N^\kappa \, \le \,
\left(\frac{3}{2}\right)^\kappa d_N^\kappa \,  < \, \frac 32
d_N^\kappa\qquad\mbox{for all $N=1,2,\ldots$.}
\end{equation}

Apply Lemma~\ref{shrink-box} and Remark~\ref{2pt-plane} with  $\theta =\delta$, $r=r_1$ and $\varphi=\varphi_1$ to choose
$H_1\in G(n,m)$ such that $y-x$ in $H_1$ and the cylinder
\[
F_1:= \{w\in B^n(x,2r_1)\colon \dist(w, x+H_1)\le 2c_4\varphi_1
r_1\}
\]
is a $(c_4\varphi_1, H_1)$-trapping box for $\Sigma$ in
$B^n(x,2r_1)$. (The plane $H_1$ will serve, roughly speaking, as a
sort of average position for all tangent planes to $\Sigma$ in
$B^n(x,r_1)$.)

\smallskip\noindent\textbf{Step 2. The choice of $H^\ast_x$.}
Since $r_2=r_1/5$, we have $2c_4\varphi_1r_1=\theta_1r_2$, and the
intersection $F_1\cap B^n(x,r_2)$ provides a
$(\theta_1,H_1)$-trapping box for $\Sigma$ in $B^n(x,r_2)$.
Invoking Lemma~\ref{shrink-box} again, we find a plane $H_2\in
G(n,m)$ such that
\[
\ang(H_2,H_1) \le 2c_2\theta_1
\]
and the cylinder $F_2:= \{w\in B^n(x,2r_2)\colon \dist(w,
x+H_2)\le 2c_4\varphi_2 r_2\}$ is a $(c_4\varphi_2,
H_2)$-trapping box for $\Sigma$ in $B^n(x,2r_2)$. Proceeding
inductively, we find a sequence of planes $H_N\in G(n,m)$ such
that for each $N=1,2,\ldots$ the cylinder
\begin{equation}
\label{FN} F_N:= \{w\in B^n(x,2r_N)\colon \dist(w, x+H_N)\le
2c_4\varphi_N r_N\}
\end{equation}
is a $(c_4\varphi_N, H_N)$-trapping box for $\Sigma$ in
$B^n(x,2r_N)$
and we have the estimate
\begin{equation}
\label{angleHNN+1} \ang (H_{N+1},H_N) \le 2 c_2 \theta_N
\qquad\mbox{for all $N=1,2,\ldots$}
\end{equation}
Since $\sum \theta_N<\infty$, the planes $H_N$ converge to some
plane $H^\ast_x\in G(n,m)$ such that
\begin{eqnarray}
\ang (H_x^\ast,H_N) & \le & \sum_{j=N}^\infty  \ang (H_{j+1},H_j)
\nonumber  \\
& \le & 20c_2c_4\sum_{j=N}^\infty \varphi_j \qquad\mbox{by
\eqref{angleHNN+1} and \eqref{thetaN}} \nonumber\\
& = & 20c_2c_4c_5^{1/(q+4m)} E^{1/(q+4m)}\,  r_N^\kappa\,
\sum_{i=0}^\infty 5^{-i\kappa} \qquad\mbox{by
\eqref{rN} and \eqref{varphiN}}\nonumber\\
& \le & A r_N^\kappa, \qquad N=1,2,\ldots, \label{H*HN}
\end{eqnarray}
with
\begin{equation}
\label{A}
A:= \frac{40c_2c_4c_5^{1/(q+4m)} E^{1/(q+4m)}}{\kappa}\,
.
\end{equation}
For the last inequality above, we have used an elementary estimate
$5^\kappa/(5^\kappa-1) \le 2/\kappa$ which holds for each
$\kappa\in (0,1)$\footnote{Indeed, $f(\kappa)=\kappa a^\kappa\le 2(a^\kappa-1)=g(\kappa)$ for all $\kappa\in (0,1)$ and $a>e$, as  $f(0)=g(0)$ and $f'<g'$ on $(0,1)$.}.

Now, note that since $y-x\in H_1$ the initial cylinder $F_1$ is
such that $F_1\cap B^n(y,r_2)$ provides a
$(\theta_1,H_1)$-trapping box for $\Sigma$ in $B^n(y,r_2)$. Thus,
replacing the roles of $x$ and $y$ from the second step on, we may
run a similar iteration and obtain a plane $H^\ast_y$ such that
\begin{equation}
\label{HyH1} \ang (H_y^\ast, H_1) \le Ar_1^\kappa,
\end{equation}
together with a sequence of planes $P_N\to H_y^\ast$ (with
$P_1=H_1$) and appropriate trapping boxes determined by those
planes. By  the triangle inequality, \eqref{H*HN} for $N=1$ and
\eqref{HyH1} yield
\begin{equation}
\label{oscH*} \ang (H_x^\ast,H^\ast_y) \le 2 A r_1^\kappa\, .
\end{equation}
Once  the uniqueness of $H^\ast_x$ is established, we identify 
$H^\ast_x$ with $T_x\Sigma$ . The estimate \eqref{oscH*} combined with \eqref{A} will
yield the desired \eqref{osc-tan} (note that $r_1\approx|x-y|$ up
to a constant factor which is less than 2).

\smallskip\noindent\textbf{Step 3. Trapping boxes around
$H^\ast_x$.} It is now easy to check that tilting the cylinders
$F_N$ and enlarging them slightly, we can obtain new trapping
boxes $K_N$ for $\Sigma$ in $B^n(x,2r_N)$.

Fix $w\in F_N$. For sake of brevity, let $Q_\ast$ and $Q_N$ denote
the orthogonal projections of $\R^n$ onto $(H_x^\ast)^\perp$ and
$H_N^\perp$.  We have
\begin{eqnarray}
|Q_\ast(w-x)| & = & \bigl|Q_N(w-x) + (Q_\ast(w-x)-Q_N(w-x)) \bigr| \nonumber \\
& \le & 2c_4\varphi_Nr_N + \ang(H_x^\ast,H_N) |w-x|\label{Q*w} \\
& \le & 2A r_N^\kappa \cdot 2r_N,\nonumber
\end{eqnarray}
as $c_4\varphi_N \le A r_N^\kappa$. Hence, by \eqref{rNdNkappa},
\[
|Q_\ast(w-x)|  < 9A d_N^\kappa \cdot d_N\, .
\]
Therefore, if $\beta_N$ is defined by \eqref{dNbetaN} with
\begin{equation}
c_6 := 10A\,  E^{-1/(q+4m)} = 400 \kappa^{-1} c_2 c_4
c_5^{1/(q+4m)},
\label{def-c6}
\end{equation}
then we have $|Q_\ast(w-x )| < \beta_Nd_N$ for each $w\in F_N$.

Thus the cylinder 
\begin{equation}
\label{cyl-KNH*} K_N :=\{w\in B^n(x,2d_N)\colon \dist(w, x+
H_x^\ast) \le \beta_N\cdot 2 d_N\}, \qquad N=1,2,\ldots
\end{equation}
contains $F_N\cap
B^n(x,2d_N)$. It follows that $\Sigma\cap B^n(x,2d_N)\subset K_N$,
as $F_N$ was a trapping box for $\Sigma$ in a larger ball
$B^n(x,2r_N)$. It is easy to see that the linking condition of
Definition~\ref{TB} is also satisfied (we just take a smaller set
of spheres  that are slightly tilted) so that $K_N$ indeed is a
$(\beta_N,H_x^\ast)$-trapping box for $\Sigma$ in $B^n(x,2d_N)$.

Let us now specify $a_2$.  We choose this constant so that
\begin{equation}
\label{def-cons-a2}
c_6a_2^\kappa < \frac 1{20} \qquad\mbox{and}\qquad 0< a_2<a_1, 
\end{equation}
where $a_1$ is the constant of Theorem~\ref{thm:UAR}. Then,  by
\eqref{xy-R2},
\begin{eqnarray}
\beta_1 & = & c_6 E^{1/(q+4m)} |x-y|^\kappa \nonumber \\
& <   & c_6 E^{1/(q+4m)} R_2^\kappa \label{a2}\\
& = & c_6 E^{1/(q+4m)} a_2^\kappa E^{-1/(q+4m)} = c_6 a_2^\kappa <
\frac{1}{20}\, .\nonumber
\end{eqnarray}
For  these choices of $c_6$ and $a_2$ all applications of Lemma~\ref{shrink-box} were justified. Now, returning to \eqref{oscH*}, we obtain
\begin{eqnarray}
\ang(H_x^\ast,H_y^\ast) & \le & 2Ar_1^\kappa\nonumber \\
&\stackrel{\eqref{def-c6}}=& \frac{c_6}5 E^{1/(q+4m)} r_1^\kappa \nonumber \\
& \stackrel{\eqref{r1}}<   & \frac{c_6}5 E^{1/(q+4m)} \Bigl(\frac 32 |x-y|\Bigr)^\kappa \label{5.28b} \\
& <   &{c_6} E^{1/(q+4m)} |x-y|^\kappa\, .\nonumber
\end{eqnarray}
In particular, as $|  x-y|<R_2$, we also have
\begin{equation}
\label{5.28c}
\ang(H_x^\ast,H_y^\ast) <    {c_6} E^{1/(q+4m)} R_2^\kappa\, . 
\end{equation}

\smallskip

To finish the whole proof, it remains to demonstrate that
$H^\ast_x$ is indeed  unique and that $\Sigma\cap\INT
B^n(x,R_2)=U$ is an open $m$-dimensional disk such that the projection 
$\bigl.\pi_{H^*_x }\!\bigr|_{U}$ is injective.

\smallskip\noindent\textbf{Step 4. Uniqueness of $H^\ast_x$.}
Since formally Lemma~\ref{shrink-box} alone does not guarantee
that the choice of each new plane $H_N$ is unique, we must now
show that $H^\ast_x=\lim H_N$ \emph{is\/} unique.

Suppose that this were not the case, and that choosing
$L_N\not=H_N$ in some steps of the iteration we could obtain a
different limiting plane $L$, with $\ang(L, H_x^\ast)>0$.

Select $w\in H_x^\ast$ with $|w|=1$ such that $|w-\pi_{L}(w)|>
\vartheta > 0$. Set $V:= (H_x^\ast)^\perp$ and without loss of
generality suppose that $x=0$. The spheres
\[
M_N := \Sphere(d_N  w, 3\beta_Nd_N; V)
\]
are contained in $\INT B^n(0,2d_N)$,  away from $\Sigma$    since
$\beta_N\le\beta_1<1/20$, and by Lemma \ref{H-inv-link}  , 
are nontrivially linked with
$\Sigma$ since $K_N$ is a  $(\beta_N,H_x^\ast)$-trapping box for $\Sigma$ in
$B^n(0,2d_N)$. Since $L$ has been obtained by an analogous
iteration process, the cylinders
\[
\tilde{K}_N:=\{w\in B^n (0,2d_N)\, \colon \dist (w,L)\le
\beta_N\cdot 2d_N\}
\]
should also provide $(\beta_N,L)$-trapping boxes for $\Sigma$ in
$B^n(0,2d_N)$. However, taking $N$ so large that $6\beta_N <
\vartheta$, we obtain $\dist (d_N w, L) = d_N |w-\pi_{L}(w)| >
d_N\vartheta > 6\beta_Nd_N$. Thus, the sphere $M_N$ is contained
in the interior of $ B^n(x,2d_N )\setminus \tilde{K_N}$ and
satisfies the assumptions of Lemma~\ref{lk=0} with
$\eps=2\beta_Nd_N$ and therefore \emph{is not\/} linked with
$\Sigma$, a contradiction which proves that $H_x^\ast$ has to be
unique.

Moreover, since $K_N$ is a $(\beta_N,H_x^\ast)$-trapping box for
$\Sigma$ in $B^n(x,2d_N)$ and $\beta_N\approx d_N^\kappa)$ one
easily concludes that for $y\in \Sigma$ we have
\[
\dist (y, x+H_x^\ast) = O(|x-y|^{1+\kappa}) \qquad \mbox{as $y\to
x$,}
\]
 which justifies the definition $T_x\Sigma:=H_x^\ast$.

\smallskip\noindent\textbf{Step 5. Injectivity of the projection.} Again,
we argue by contradiction. Suppose that there exist $y\not=y_1\in
U\equiv \Sigma \cap \INT B^n(x,R_2)$ such that
$\pi_{T_x\Sigma}(y)=\pi_{T_x\Sigma}(y_1)$. Without loss of generality suppose that
\[
|x-y_1|\le |x-y|=d< R_2\, ;
\]
and let $d_N,\beta_N$ be defined by \eqref{dNbetaN}. Set $v=y_1-y$
and let  $Q_{T_x\Sigma}$, $Q_{T_y\Sigma}$ denote the projections onto 
$(T_x^\Sigma)^\perp=(H_x^\ast)^\perp$, $(T_y\Sigma)^\perp=(H_y^\ast)^\perp$, respectively. As $v\perp
T_x\Sigma$, we have $Q_{T_x\Sigma}(v)=v$ and
\begin{eqnarray*}
|Q_{T_y\Sigma}(v) | & = & |Q_{T_x\Sigma}(v) + (Q_{T_y\Sigma}(v)-Q_{T_x\Sigma}(v))| \\
& \ge & |v|\bigl(1- \|Q_{T_x\Sigma}-Q_{T_y\Sigma}\|\bigr) \\
& = & |v| \bigl(1 - \ang (H_x^\ast, H_y^\ast)\bigr) \\
& \ge & |v| (1- c_6 E^{1/{(q+4m)}} R_2^\kappa) \qquad\mbox{ by \eqref{5.28c}}\\
& = & |v| (1-c_6a_2^\kappa) > \frac{19}{20}|v|
\end{eqnarray*}
by  \eqref{def-cons-a2}. However, fixing $N$ so that $d_{N+1}<
|v|=|y-y_1|\le d_N$, we could use the trapping boxes constructed
along with $H_y^\ast$, i.e. the cylinders
\[
\{w\in B^n(y,2d_N)\colon \dist (w,y+H_y^\ast) \le \beta_N \cdot
2d_N\}
\]
which contain $\Sigma \cap B^n(y,2d_N)$, to estimate by virtue of \eqref{a2} 
\[
|Q_{T_y\Sigma}(v)|\le 2\beta_Nd_N = 10\beta_N d_{N+1} \le 10\beta_1 d_{N+1}
< \frac 12|v|,
\]
a contradiction.

 Since for each $d_1<R_2$
the cylinder
\begin{equation}
K_1 = \{w \in B^n(x,2d_1)\colon \dist(w,x+H_x^\ast)\le 2 \beta_1 d_1\}
\label{final-decay-b}
\end{equation}
is a $(\beta_1,H_x^\ast)$-trapping box for $\Sigma$ in
$B^n(x,2d_1)$, and $\beta_1<1/20$, we know by now -- as $d_1$ can be taken very close to $R_2$ -- that the image of $\pi_{T_x\Sigma}$ restricted to, say,
$\Sigma \cap B^n(x,3R_2/2)$, certainly contains the disk with
center at $\pi_{T_x\Sigma}(x)$ and radius $R_2$. It follows that
$U=\INT B^n(x,R_2)\cap \Sigma$ is a topological disk,  since $\pi_{T_x\Sigma}$ was also shown to be injective.\hfill $\Box$

\medskip\noindent\textbf{Proof of Theorem~\ref{thm:betas}.} As $\beta_1=c_6 E^{1/(q+4m)}d_1^\kappa$, it can be checked that Theorem~\ref{thm:betas} stated in the introduction  follows from \eqref{final-decay-b} and the definition of a trapping box. One can use the plane $H^*_x  \in G(n,m)$ to 
estimate the infimum in the definition \eqref{betasmall} of $\beta$-numbers. \hfill $\Box$

\subsection{Local graph representations of $\Sigma$}

\label{sec:5.4}

We shall now use Theorem~\ref{thm:5.6} to construct the graph
representations of an admissible surface $\Sigma$ with $\E_q(\Sigma)
<\infty$ for some $q>2m$. In Section~\ref{slicing}, this will be used to show that $\Sigma$ is in fact a manifold of class $C^{1,\mu}$ for $\mu= 1- 2m/q>\kappa$.

\begin{corollary} \label{graphpatch}
Suppose that $\Sigma\in \A(\delta)$ for some 
$\delta\in (0,\delta(m)]$, $\E_q(\Sigma) \le E$, $q>2m$. Let  $a_2>0$ and $R_2=a_2E^{-1/(q-2m)}$ denote the constants introduced in Theorem~\ref{thm:5.6}. Set $R_3=\frac 12 R_2$. Then, for each $x\in \Sigma$, the following is true.

There exists a  function 
\[
f\colon T_x\Sigma=:P\simeq  \R^m\to P^\perp\simeq\R^{n-m}
\]
of class $C^{1,\kappa}$, $\kappa=\frac{q-2m}{q+4m}$, such that
$f(0)=0$ and $\nabla f(0)=0$, and
\[
\Sigma\cap B^n(x,R_3) =   x+ \Bigl(\mathrm{graph}\, f\cap B^n(0,R_3)\,\Bigr)\, ,
\]     
where $\mathrm{graph}\, f \subset P\times P^\perp =\R^n$ denotes the graph of $f$, and 
\begin{equation}
|\nabla f(z)-\nabla f(w)|\le c_7  E^{1/({q+4m})}  |z-w|^\kappa\le c_7E^{1/({q+4m})}\, (2R_3)^\kappa,  \qquad z,w\in P \cap B^n(x,R_3), 
\label{1-graphosc}  
\end{equation}
for some constant $c_7$ depending only on $n,m, q$. 
\end{corollary}

\medskip\noindent\textbf{Proof.} Without loss of generality suppose that $x=0\in \R^n$ and $T_x\Sigma=P=\subsp(e_1,\ldots, e_m)$, where $e_j$, ${j=1,\ldots,n},$ form the standard orthonormal basis of $\R^n$. 
By Theorem~\ref{thm:5.6} we know that 
$$\bigl.\pi_P\bigr|_U \colon {\Sigma\cap B^n(x,R_2)}\to \pi(U)\subset P\, , 
\qquad U:=\Sigma\cap \INT B^n(x,R_2)\, , $$ is invertible. By \eqref{cyl-KN} and \eqref{dNbetaN} for $N=1$, the image  of this map contains an $m$-dimensional disk of radius $R_2'=R_2^2-(R_2/10)^2>\frac{9}{10}R_2$. 

\smallskip\noindent\textbf{Step 1.} We now let
\[
f\colon= \Bigl. Q_P\circ \Bigl(\bigl.\pi_P\bigr|_U\Bigr)^{-1} \Bigr|_D\colon D\to P^\perp, \qquad D=\INT\, D^m(0,R_2')\subset P,
\]
so that
\[
D\ni z\, \longmapsto\, F(z):=(z,f(z))\ \in\ P\times P^\perp =\R^n
\]
is a natural parametrization of $\Sigma$. Note that $F(D)$ contains $\Sigma\cap B^n(x,R_3)$ and that $f(0)=0$. Both $f$ and $F$ are continuous.

\smallskip\noindent\textbf{Step 2.} To prove that $\nabla f(0)$ exists and equals $0$,  use now the definition of $f$ to see that \eqref{xy-R2}, \eqref{cyl-KN} and \eqref{dNbetaN} of Theorem~\ref{thm:5.6} yield
\[
|f(z)|\ \stackrel{\eqref{cyl-KN}}\le \ 2\beta_Nd_N  \ \stackrel{\eqref{dNbetaN}}\le \ C(n,m,q,E)\, d_N^{1+\kappa} \qquad\mbox{for all  $N\in\N$} ,
\]
whenever $F(z)=(z,f(z))\in B^n(0,2d_N)$. (Recall that $d_N=d_1\cdot 5^{1-N}$; we are free to use any $d_1<R_2$ here.) Set $\rho_N:=d_N(1-\beta_N^2)^{1/2}$; by \eqref{dNbetaN}, $\frac{19}{20}d_N<\rho_N\le d_N$.  Thus, we also have $|f(z)|\le \mathrm{const}\cdot \rho_N^{1+\kappa}$ whenever $z\in D^m(0,2\rho_N)\subset P$. As $\rho_N\approx d_N=d_15^{1-N}$ for $N=1,2,\ldots $, this gives $|f(z)|=O(|z|^{1+\kappa})$ near $0$ and consequently $\nabla f(0)=0$.     

We shall now show that $F$ (and hence $f$) is differentiable at each $z\in D$. Fix $z\in D$ and $h\in P$ with $|h|$ small. Set  
\[
  L:= \Bigl(\Bigl.\pi_P\Bigr|_{T_{F(z)}\Sigma}\, \Bigr)^{-1}\colon P\to T_{F(z)}\Sigma\hookrightarrow \R^n\, .
\]
We have $F(z+h)-F(z)=L(h) + e,$ where the error $e=F(z+h)-F(z)-L(h)$ satisfies, by definition of $L$ and $F$, $\pi_P(e)=0$. Thus, $e=Q_P(e)$, so that     
\[
	|e|  \le  \left|\left( Q_P-Q_{T_{F(z)}\Sigma}\right) (e) \right| + \left|Q_{T_{F(z)}\Sigma}  (e)\right|  
	 \le  \frac{1}{20}|e| +  \left|Q_{T_{F(z)}\Sigma}  (e)\right|\qquad\mbox{by \eqref{osc-tan} and \eqref{dNbetaN} for $N=1$.}
\]
Absorbing the first term and using now Theorem~\ref{thm:5.6} at $x=F(z)$, we obtain
\begin{equation}
  \label{e-est}
|e|\le \frac {20}{19} \left|Q_{T_{F(z)}\Sigma}  (e)\right| =   \frac {20}{19} \mathrm{dist}\, \bigl(F(z+h), F(z) + T_{F(z)}\Sigma\bigr)  = O (|F(z+h)-F(z)|^{1+\kappa})\, .
\end{equation}
To finish the estimates, note that 
\[
|L(h)-h|= |\pi_{T_{F(z)}\Sigma} (L(h))- \pi_P(L(h)) |\ \stackrel{\eqref{osc-tan}}\le \ \frac{1}{20}|L(h)|\, ;
\]
therefore, $\frac{19}{20}|L(h)|\le |h|\le \frac{21}{20}|L(h)|$. Using this and  \eqref{e-est}, we now write
\[
|F(z+h)-F(z)|\le |L(h)|+|e| \le   \frac {20}{19} \left(|h|+   \left|Q_{T_{F(z)}\Sigma} (e)\right| \right) \le    \frac {20}{19} \left(|h|+ \mathrm{const}\cdot|F(z+h)-F(z)|^{1+\kappa} \right)\, .
\]                               
Now, for all $|h|$ sufficiently small we have $\frac {20}{19}\mathrm{const}\cdot|F(z+h)-F(z)|^{1+\kappa}  < \frac{1}{2} |F(z+h)-F(z)|$, as $F$ is continuous at $z$. Thus, the second term can be absorbed, yielding $|F(z+h)-F(h)|=O(|h|)$ as $h\to 0$. Plugging this into the right hand side of \eqref{e-est}, we obtain the desired error estimate $|e|=O(|h|^{1+\kappa})=o(|h|)$ as $h\to 0$. Therefore, $F$ is differentiable at $z$ with $DF(z)=L$. 

\medskip 

The uniform H\"{o}lder bound for $\nabla f$ results now from one more application of the oscillation estimate \eqref{osc-tan} for tangent planes:

\smallskip\noindent\textbf{Step 3.}  With
$$
|\partial_if(w)-\partial_if(z)|=\left|\Big[\begin{array}{c}
e_i\\
\partial_if(w)
\end{array}\Big]-
\Big[\begin{array}{c}
e_i\\
\partial_if(z)
\end{array}\Big]\right|=
\left|\pi_{T_{F(w)}\Sigma}\left(\Big[\begin{array}{c}
e_i\\
\partial_if(w)
\end{array}\Big]\right)-
\pi_{T_{F(z)}\Sigma}\left(\Big[\begin{array}{c}
e_i\\
\partial_if(z)
\end{array}\Big]\right)\right|
$$
we can estimate
\begin{eqnarray*}
|\partial_if(w)-\partial_if(z)| & \le &
\left|\pi_{T_{F(w)}\Sigma}\left(\Big[\begin{array}{c}
e_i\\
\partial_if(w)
\end{array}\Big]\right)-
\pi_{T_{F(z)}\Sigma}\left(\Big[\begin{array}{c}
e_i\\
\partial_if(w)
\end{array}\Big]\right)\right|\\
&& \qquad +
\left|\pi_{T_{F(z)}\Sigma}\left(\Big[\begin{array}{c}
e_i\\
\partial_if(w)
\end{array}\Big]\right)-
\pi_{T_{F(z)}\Sigma}\left(\Big[\begin{array}{c}
e_i\\
\partial_if(z)
\end{array}\Big]\right)\right|\\
&& \hspace{-2cm} \le 
\ang (T_{F(w)}\Sigma,T_{F(z)}\Sigma)(1+|\nabla f(w)|^2)^{1/2}+
\left|(\pi_{T_{F(z)}\Sigma}-\pi_{T_0\Sigma})\left(\Big[\begin{array}{c}
0\\
\partial_if(w)-\partial_if(z)
\end{array}\Big]\right)\right|\\
& &\hspace{-2.2cm}\overset{\eqref{osc-tan}}{<}  c_6E^{1/(q+4m)}|w-z|^\kappa + \ang(T_{F(z)}\Sigma,T_0\Sigma)|\partial_if(w)-\partial_if(z)|.
\end{eqnarray*}
Since  $\ang(T_{F(z)}\Sigma,T_0\Sigma)<1/2$ by \eqref{osc-tan} and our choice
of constants, we can absorb the right term on the left-hand side to conclude.

%
Now, using a standard cutoff technique, we leave $f$ unchanged on $D^m(0,2R_2/3)$, and extend it to the whole plane $P$, so that the extension vanishes off $D^m(0,3R_2/4)$. The corollary follows.\hfill $\Box$

\section{Slicing and bootstrap to optimal H\"older exponent }

\label{sec:6}

\setnumbers

\label{slicing}

In this  section we assume that $\Sigma$ is a flat $m$-dimensional
graph of class $C^{1,\kappa}$ having finite tangent-point energy
$\E_q(\Sigma)$. The goal is to show how to bootstrap the
H\"{o}lder exponent  $\kappa$ to $\mu = 1- 2m/q$.

Relying on Corollary~\ref{graphpatch}, without loss of generality we can assume that
\[
\Sigma \cap B^n(0,5R) = \mathrm{Graph}\, f \cap B^n(0,5R)
\]
for a fixed number $R>0$, where
\[
f\colon P\cong \R^m \to P^\perp \cong \R^{n-m}
\]
is of class $C^{1,\kappa}$ and satisfies $\nabla f(0)=0$,  $f(0)=0$,
\begin{equation}
\label{flatgraph} |\nabla f|\le \eps_0:=\frac{\eps_1}{800m\, c_2}=2^{-5}10^{-3}m^{-2}(10^m+1)^{-2} \qquad\mbox{on $P$. }
\end{equation}
To achieve \eqref{flatgraph}, we use \eqref{1-graphosc} of Corollary~\ref{graphpatch} and shrink $R_3$ by a constant factor if necessary.
The number $\eps_0$ is chosen so that $\eps_0< \eps_1 /(400m c_2)$
for the constants $\eps_1$ and $c_2$ used in Lemma~\ref{bases} and
other auxiliary estimates in Section~\ref{sec:linear}. We let
$F\colon P\to \R^n$ be the natural parametr\-ization of
$\Sigma\cap B(0,5R)$, given by $F(x) = (x,f(x))$ for $x\in P$;
outside $B^n(0,5R)$ the image of $F$ does not have to coincide
with $\Sigma$. The choice of $\eps_0$ guarantees that, due to
Lemma~\ref{bases}~(ii),
\begin{equation}
\ang (T_{F(x)}\Sigma, T_{F(0)}\Sigma ) \le c_2\eps_0 <
\frac{\eps_1}{400m} \label{angle-0}
\end{equation}
whenever $x\in B^n(0,5R)\cap P$. Thus,
\begin{equation}
\ang (T_{F(x_1)}\Sigma, T_{F(x_2)}\Sigma ) < \frac{\eps_1}{200m}
<\frac{1}{m 4^{m+1}} \qquad \mbox{for all $x_1,x_2\in
B^n(0,5R)\cap P$.} \label{tan-angles}
\end{equation}

 As in our paper
\cite[Section~6]{svdm-surfaces}, we introduce the maximal
functions controlling the oscillation of $\nabla f$ at various
places and scales,
\begin{equation}
\label{Phi-ast} \Phi^\ast(\varrho, A) = \sup_{{B_\varrho\subset A}}
\left(\osc_{B_\varrho} \nabla f\right)
\end{equation}
where the supremum is taken over all  possible closed $m$-dimensional balls
$B_\varrho$ of radius $\varrho$ that are contained in a subset  $A
\subset B^n(0,5R) \cap P$, with $\varrho\le 5R$. Since $f\in
C^{1,\kappa}$, we have a priori
\begin{equation}
\label{apriori} \Phi^\ast(\varrho, A)\le C\varrho^\kappa
\end{equation}
for some constant $C$ which does not depend on $\varrho,A$.

To show that $f\in C^{1,\mu}$ for $\mu=1-2m/q$, we check
that locally,  on each scale, the oscillation of $\nabla f$ is
controlled by a main term  which involves the local energy
and resembles the right hand side of \eqref{optimal}, up to
a small error,  which itself is controlled by the oscillation of $\nabla f$  on a
much smaller scale.

\begin{lemma}\label{key}
Let $f$, $F$, $\Sigma$, $R>0$ and $P$ be as above. If $z_1,z_2\in
B^n(0,2R)\cap P$ with $|z_1-z_2|=t>0$, then for any $N>2$ we have
\begin{equation}
\label{osc-error} |\nabla f(z_1)-\nabla f(z_2)| \le 2 \Phi^\ast
(t/N,B) + C(N, m,q) \,  E_B^{1/q} \, t^\mu
\end{equation}
where $B:=B^m(\frac{z_1+z_2}2,t)$ is an $m$-dimensional disc in
$P$, $\mu:=1-2m/q$, and
\begin{equation}
\label{EonB} E_B=\int\!\! \int_{F(\! B)\times F(\! B)} \rtp^{-q}
\,\, d\H^m\otimes d\H^m
\end{equation}
is the local energy of $\Sigma$ over $B$.
\end{lemma}

\noindent\textbf{Remark.} Once this lemma is proved, one can fix
an $m$-dimensional disk $B^m(b,s)\subset B^n(0,R)\cap P$ and use
\eqref{osc-error} to obtain for $ t\le s$
\begin{equation}
\label{pre-morrey} \Phi^\ast(t, B^m(b,s))\le 2 \Phi^\ast
\bigl(2t/N, B^m(b,s+2t) \bigr) + C( N,m,q)\, M_q(b,s+2t)\,
t^\mu\,
\end{equation}
where
\[
M_q(b,r):= \left(\int\!\! \int_{F(\! B(b,r))\times F(\! B(b,r))}
\rtp^{-q} \,\, d\H^m\otimes d\H^m\right)^{1/q}\, .
\]
Fixing $N>2$ such that  $2^\kappa/N^\kappa<\frac 12$ we obtain  $2^j\cdot (2/N)^{j\kappa}\to
0$ as $j\to \infty$. Using this, one can iterate
\eqref{pre-morrey} and show that
\[
\osc_{B^m(b,s)} \nabla f \le C'( m,q) M_q(b, 5s) s^\mu\, .
\]
Combining this estimate with Corollary~\ref{graphpatch}, 
we obtain Theorem~\ref{thm:bootstrap} stated in the introduction. 
Note that in fact the result holds for all surfaces $\Sigma\in \A(\delta)$ 
for $\delta\in (0,\delta(m)]  $, where $\delta(m)$ is the constant of Lemma~\ref{mainlemma}.

The remaining part of this  section is devoted to the

\medskip\noindent\textbf{Proof of Lemma~\ref{key}.} Fix $z_1,z_2$
and the disk $B$ as in the statement of the lemma; we have
$\H^m(B)=\omega(m)t^m$. Pick $N>2$ and let  $E_B$ be the local energy
of $\Sigma$ over $B$, defined by \eqref{EonB}. Assume that $\nabla
f\not\equiv \mathrm{const}$ on $B$, for otherwise there is nothing
to prove.

\smallskip\noindent\textbf{Step 1.} Take
\begin{equation}
\label{K0} K_0:= \left( E_B\cdot N^{2m} \omega(m)^{-2}\right)^{1/q} >
0
\end{equation}
and set
\begin{eqnarray}
Y_1 & :=& \{ x_1\in B \ \colon \H^1(Y_2(x_1))
           \ge N^{-m}\H^m(B) \}\, ,\label{bad1}\\
Y_2(x_1) & := & \Bigl\{x_2\in B \ \colon \frac{1}{\rtp(F(x_1),F(x_2))}> K_0 \, t^{-2m/q}
\Bigr\}\, .\label{bad2}
\end{eqnarray}

We now estimate the local energy to obtain a bound for
$\H^m(Y_1)$, shrinking the domain of integration, as follows:
\begin{eqnarray*}
E_B & =& \int\!\! \int_{F(\! B)\times F(\! B)} \rtp^{-q} \,\,
d\H^m\otimes d\H^m \\
&\ge & \int\!\! \int_{B\times B} \biggl(\frac{1}{\rtp(F(x_1),F(x_2))}\biggr)^q\, dx_1\, dx_2  \\
&\ge & \int_{Y_1} \biggl(\int_{Y_2(x_1)} \biggl(\frac{1}{\rtp(F(x_1),F(x_2))}\biggr)^q\, dx_2\,\biggr) dx_1 \\
&\stackrel{\eqref{bad1},\ \eqref{bad2}}{>} & \H^m(Y_1)
N^{-m}\H^m(B) \, K_0^q t^{-2m} \ = \ E_B \H^m(Y_1) N^m
\bigl(\H^m(B)\bigr)^{-1}\, .
\end{eqnarray*}
The last equality follows from \eqref{K0}. Thus, we obtain
\[
\H^m(Y_1) < \frac{1}{N^m} \H^m(B),
\]
and since the radius of $B$ equals $t$, we obtain
\begin{equation}
\label{close}
B^m(a_i,t/N) \setminus Y_1 \not=\emptyset
\qquad\mbox{for $i=1,2$.}
\end{equation}
Now, select two points $u_i \in B^m(a_i,t/N) \setminus Y_1$
($i=1,2$). By  the triangle inequality,
\begin{eqnarray*}
|\nabla f(z_1)-\nabla f(z_2)| & \le & |\nabla f(z_1)-\nabla f
(u_1)| + |\nabla f(u_2)-\nabla f(z_2)| + |\nabla f(u_1)-\nabla
f(u_2)| \\
& \le & 2 \Phi^\ast(t/N,B) + |\nabla f(u_1)-\nabla f(u_2)|\, .
\end{eqnarray*}
Thus, it remains to show that the last term, $|\nabla
f(u_1)-\nabla f(u_2)|$, does not exceed a constant multiple of
$E_B^{1/q}\, t^\mu$. To achieve this goal, we assume that
$\nabla f(u_1)\not= \nabla f(u_2)$ and work with the portion of
the surface parametrized by the points in
\begin{equation}
\label{G} G:= B\setminus \bigl(Y_2(u_1) \cup Y_2(u_2)\bigr)\, .
\end{equation}
By \eqref{bad1}, $G$ satisfies
\begin{equation}
\label{G-below} \H^m(G) > (1-2N^{-m}) \H^m(B) =: C_1( q,m)\,
t^m\, .
\end{equation}
To conclude the whole proof, we shall derive an upper estimate for
the measure of $G$,
\begin{equation}
\label{G-above} \H^m(G) \le  C_2(q,m)\, K_0\,
\frac{t^{m+\mu}}{\alpha}, \,
\end{equation}
where $\alpha:=\ang (H_1,H_2)\not= 0$ and $H_i:=T_{F(u_i)}\Sigma$
denotes the tangent plane to $\Sigma$ at $F(u_i)\in \Sigma$ for $i=1,2.$ 
Combining \eqref{G-above} and \eqref{G-below}, we will then obtain
\[
\alpha <   (C_1)^{-1} C_2 K_0 t^\mu =:  C_3 E_B^{1/q} t^\mu\,
.
\]
(By a reasoning analogous to the proof of Corollary~\ref{graphpatch}, this also yields an estimate for the oscillation of $\nabla f$.)

\smallskip\noindent\textbf{Step 2. Proof of \eqref{G-above}.} By
\eqref{tan-angles}, we have  $\alpha=\ang(H_1,H_2) <
m^{-1}4^{-m-1}$. By Lemma~\ref{proj-meas}  applied to $\eps=m^{-1}4^{-m-1}$, we obtain
\[
\H^m(G)\le \H^m(F(G))  < 2 \H^m\bigl(\pi_{H_1}(F(G))\bigr),
\]
so that \eqref{G-above} would follow from
\begin{equation}
\H^m\bigl(\pi_{H_1}(F(G))\bigr)\le C_4\, K_0\,
\frac{t^{m+\mu}}{\alpha}\, . \label{proj-G-H1}
\end{equation}
Now, for $\zeta \in G$ and $i=1,2$ we have by \eqref{bad2}
\[
\frac{1}{\rtp(F(u_i),F(\zeta))} =
\frac{2\bigl|Q_{H_i}(F(\zeta)-F(u_i))\bigr|}{|F(\zeta)-F(u_i)|^2}\le
K_0 t^{-1+\mu}\, .
\]
Let $P_i=F(u_i)+H_i$ be the affine tangent plane to $\Sigma$ at
$F(u_i)$. Since $F$ is Lipschitz with constant $(1+\eps_0)$ and
$|z-u_i|\le  2 t$,
\begin{eqnarray}
\dist (F(\zeta),P_i) & = & \dist (F(\zeta)-F(u_i), H_i) \label{z-Pi} \\
& = & \bigl|Q_{H_i}(F(\zeta)-F(u_i))\bigr|\  \le \  8K_0
t^{1+\mu} =: h_0 \nonumber
\end{eqnarray}
for $\zeta\in G$, $i=1,2$.  Select the points $p_i\in P_i$, $i=1,2$,
so that $|p_1-p_2|=\dist (P_1,P_2)$. The vector $p_2-p_1$ is then
orthogonal to $H_1$ and to $H_2$, and since $G$ is nonempty by
\eqref{G-below}, we have $|p_1-p_2|\le 2h_0$ by \eqref{z-Pi}.

Set $p=(p_1+p_2)/2$, pick a parameter $\zeta\in G$ and consider
$y=F(\zeta)-p$. We have
\[
y= (F(\zeta) - F(u_1)) + (F(u_1)-p_1) + (p_1-p),
\]
so that $\pi_{H_1}(y) = \pi_{H_1} (F(\zeta)-F(u_1)) + (F(u_1) - p_1)
$, and
\begin{eqnarray*}
|y-\pi_{H_1}(y)| & = & |(p_1-p) + F(\zeta) - F(u_1) -  \pi_{H_1}
(F(\zeta)-F(u_1))|\\
& = & |(p_1-p) + Q_{H_1}(F(\zeta) - F(u_1)) |
\, .
\end{eqnarray*}
Therefore, since $|p-p_1|\le h_0$ and by  \eqref{z-Pi}, $
|y-\pi_{H_1}(y)| \le h_0 + h_0 = 2h_0$. In the same way, we obtain
$|y-\pi_{H_2}(y)| \le 2h_0$. Thus,
\[
\frac{y}{2h_0}=\frac{F(\zeta)-p}{2h_0} \in S(H_1,H_2),
\]
where $S(H_1,H_2)=\{x\in \R^n \colon \dist(x,H_i)\le 1\mbox{ for
$i=1,2$}\}$ is the intersection of two slabs considered in
Section~\ref{sec:linear}. Applying Lemma~\ref{proj-slab} which is possible due to the estimate \eqref{tan-angles} for $\ang(H_1,H_2)$, we
conclude that there exists an $(m-1)$-dimensional subspace
$W\subset H_1$ such that
\begin{equation}
\label{strip} \pi_{H_1}(F(G)-p) \subset \{x\in H_1\, \colon\,
\dist (x,W)\le 2h_0\cdot 5c_2/\alpha\}\, .
\end{equation}
On the other hand, since $F$ is Lipschitz, we certainly have
$F(G)\subset  B^n(F(\frac{a_1+a_2}2),2t)$ and therefore
\begin{equation}
\label{ball} \pi_{H_1}(F(G)-p)\subset B^n(a,2t), \qquad
a:=\pi_{H_1}(F( \frac{a_1+a_2}2)-p).
\end{equation}
Combining \eqref{strip} and \eqref{ball}, we invoke
Lemma~\ref{strip-ball} to $H:=H_1$, $S':=\pi_{H_1}(F(G)-p)$, and
$d:=2h_05c_2/\alpha$, to obtain
\[
\H^m\bigl(\pi_{H_1}(F(G))\bigr) \le 4^{m-1}t^{m-1} \cdot 20h_0
c_2/\alpha =: C_2 (m) K_0\frac{t^{m+\mu}}{\alpha}\, ,
\]
which is \eqref{proj-G-H1}, implying \eqref{G-above} and thus
completing the proof.


\small
\vspace{1cm}
\begin{minipage}{56mm}
{\sc Pawe\l{} Strzelecki}\\
Instytut Matematyki\\
Uniwersytet Warszawski\\
ul. Banacha 2\\
PL-02-097 Warsaw \\
POLAND\\
E-mail: {\tt pawelst@mimuw.edu.pl}
\end{minipage}
\hfill
\begin{minipage}{56mm}
{\sc Heiko von der Mosel}\\
Institut f\"ur Mathematik\\
RWTH Aachen\\
Templergraben 55\\
D-52062 Aachen\\
GERMANY\\
Email: {\tt heiko@}\\{\tt instmath.rwth-aachen.de}
\end{minipage}

\end{document}